\documentclass[11pt]{article}
\usepackage{graphicx}
\RequirePackage[numbers,sort&compress]{natbib}
\RequirePackage[colorlinks,citecolor=blue,urlcolor=blue]{hyperref}
\usepackage{geometry}
\usepackage{setspace}
\usepackage{rotating}
\usepackage[T1]{fontenc}
\usepackage[centertags]{amsmath}
\usepackage{amsfonts}
\usepackage{lscape}
\usepackage{amssymb}
\usepackage{amsmath}
\UseRawInputEncoding
\usepackage{epstopdf}
\usepackage{color}

\numberwithin{equation}{section}
\newtheorem{theorem}{Theorem}[section]

\newtheorem{remark}{Remark}[section]
\newtheorem{corollary}{Corollary}[section]

\newtheorem{lemma}{Lemma}[section]

\newtheorem{assumption}{Assumption}[section]

\begin{document}
\title{Composite $L^p$-quantile regression, near quantile regression and the oracle model selection theory}
\author{
\thanks{\hspace*{\parindent}
Corresponding author, E-mail address: linfuming@suse.edu.cn(F. Lin). } Fuming Lin \ Weilin Mou \\
{\small School of Mathematics and Statistics, Sichuan University
of Science \&
Engineering,}\\ {\small Sichuan Zigong, China} \\
}
\date{}
\maketitle
\begin{quote}
{\bf Abstract:} High-dimensional quantile regression and asymmetric
least squares regression have wide applications in statistics,
econometrics, finance, etc. But these two types of methods are apt to incur some shortcomings, such as
quantile regression having low efficiency in many cases, linear program and interior point algorithms
it heavily computationally depends being almost frozen in a typical desktop computer; the second one
theoretically requiring much higher moments.
In this paper, we consider high-dimensional $L^p$-quantile regression which only requires
a finite $2(p-1)$th ($1<p\leq 2$) moment of the error and is also a natural generalization of the above methods
and $L^p$-regression as well. The loss function of $L^p$-quantile regression circumvents
the non-differentiability of the absolute loss function and the difficulty of
the squares loss function requiring the finiteness of error's variance
and thus promises excellent properties of $L^p$-quantile regression.
Specifically, we first develop a new method called composite $L^p$-quantile
regression(CLpQR). We study the oracle model selection theory
based on CLpQR (call the estimator CLpQR-oracle) and show in some cases of $p$ ($p>1$) CLpQR-oracle behaves
better than CQR-oracle (based on composite quantile regression) when error's variance is infinite.
Moreover, CLpQR has high efficiency and can be sometimes arbitrarily more efficient than
both CQR and the least squares regression. Second, in order to deal with quantiles issues using CLpQR
we further propose another new regression method, i.e. near quantile regression and prove the asymptotic
normality of the estimator when $p\rightarrow1+$ and the sample size
$T\rightarrow\infty$ simultaneously in any way. As by-products, a new thought of
smoothing quantile objective functions and a new estimation are provided
for the asymptotic covariance matrix of quantile regression.
Third, we develop a unified efficient algorithm for fitting high-dimensional
$L^p$-quantile regression ($p\geq1$) by combining the cyclic coordinate
descent and an augmented proximal gradient algorithm.
Remarkably, the algorithm turns out to be a favourable alternative of the
commonly used liner programming and interior point algorithm when fitting
quantile regression.

{\bf MSC2020 subject classifications.} Primary 62J07; Secondary 62G08.

\end{quote}
\section{Introduction}
As ever widely used variable and model selecting methods, the AIC and BIC
criteria are overwhelmed with the ever increasing high-dimensional
and even ultra-high dimensional data. In this regard, various methods have been proposed for
analyzing high-dimensional data, among which sparse estimation is a dominant approach due to
its selecting variables and estimating coefficients simultaneously.
Using $L_{1}$-penalized least squares loss, Tibshirani (1996)\cite{TIBSHIRANI} developed
the least absolute shrinkage and selection operator, i.e. Lasso. Fan and Li (2001)\cite{Fan} observed
Lasso's $L_{1}$-penalty yielding biased estimates and suggested instead using the SCAD penalty
which yields a unbiased coefficient estimation and its desired oracle properties as well.
Zou (2006)\cite{ZOU} introduced the adaptive Lasso which also enjoys the oracle properties.

Although sparse regression has convenient theoretical derivation and efficient algorithms based on the
squared loss function, it easily suffers some drawbacks such as the breakdown issue when the error
variance is infinite and over-sensitivity to outliers. Hence, Zou and Yuan (2008)\cite{ZOUY}  and Wu and Liu (2009)\cite{Wu}  had
recourse to quantile regression first introduced by Koenker and Bassett (1978)\cite{Koenker Gilbert} and developed penalized
composite quantile regression and penalized quantile regression, respectively. Zou and Yuan (2008)\cite{ZOUY} also
considered the oracle estimator, namely CQR-oracle that estimates the coefficient vector by their method
and referred to as LS-oracle the corresponding estimator by the least squares. Thanks to its asymmetric absolute loss function the quantile regression theory has no moment assumptions on the error and quantile regression allows modeling of the entire conditional distribution of the response variable $y$ given covariate $X$ which can show heterogeneity in the relationship between $X$ and $y$. However, quantile regression may have non-unique solutions (Koenker and Bassett (1978)\cite{Koenker Gilbert}), is inefficiency for Gaussian-like errors and has estimation difficulty with the asymptotic covariance matrix. There is a greater concern in its computational aspect.
It is known that linear program and interior point algorithms are usually used to solve regression quantile
optimization problems. But these two algorithms tend to be slow or too memory-intensive in deal with high-dimensional data on a ordinary computer and thus quantile regression may lack attraction compared to other machine learning tools (Gu and Zou (2016)\cite{Gu}, He et al. (2023)\cite{He}). Based on asymmetric least squares regression (proposed by Newey and Powell (1987)\cite{Newey} and also called expectile regression) in stead of quantile regression, Gu and Zou (2016)\cite{Gu} developed two methods: sparse asymmetric least squares regression and coupled sparse asymmetric least squares regression
to consider heteroscedasticity detection in the high-dimensional data. Embracing the sparse least squares regression as a special case the sparse asymmetric least squares regression still gets stuck with the difficulties of the latter as mentioned at the beginning of this paragraph.

So recently, originating from Efron (1991)\cite{Efron} and Chen (1996)\cite{Chen}, $L^{p}$-quantile regression has received ever growing attention due to it relieving insufficiencies of both of quantile and expectile regression.
Hu et al. (2021)\cite{Hu} considered high-dimensional $L^{p}$-quantile regression and investigated its oracle properties. $L^{p}$-quantile regression usually sets $p\geq1$. When letting $p$ take 1 or 2 in the $L^{p}$-quantile regression loss function, quantile or expectile regression is restored and when the weight being 0.5, $L^{p}$ regression appears.

In this paper, we systematically study some problems about $L^p$-quantile
regression, which only requires a finite $2(p-1)$th ($1<p\leq 2$) moment of the error and thus can be used to analyze heavy-tailed data. We prove the asymptotic theory for the composite $L^p$-quantile
regression (CLpQR for short) under mild conditions. For dealing with high-dimensional data, we define the CLpQR oracle estimator (CLpQR-oracle), analyse its asymptotic relative efficiency in detail, and develop the oracle model selection theory.
Since quantile regression will confront the computational difficulty mentioned above we attempt to use CLpOR to get over it.
To this end we develop a new regression method, i.e. near quantile regression and prove the asymptotic normality of the estimator when $p\rightarrow1+$ and the sample size $T\rightarrow\infty$ simultaneously in any way. The near quantile regression has many other important applications. Here are two application scenarios that come to mind immediately. One of them is that we can obtain a new estimation for the asymptotic covariance matrix of quantile regression without involving the estimation of the density function of the error as current methods do. The other one is concerned with an intriguing issue all the time, i.e. smoothing the objective function of quantile regression. While current methods mainly apply smooth kernel functions to modify the objective function of quantile regression, see Horowitz (1998)\cite{Horowitz}, Fernandes et al. (2021)\cite{Fernandes}, He et al. (2023)\cite{He} among others, near quantile regression acts as a natural choice as its objective function itself is differentiable and facilitate gradient-based optimization methods. Finally, we develop a unified efficient algorithm for fitting high-dimensional $L^p$-quantile regression ($p>1$) by combining the cyclic coordinate
descent and an augmented proximal gradient algorithm. Surprisingly, the algorithm can also fit high-dimensional quantile regression very well in our random simulation and empirical analysis. The study on asymptotic relative efficiency illustrates that CLpQR-oracle has high efficiency and can be sometimes arbitrarily more efficient than both CQR-oracle and LS-oracle. Simulation results show that
the proposed algorithm can be efficiently applied to deal with the quantile issues (this manifests CCPA is a useful alternative of linear program and interior point algorithms, especially in high-dimensional cases),
in some cases of $p$ ($p>1$) CLpQR-oracle behaves better than CQR-oracle in terms of estimation accuracy even when the error variance is infinite, and the asymptotic theorems for near quantile regression works very well in several distribution cases. In the empirical analysis, we provide a method for choosing the suitable values of $p$.

The paper proceeds as follows. Section 2 is devoted to the definition of the CLpQR estimator, its asymptotic normality and asymptotic relative efficiency. Section 3 contains the CLpOR-oracular estimation theory. Near quantile regression is expounded in Section 4. In Section 5, we describe the algorithm for fitting CLpQR, CLpQR-oracle and quantile regression. Simulation and empirical analysis are contained in Section 6 and Section 7. All proofs and lemmas are presented in Section 8. We conclude the paper with Section 9.
\section{CLpQR and asymptotic relative efficiency}
\subsection{Estimator's definition and its asymptotic normality}
Suppose the data come from the following linear model
\begin{equation}\label{eq1.1}
y=\textbf{x}'\boldsymbol{\beta}^{*}+\varepsilon,
\end{equation}
where $\textbf{x}$ is the centered predictor, $\boldsymbol{\beta}^{*}$ the unknown $m$-dimensional parameter vector and
$\varepsilon$ the error term. Consider the loss function associated with $L^{p}$-quantiles ($p>1$) as follows:
\begin{equation}\label{eq1.2}
\boldsymbol{\eta}_{\tau,p}(s)=|\tau-I(s<0)||s|^{p},
\end{equation}where $\tau$ is called weight.
According to its linear transformation invariance the $\tau$th $L^{p}$-quantile of $y$ can be written as
\begin{equation}\label{eq1.3}
\textbf{x}'\boldsymbol{\beta}^{*}+b^{*}_{\tau},
\end{equation}
where $b^{*}_{\tau}$ is the $\tau$th $L^{p}$-quantile of $\varepsilon$.

Setting various weights such that $0<\tau_{1}<\tau_{2}<\cdots<\tau_{K}<1$, define $\hat{\boldsymbol{\beta}}^{clp}$ as the composite $L^{p}$-quantile regression estimator of $\boldsymbol{\beta}^{*}$ calculated by
\begin{eqnarray}\label{eq1.4}
(\hat{b}_{1}, \cdots, \hat{b}_{K}, \hat{\boldsymbol{\beta}}^{clp})
=\arg\min_{b_{1}, \cdots, b_{K}, \boldsymbol{\beta}}\sum^{K}_{k=1}\sum^{T}_{t=1}\boldsymbol{\eta}_{\tau_{k},p}(y_{t}-b_{k}-\textbf{x}'_{t}\boldsymbol{\beta}).
\end{eqnarray}
Here, $\hat{b}_{i}$ is the estimator of $b^{*}_{\tau_{i}}$ and $y_{t}=\textbf{x}'_{t}\boldsymbol{\beta}^{*}+\varepsilon_{t}$, $t=1, \cdots, T$ with $\varepsilon_{t}$ being i.i.d. and having the same distribution as $\varepsilon$.

The asymptotic normality of $\hat{\boldsymbol{\beta}}^{clp}$ depends on the following conditions.

\begin{assumption}\label{ass2.1}
There is a $m\times m$ positive definite matrix $\textbf{C}$ such that
\begin{eqnarray}\label{eq1.5}
\lim_{T\rightarrow\infty}\frac{1}{T}\textbf{X}'\textbf{X}=\textbf{C},
\end{eqnarray}
where $\textbf{X}=(\textbf{x}_{1}, \cdots, \textbf{x}_{T})'$ is the $T\times m$ design matrix.
\end{assumption}

\begin{assumption}\label{ass2.2}
$E(|\varepsilon_{t}|^{2(p-1)})<\infty$, for $1<p\leq2$.
\end{assumption}

\begin{assumption}\label{ass2.3}
For $1<p\leq2$, there exists a positive constant $\delta>0$ such that $E(|\varepsilon_{t}-b|^{p-2})<\infty$ when $b\in U(b^{*}_{\tau_{k}}, \delta)$.
\end{assumption}

\begin{remark}\label{rem2.1}
It is well known that the conditions on which the quantile regression theory, such as the asymptotic normality, is built is extremely mild. Compared with those conditions, we remark that our conditions are really less restrictive. Indeed, Assumption \ref{ass2.1} is common, which is widely used in all kinds of regression methods
including quantile regression. The essential Assumption \ref{ass2.2} seems a little stronger but becomes negligible when $p$ approaches to 1 from above. At first glance, as a technical assumption, Assumption \ref{ass2.3} looks weird and strong but is valid at least when the true distribution of $\varepsilon_{t}$ has bounded density function near $b^{*}_{\tau_{k}}$. The boundedness of the probability density function of the error term is implicitly required in the quantile regression theory, see
Koenker (2005)\cite{Koenker}, Zou and Yuan (2008)\cite{ZOUY}, among others.
\end{remark}

\begin{theorem}\label{thm2.1}
Suppose $1<p\leq2$ and Assumptions \ref{ass2.1}-\ref{ass2.3} hold, then
$\sqrt{T}(\hat{\boldsymbol{\beta}}^{clp}-\boldsymbol{\beta}^{*})$ is asymptotically normal with mean
0 and covariance matrix
\begin{eqnarray}\label{eq1.6}
\boldsymbol{\Sigma}_{clp}=\textbf{C}^{-1}\frac{\sum^{K}_{k^{'}=1}\sum^{K}_{k=1}
E[\boldsymbol{\varphi}_{\tau_{k'},p}(\varepsilon-b^{*}_{\tau_{k'}})\boldsymbol{\varphi}_{\tau_{k},p}(\varepsilon-b^{*}_{\tau_{k}})]}
{(\sum_{k=1}^{K}E\boldsymbol{\psi}_{\tau_{k},p}(\varepsilon-b^{*}_{\tau_{k}}))^{2}},
\end{eqnarray}
where $\boldsymbol{\varphi}_{\tau,p}(s)=p|\tau-I(s<0)||s|^{p-1}\mbox{sign}(s)$ and $\boldsymbol{\psi}_{\tau,p}(s)=p(p-1)|\tau-I(s<0)||s|^{p-2}$.
\end{theorem}
\subsection{Asymptotic relative efficiency}
In order to consider the asymptotic relative efficiency of CLpQR, we need the limit version of the asymptotic variance matrix in \eqref{eq1.6} when the partition thinness for $(0,1)$, the range of $\tau$, converges to 0.
For the sake of convenience, we use the equally spaced weights, namely $\tau_{k}=k/(K+1)$, $k=1, 2, \cdots, K$, and get the following theorem.
\begin{theorem}\label{thm2.2}
We have, as $K\rightarrow\infty$,
\begin{eqnarray*}\label{eq1.20}
\nonumber & & \frac{\sum^{K}_{k^{'}=1}\sum^{K}_{k=1}
E[\boldsymbol{\varphi}_{\tau_{k'},p}(\varepsilon-b^{*}_{\tau_{k'}})\boldsymbol{\varphi}_{\tau_{k},p}(\varepsilon-b^{*}_{\tau_{k}})]}
{(\sum_{k=1}^{K}E\boldsymbol{\psi}_{\tau_{k},p}(\varepsilon-b^{*}_{\tau_{k}}))^{2}}\\
\nonumber &\longrightarrow &\frac{E_{\varepsilon_{b}}E_{\varepsilon_{c}}
E_{\varepsilon}((F_{\varepsilon, p}(\varepsilon_{c})-I(\varepsilon<\varepsilon_{c}))
(F_{\varepsilon, p}(\varepsilon_{b})-I(\varepsilon<\varepsilon_{b}))
|\varepsilon-\varepsilon_{b}|^{p-1}|\varepsilon-\varepsilon_{c}|^{p-1})}
{(p-1)^{2}(E_{\varepsilon_{a}}E_{\varepsilon}
(|F_{\varepsilon, p}(\varepsilon_{a})-I(\varepsilon<\varepsilon_{a})||\varepsilon-\varepsilon_{a}|^{p-2}))^{2}},
\end{eqnarray*}
where $\varepsilon_{a}$, $\varepsilon_{b}$ and $\varepsilon_{c}$ are three independent random variables with the identical cdf $F_{\varepsilon, p}$ such that its inverse function satisfies
\begin{eqnarray*}\label{eq1.21}
\frac{\int^{F^{-1}_{\varepsilon, p}(\tau)}_{-\infty}|r-F^{-1}_{\varepsilon, p}(\tau)|^{p-1}dF_{\varepsilon}(r)}
{\int^{\infty}_{-\infty}|r-F^{-1}_{\varepsilon, p}(\tau)|^{p-1}dF_{\varepsilon}(r)}=\tau,
\end{eqnarray*}
namely, $F^{-1}_{\varepsilon, p}(\tau)$ is the $\tau$th $L^{p}$-quantile of $\varepsilon$.
The subscript in the expectation sign $E$ indicates with respect to which random variable expectation is calculated.
\end{theorem}

We consider the asymptotic relative efficiency (ARE) of the CLpQR with respect to least square regression (LS). In order to compare CLpQR with CQR (the composite quantile regression developed by Zou and Yuan (2008)\cite{ZOUY}), a similar ARE of the CQR also be calculated. The asymptotic variance matrix of the LS is $\sigma^{2}C^{-1}$ when the error variance $\sigma^{2}<\infty$. So the
ARE of the CLpQR and CQR with respect to the LS can be calculated as follows.
\begin{eqnarray}\label{eq1.024}
\nonumber& & \mbox{ARE}_{CLpQR}\\
&=&\frac{\sigma^{2}(p-1)^{2}(E_{\varepsilon_{a}}E_{\varepsilon}
(|F_{\varepsilon, p}(\varepsilon_{a})-I(\varepsilon<\varepsilon_{a})||\varepsilon-\varepsilon_{a}|^{p-2}))^{2}}
{E_{\varepsilon_{b}}E_{\varepsilon_{c}}
E_{\varepsilon}((F_{\varepsilon, p}(\varepsilon_{c})-I(\varepsilon<\varepsilon_{c}))
(F_{\varepsilon, p}(\varepsilon_{b})-I(\varepsilon<\varepsilon_{b}))
(|\varepsilon-\varepsilon_{b}||\varepsilon-\varepsilon_{c}|)^{p-1})}.
\end{eqnarray}According to Theorem 3.1 in Zou and Yuan (2008)\cite{ZOUY}
\begin{eqnarray}\label{eq1.025}
\mbox{ARE}_{CQR}=\frac{1}{12(E(f(\varepsilon)))^2},
\end{eqnarray}
where $f$ is the density function of $\varepsilon$.
It is obvious according to the result in the next section that $\mbox{ARE}_{CLpQR}$ ($\mbox{ARE}_{CQR}$) are also the ARE of the CLpQR-oracle (CQR-oracle) with respect to the LS-oracle.

We consider two commonly-used distributions: a mixture of two normals and the generalized error distribution (GED).

\begin{figure}
\centering
\includegraphics[width=2.6 in,height=2.5in , angle=0]{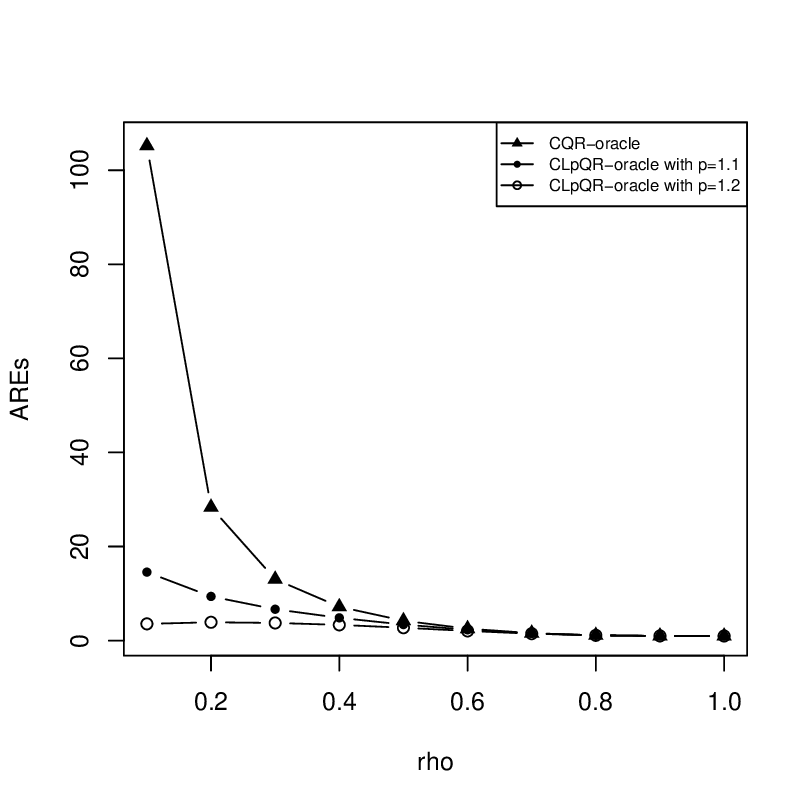}
\includegraphics[width=2.6 in,height=2.5in , angle=0]{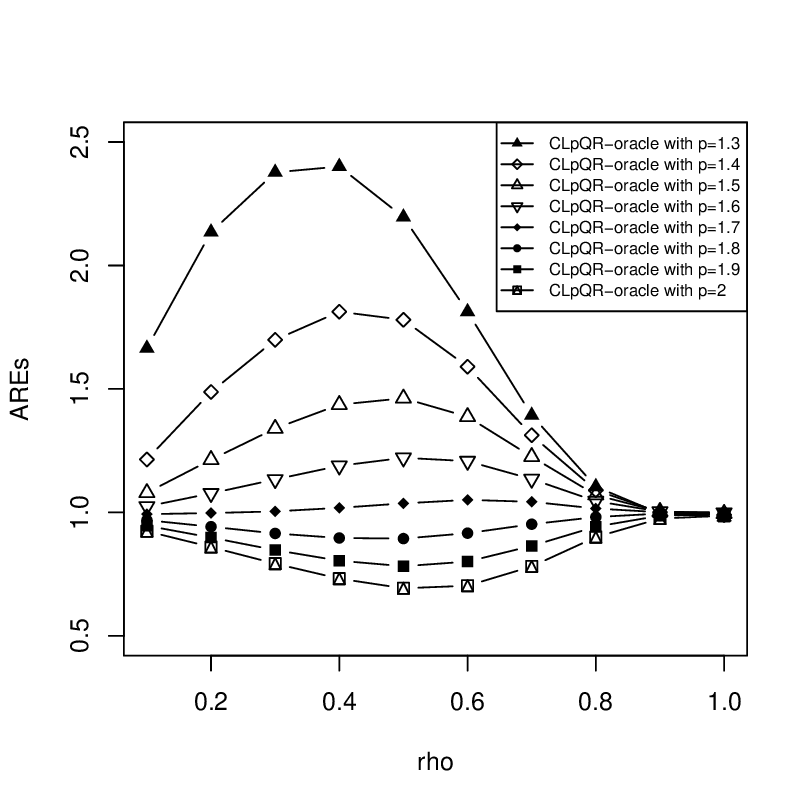}
\includegraphics[width=2.6 in,height=2.5in , angle=0]{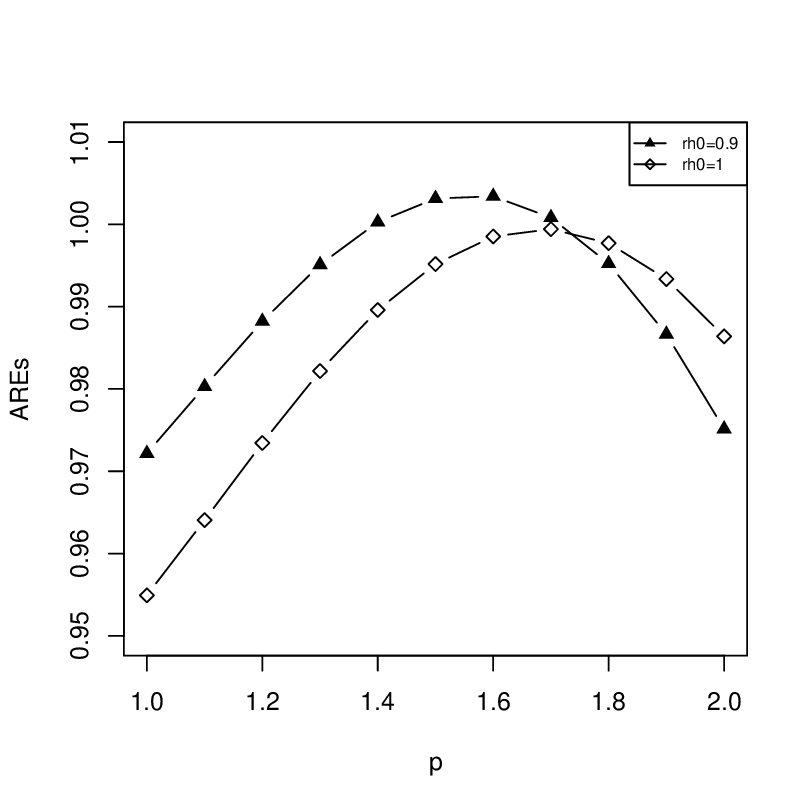}
\includegraphics[width=2.6 in,height=2.5in , angle=0]{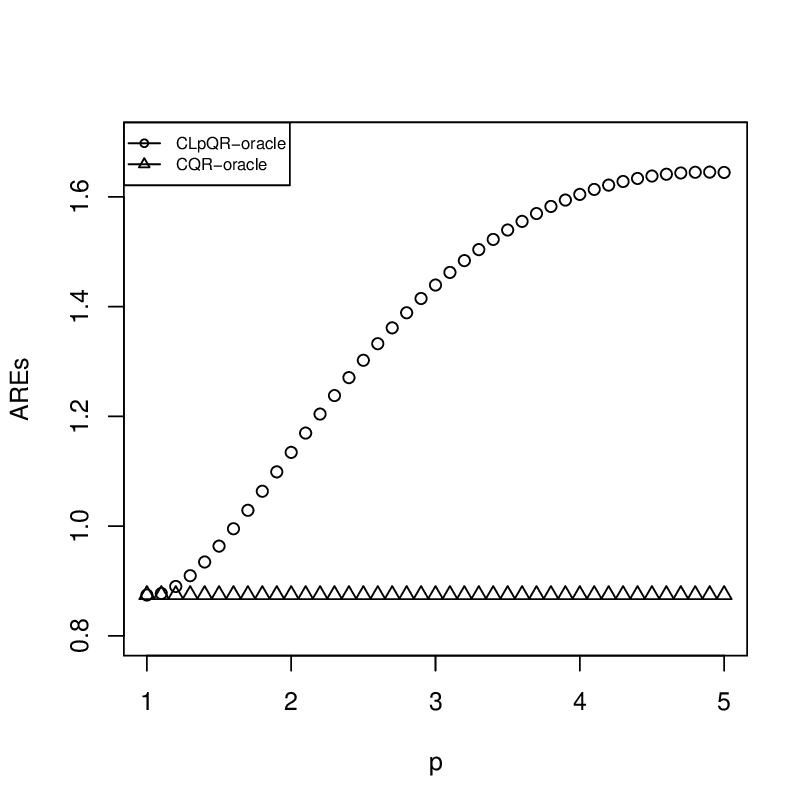}
\caption{Upper left panel: $ARE_{CQR}$ and $ARE_{CLpQR}$ ($p=1.2$ and 1.2) as the functions of
$\rho$ the mixture parameter of the mixture of two normals. Upper right panel: $ARE_{CLpQR}$ ($p\geq1.3$) as the functions of $\rho$. Lower left panel: $ARE_{CLpQR}$ ($\rho=0.9$ and 1) as the functions of
$p$. When $p=1$ $ARE_{CLpQR}$ is just $ARE_{CQR}$. Lower right panel: $ARE_{CLpQR}$ as the function of
$p$ when the error obeys the GED. The horizontal line marked by triangular indicates $ARE_{CQR}=0.8748277$ always.
} \label{fig1}
\end{figure}

{\bf Case 1} (a mixture of two normals) The error $\varepsilon$ has the density function
\begin{eqnarray*}\label{eq1.026}
\frac{1-\rho}{\sqrt{2\pi}}\exp\Big(-\frac{x^2}{2}\Big)+\frac{1}{\rho^2\sqrt{2\pi}}\exp\Big(-\frac{x^2}{2\rho^6}\Big)
\end{eqnarray*}
for $0<\rho<1$.
According to \eqref{eq1.025}, a precise function for $\mbox{ARE}_{CQR}$ is obtained
\begin{eqnarray*}\label{eq1.027}
\mbox{ARE}_{CQR}=\frac{3(1-\rho+\rho^{7})}{\pi}\Big((1-\rho)^{2}
+\frac{1}{\rho}+\frac{2\sqrt{2}\rho(1-\rho)}{\sqrt{1+\rho^{6}}}\Big)^{2}.
\end{eqnarray*}
A notable property about $\mbox{ARE}_{CQR}$ is that $\mbox{ARE}_{CQR}\rightarrow\infty$ as $\rho\rightarrow0$.
Using \eqref{eq1.024}, we calculate $\mbox{ARE}_{CLpQR}$ for $p\leq1.1$ and find it could also converge to infinity as $\rho\rightarrow0$ at a slower rate than $\mbox{ARE}_{CQR}$, see the upper left panel in Figure \ref{fig1}. But when in some value cases of $\rho$ for example $\rho=0.9, 1$, $\mbox{ARE}_{CLpQR}$ is larger than $\mbox{ARE}_{CQR}$ (the case $p=1$ corresponds to $\mbox{ARE}_{CQR}$), see the lower left panel in Figure \ref{fig1}. For $p\leq1.3$, the smaller the value of $p$ is, the smaller $\mbox{ARE}_{CLpQR}$ becomes, see the upper right panel in Figure \ref{fig1} for details.

{\bf Case 2} (the generalized error distribution) The density function of the GED is
\begin{eqnarray*}\label{eq1.028}
\frac{\beta}{2\alpha\Gamma(1/\beta)}\exp\Big(-\Big(\frac{|x|}{\alpha}\Big)^{\beta}\Big).
\end{eqnarray*}
Based on \eqref{eq1.025}, a precise function for $\mbox{ARE}_{CQR}$ is yielded
\begin{eqnarray*}\label{eq1.029}
\mbox{ARE}_{CQR}=\frac{3\beta^{2}}{4^{1/\beta}}\frac{\Gamma(3/\beta)}{(\Gamma(1/\beta))^3}.
\end{eqnarray*}
We set $\alpha=1$ and $\beta=5$. The lower right panel in Figure \ref{fig1} depicts that
$\mbox{ARE}_{CLpQR}$ keeps increasing with $p$ and is larger than $\mbox{ARE}_{CQR}=0.8748277$ uniformly
in $p\in(1,5)$.

While the above analysis is based on the limit version of the asymptotic relative efficiency
when $K\rightarrow\infty$, the ARE is empirically the same as its limit when $K=19$. So in the latter simulation and empirical analysis, we set $K=19$.
\section{The CLpQR-oracular estimation}
When the high-dimensional covariant has sparsity, Tibshirani (1996)\cite{TIBSHIRANI} invested the Lasso regression to select variables and estimate coefficients simultaneously. Fan and Li (2001)\cite{Fan} considered the SCAD-penalized least square regression and discussed its oracle properties. Zou (2006)\cite{ZOU} used the reciprocal of LS estimates to differentially tune the penalization intensity, called it the adaptive lasso and proved its oracle properties.

In this section, following the tack of Zou (2006)\cite{ZOU} we develop a penalized composite $L^{p}$-quantile regression method. Define the estimator as follows.
\begin{eqnarray}\label{eq1.24}
(\hat{b}_{1}, \cdots, \hat{b}_{K}, \hat{\boldsymbol{\beta}}^{Aclp})
=\arg\min_{b_{1}, \cdots, b_{K}, \boldsymbol{\beta}}\sum^{K}_{k=1}\sum^{T}_{t=1}\boldsymbol{\eta}_{\tau_{k},p}(y_{t}-b_{k}-\textbf{x}'_{t}\boldsymbol{\beta})+\lambda\sum^{m}_{j=1}
\frac{|\beta_{j}|}{|\hat{\beta}^{clp}_{j}|^{2}},
\end{eqnarray}
where these $\hat{\beta}^{clp}_{j}$ are the non-penalized CLpQR estimators.
The following theorem shows that the adaptively penalized CLpQR estimator
also enjoys the oracle properties.
\begin{theorem}\label{thm3.1}
Suppose the conditions in Theorem \ref{thm2.1} are satisfied.
Let $\lambda$ be the function of $T$, namely $\lambda=\lambda(T)$. If  $\frac{\lambda(T)}{\sqrt{T}}\rightarrow0$, and $\lambda(T)T^{\frac{p-2}{2}}\rightarrow\infty$ as $T\rightarrow\infty$,
then we have, for $\hat{\boldsymbol{\beta}}^{Aclp}$,

1. Consistency in selection: $P(\{j: \hat{\boldsymbol{\beta}}^{Aclp}\neq0\}=\mathcal{A})\rightarrow1$.

2. Asymptotic normality: $\sqrt{T}(\hat{\boldsymbol{\beta}}_{\mathcal{A}}^{Aclp}-\boldsymbol{\beta}_{\mathcal{A}}^{*})\rightarrow N(0, \boldsymbol{\Sigma}_{Clp_{oracle}})$.

\noindent Here, $\mathcal{A}=\{j: \beta^{*}_{j}\neq0\}$ and the vector $\boldsymbol{\beta}_{\mathcal{A}}^{*}$ consists of those nonzero components of $\boldsymbol{\beta}^{*}$.
\end{theorem}

\begin{remark}\label{rem3.1}
In \eqref{eq1.24}, we can also consider the SCAD penalty and the oracle properties similar to those in Theorem 3.1 should also hold. The main reason why choosing the adaptive lasso is that a unified algorithm for $p\geq1$ is easy to construct in the case.
\end{remark}
\section{Nearly quantile regression}
In this section we instead consider the data model as follows.
\begin{eqnarray}\label{eq1.37}
y_{t}=\textbf{x}'_{t}\boldsymbol{\beta}_{0}+u_{t}, t=1, 2, \cdots, T,
\end{eqnarray}
for observed $\{\textbf{x}_{t}\}$, unknown $\boldsymbol{\beta}_{0}\in R^{m}$ and i.i.d. unknown errors $\{u_{t}\}$ with the distribution density $f(u)$ being continuous in a neighborhood of 0. Let $u_{t}$'s $\tau$th quantile be zero and thus the conditional $\tau$th quantile of $y$ denoted by $\textbf{x}'_{t}\boldsymbol{\beta}(\tau)$ is just $\textbf{x}'_{t}\boldsymbol{\beta}_{0}$. The $\tau$th $L^{p}$-quantile of $u_{t}$ is denoted by $q^{lp}_{u}(\tau)$. The near quantile regression estimator of $\boldsymbol{\beta}(\tau)$ is
\begin{eqnarray}\label{eq1.38}
\hat{\boldsymbol{\beta}}_{T, p}(\tau)
=\arg\min_{\boldsymbol{\beta}}\Bigg\{\frac{1}{T}\sum^{T}_{t=1}(\boldsymbol{\eta}_{\tau,p}(y_{t}-\textbf{x}'_{t}\boldsymbol{\beta})
-\boldsymbol{\eta}_{\tau,p}(y_{t}))\Bigg\},
\end{eqnarray}
for $p\in(1, \epsilon)$, $\epsilon$ is a small positive number.
The objective function of near quantile regression in \eqref{eq1.38} is smooth, which is an
appealing property. Although quantile regression (Koenker and Bassett (1978)\cite{Koenker Gilbert}) is a powerful statistical learning tool, the un-smoothness of its objective function is sometimes
an obstacle in developing statistics theory and methods. The literature are devoted to smooth the objective function of quantile regression. Horowitz (1998)\cite{Horowitz} used an analogous to the integral of a kernel function to smooth the objective function of $L^{1}$ regression in order to apply the standard theory of the bootstrap. Fernandes et al. (2021)\cite{Fernandes} proposed  a convolution-type smoothing method to produce a continuous QR estimator which saves from the curse of
dimensionality. He et al. (2023)\cite{He} considered a convolution smoothed approach that
achieves adequate approximation to computation and inference for high-dimensional quantile regression. These ``kernel function" approaches are sophisticated as for example they involve the intractable bandwidth selection.
Built on this growing literature, the near quantile regression is more manageable and serves as a natural smoothness scheme.

The asymptotic property of the estimator is established on the following conditions.
In the descriptions of these conditions, $\triangle$ is an arbitrarily small positive constant.
\begin{assumption}\label{ass4.1}
There is a $m\times m$ positive definite matrix $\textbf{D}_{0}$ such that
\begin{eqnarray*}\label{eq1.39}
\lim_{T\rightarrow\infty}\frac{1}{T}\sum^{T}_{t=1}\textbf{x}_{t}\textbf{x}'_{t}=\textbf{D}_{0}.
\end{eqnarray*}
\end{assumption}

\begin{assumption}\label{ass4.2}
$E(|u_{t}|^{2(p-1)})<\infty$, for $p\in (1, \triangle)$.
\end{assumption}

\begin{assumption}\label{ass4.3}
For $p\in (1, \triangle)$, there exists a positive constant $\delta>0$ such that when $b\in U(q^{lp}_{u}(\tau), \delta)$, $E(|u_{t}-b|^{p-2})<\infty$.
\end{assumption}

\begin{assumption}\label{ass4.4}
For $p\in (1, \triangle)$, $f(x)$ has the first derivative function $f^{(1)}(x)$ such that
$$\int_{-\infty}^{+\infty}|x|^{p-1}|f^{(1)}(x)|dx<\infty.$$
\end{assumption}

Assumptions \ref{ass4.1}-\ref{ass4.3} are similar to Assumptions \ref{ass2.1}-\ref{ass2.3} and Assumption \ref{ass4.4} is a key technical assumption.

\begin{theorem}\label{thm4.1}
Under the model \eqref{eq1.37} and Assumptions \ref{ass4.1}-\ref{ass4.4}, we have
\begin{eqnarray*}\label{eq1.40}
\lim_{T\rightarrow\infty
\atop p\rightarrow1+}\sqrt{T}(\hat{\boldsymbol{\beta}}_{T, p}(\tau)-\boldsymbol{\beta}(\tau))\stackrel{\textit{D}}{\longrightarrow}
N(0, \boldsymbol{\Sigma}_{0}),
\end{eqnarray*}
with $\boldsymbol{\Sigma}_{0}=\tau(1-\tau)f^{-2}(0)\textbf{D}_{0}^{-1}$.
\end{theorem}

\begin{remark}\label{rem4.1}
Theorem \ref{thm4.1} shows that the near quantile regression estimator is asymptotically
equivalent to the standard QR estimator. The challenge of Theorem \ref{thm4.1} is that it makes sure
that this convergence holds when $p\rightarrow1+$ and $T\rightarrow\infty$ in any way, not just
successively for example first with respect to $T$ and then $p$. The convergence in the latter way is easier to prove.
\end{remark}

For the least absolute deviations estimator, i.e. the $L^{1}$ regression estimator, we have its smoothed version: the near $L^{1}$ regression estimator, namely $\hat{\boldsymbol{\beta}}_{T, p}(0.5)$ obtained by \eqref{eq1.38} with $\tau=0.5$ for which the asymptotic property is the following corollary of Theorem \ref{thm4.1}.

\begin{corollary}\label{coro1}
When $\tau=0.5$, under the model \eqref{eq1.37} and Assumptions \ref{ass4.1}-\ref{ass4.4}, we have
\begin{eqnarray*}\label{eq1.40}
\lim_{T\rightarrow\infty
\atop p\rightarrow1+}\sqrt{T}(\hat{\boldsymbol{\beta}}_{T, p}(0.5)-\boldsymbol{\beta}(0.5))\stackrel{\textit{D}}{\longrightarrow}
N(0, \boldsymbol{\Sigma}_{0,0.5}),
\end{eqnarray*}
with $\boldsymbol{\Sigma}_{0,0.5}=0.25 f^{-2}(0)\textbf{D}_{0}^{-1}$.
\end{corollary}

We next consider an new estimate of the asymptotic variance matrix $\boldsymbol{\Sigma}_{0}$ in Theorem \ref{thm4.1} and define the estimator as follows.
\begin{eqnarray}\label{eq1.70}
\hat{\boldsymbol{\Sigma}}_{0}=\frac{\tau(1-\tau)}{(\frac{1}{T}\sum^{T}_{t=1}\boldsymbol{\psi}_{\tau, p}(y_{t}-x'_{t}\hat{\boldsymbol{\beta}}_{T, p}(\tau)))^{2}}\Big(\frac{1}{T}\sum^{T}_{t=1}\textbf{x}_{t}\textbf{x}'_{t}\Big)^{-1},
\end{eqnarray}
for $p$ very close to 1 from the above, where $\boldsymbol{\psi}_{\tau,p}(s)=p(p-1)|\tau-I(s<0)||s|^{p-2}$.
The consistency of the estimator is proved in Theorem \ref{thm4.2} under the following conditions.

\begin{assumption}\label{ass4.5}
Let $\boldsymbol{\beta}_{p}(\tau)$ be the population counterpart of $\hat{\boldsymbol{\beta}}_{T, p}(\tau)$. There exists a close neighbourhood of $\boldsymbol{\beta}_{p}(\tau)$, i.e. $U[\boldsymbol{\beta}_{p}(\tau), \textbf{r}_{1}]$ such that, for $T\rightarrow\infty$
\begin{eqnarray}\label{eq1.70+1}
\sup_{\boldsymbol{\delta}\in U[\boldsymbol{\beta}_{p}(\tau), \textbf{r}_{1}]}\frac{1}{T}\Big|\sum^{T}_{t=1}\boldsymbol{\psi}_{\tau, p}(y_{t}-\textbf{x}'_{t}\boldsymbol{\delta})-E\boldsymbol{\psi}_{\tau, p}(y_{t}-\textbf{x}'_{t}\boldsymbol{\delta})\Big|\stackrel{\textit{P}}{\longrightarrow}0.
\end{eqnarray}
\end{assumption}

\begin{assumption}\label{ass4.6}
$E(|u_{t}+\textbf{x}'_{t}\boldsymbol{\delta}|^{p-2})$ is continuous with respect to $\boldsymbol{\delta}$ in a close neighbourhood $U[\boldsymbol{\beta}_{p}(\tau), \textbf{r}_{2}]$ uniformly for all $\textbf{x}'_{t}$, where $\boldsymbol{\beta}_{p}(\tau)$
is equal to $\boldsymbol{\beta}_{p}(\tau)=\boldsymbol{\beta}_{0}+q^{lp}_{u}(\tau)\textbf{e}$ with $\textbf{e}$ being a vector with its first component being 1 and the others 0.
\end{assumption}

Assumption \ref{ass4.5} is a uniform version of Khinchin's law of large numbers.
In Assumption \ref{ass4.6}, the uniform continuity will hold at least when $\textbf{x}_{t}$ are bounded and the bounedness of $\textbf{x}_{t}$ is often used to consider the regression with non-random covariates.

\begin{theorem}\label{thm4.2}
Under the model \eqref{eq1.37} and Assumptions \ref{ass4.1}-\ref{ass4.6}, we have
\begin{eqnarray*}\label{eq1.71}
\lim_{p\rightarrow1+}\lim_{T\rightarrow\infty}\hat{\boldsymbol{\Sigma}}_{0}=\boldsymbol{\Sigma}_{0}\ \mbox{in \ probability.}
\end{eqnarray*}
\end{theorem}

While the existing estimation methods are almost non-parametric, this theorem provides a new consistent parametric estimation for the asymptotic covariance matrix of quantile regression.

\section{Algorithm}
In this section, we consider the computational aspect of the proposed regression methods. We minimize the following objective function of composite $L^{p}$-quantile regression with penalty.
\begin{eqnarray}\label{eq1.24+1}
\min_{b_{1}, \cdots, b_{K}, \boldsymbol{\beta}}\frac{1}{T}\sum^{K}_{k=1}\sum^{T}_{t=1}\eta_{\tau_{k},p}(y_{t}-b_{k}-\textbf{x}'_{t}\boldsymbol{\beta})
+\sum^{m}_{j=1}w_{j}|\beta_{j}|,
\end{eqnarray}
where the penalty coefficients $w_{j}\geq0, j=1,2,\cdots,m$. The setting of penalty terms is very general: Setting $w_{j}=0$ leaves $\beta_{j}$ unpenalized and identical $w_{j}$ corresponds to an analogue of Lasso. Since the nondifferentiability of the penalty term makes the gradient descent infeasible we apply a combination of the cyclic coordinate descent (Tseng(2001)\cite{TSENG}) and proximal gradient algorithms (Parikh and Boyd(2013)\cite{PARIKH})(CCPA for short). A similar thought was utilized by Gu and Zou(2016)\cite{Gu} but their algorithm cannot apply to generic $L^p$ regression. Below is a description of the proposed algorithm in detail.

Let $\boldsymbol{\alpha}=(\alpha_{1},\alpha_{2},\cdots,\alpha_{K+m})'$ with
$(\alpha_{1},\alpha_{2},\cdots,\alpha_{K})'=(b_{1}, b_{2},\cdots, b_{K})'$ and\\
$(\alpha_{K+1},\alpha_{2},\cdots,\alpha_{K+m})'=(\beta_{1}, \beta_{2},\cdots, \beta_{m})'$.
Rewrite \eqref{eq1.24+1} as
\begin{eqnarray}\label{eq1.24+2}
\min_{\boldsymbol{\alpha}}\frac{1}{T}O_{T,K}(\boldsymbol{\alpha})
+\sum^{m+K}_{j=1}w_{j}|\alpha_{j}|.
\end{eqnarray}
Let $\boldsymbol{\alpha}^{q}=(\alpha^{q}_{1},\alpha^{q}_{2},\cdots,\alpha^{q}_{K+m})'$ stand for the update of $\boldsymbol{\alpha}$ after the $q$-th cycle of the coordinate descent algorithm. For convenience, write
\begin{eqnarray*}\label{eq1.24+3}
\boldsymbol{a}^{q+1}_{-i}&=&(\alpha_{1}^{q+1},\cdots,\alpha_{i-1}^{q+1},
\alpha_{i+1}^{q},\cdots,\alpha_{K+m}^{q})', 1\leq i\leq K+m, q\geq0.\\
\boldsymbol{\beta}^{q+1}_{-i}&=&(\beta_{1}^{q+1},\cdots,\beta_{i-1}^{q+1},
\beta_{i+1}^{q},\cdots,\beta_{m}^{q})', 1\leq i\leq m, q\geq0.
\end{eqnarray*}
According to the coordinate descent algorithm, updating $\alpha_{i}$ is equivalent to minimizing the objective function:
\begin{eqnarray}\label{eq1.24+4}
\min_{\alpha_{i}}\frac{1}{T}O_{T,K}(\alpha_{i},\boldsymbol{a}^{q+1}_{-i})
+w_{i}|\alpha_{i}|,
\end{eqnarray}
where
\begin{equation*}\label{eq1.1}
O_{T,K}(\alpha_{i},\boldsymbol{a}^{q+1}_{-i})=
\begin{cases}
\sum^{i-1}_{k=1}\sum^{T}_{t=1}\boldsymbol{\eta}_{\tau_{k},p}(y_{t}-b^{q+1}_{k}-\textbf{x}'_{t}\boldsymbol{\beta}^{q})\\
\ \ \ +\sum^{T}_{t=1}\boldsymbol{\eta}_{\tau_{i},p}(y_{t}-b_{i}-\textbf{x}'_{t}\boldsymbol{\beta}^{q})\\
\ \ \ +\sum^{K}_{k=i+1}\sum^{T}_{t=1}\boldsymbol{\eta}_{\tau_{k},p}(y_{t}-b^{q}_{k}-\textbf{x}'_{t}\boldsymbol{\beta}^{q}) &i\leq K\\
\sum^{K}_{k=i}\sum^{T}_{t=1}\boldsymbol{\eta}_{\tau_{k},p}(y_{t}-b^{q+1}_{k}-\textbf{x}'_{t,-i}
\boldsymbol{\beta}_{-i}^{q+1}-x_{t,i}\beta_{i}) &i>K
\end{cases}
\end{equation*}
and $w_{i}=0$ when $i\leq K$. Denote $O'_{T,K}(\alpha_{i},\boldsymbol{a}^{q+1}_{-i})$ the first derivative of
$O_{T,K}(\alpha_{i},\boldsymbol{a}^{q+1}_{-i})$ with respect to $\alpha_{i}$ and let $S_{i}=c_{1}$ when $i\leq K$ or $S_{i}=c_{2}T^{-1}\|(x_{i,1},\cdots x_{i,T})\|^{2}$ otherwise. The proximal gradient method solves problem \eqref{eq1.24+4} by the iteration formula as follows.
\begin{eqnarray}\label{eq1.24+5}
\alpha_{i}^{q,0}=\alpha_{i}^{q}, \ \alpha_{i}^{q,d+1}
=L_{S^{-1}_{i}w_{i}}(\alpha_{i}^{q,d}-S^{-1}_{i}O'_{T,K}(\alpha^{q,d}_{i},\boldsymbol{a}^{q+1}_{-i}))
\end{eqnarray}
where $L_{u}(v)=\mbox{sign}(v)((|v|-u)I(|v|-u>0))$ serves as the soft threshold operator, $w_{i}=0$ for $i\leq K$ and $w_{i}=\lambda /|\hat{\beta}^{clp}_{i}|^2$ otherwise. We run \eqref{eq1.24+5} for $S$ iterations till meeting precision requirement and have $\alpha_{i}^{q+1}=\alpha_{i}^{q,S}$.

We have two remarks on the algorithm.
\begin{remark}\label{rem3}
During the implementation of the CCPA, setting constants $c_{1}$ and $c_{2}$ is very crucial. Empirically, we found that $c_{1}$ and $c_{2}$ taking values near 1.6 and 10 is a good choose.
Moreover, when $p<1.5$ in the CLpQR lose function, we need a adaptive step width, i.e. multiplying $S_{i}$ by $c_{3}$ in each iteration for iteration becomes slower at this moment. Empirically, we found that $c_{3}=0.9^{-1}$ or so works very well.
\end{remark}
\begin{remark}\label{rem3}
As special cases included in CLpQR, CQR and QR are robust against outliers and can be implemented for heavy-tailed or skewed response distributions without correctly specifying the likelihood. However, when applied to large-scale problems: large sample size and high dimension, QR computation via the linear program and interior point algorithm is prone to be slow or too high memory-consuming, which makes QR computation infeasible in a personal computer and could make QR less attractive compared to other machine learning tools (He et al.(2023)\cite{He}). In the simulation and empirical analysis, our proposed algorithm can be used to fit CQR and QR effectively. The algorithm turns out to be an practicable alternative of the commonly used liner programming and interior point algorithm when fitting quantile regression, especially in the high-dimension regime.
\end{remark}
\section{Simulation study}
\subsection{Simulation for composite Lp-quantile regression}
In the section we provide a comparison of CLpQR-oracle with CQR-oracle by Monte Carlo simulation and a comparison of our proposed algorithm with liner programming algorithm when calculating CQR-oracle as well.
The data generating process is
\begin{eqnarray}\label{eq1.24+6}
\textbf{y}=\boldsymbol{X}'\boldsymbol{\beta}^{*}+\textbf{u}
\end{eqnarray}
where $\boldsymbol{\beta}^{*}=(3,1.5,0,0,2,0,0,0)'$ and predictor vector $\textbf{x}$ comes from a multivariate normal distribution $N(\boldsymbol{0}, (0.5^{|i-j|})_{8\times8})$. The model was often used to example high-dimension statistics modelling by many authors, such as Tibshirani (1996)\cite{TIBSHIRANI} and Fan and Li (2001)\cite{Fan}.
We consider four common error distribution examples: E1. $N(0,9)$, E2. T-distribution with 3 degrees of freedom, E3. Cauchy, and E4. the generalized error distribution with the density function $(1/(2\Gamma(1+1/r)))\exp(-|x|^r) $ with $r=4$.
In each distribution case we generate 200 observations consisting of 100 observations for training model and another 100 ones for selecting the penalty parameters. In each case we repeat 100 times to evaluate the performance of the various methods and algorithms through comparing their estimation errors and variable selection results.
The estimation error is defined as
\begin{eqnarray*}\label{eq1.24+7}
EE=E((\hat{\boldsymbol{\beta}}^{Aclp}-\boldsymbol{\beta}^{*})'(0.5^{|i-j|})_{8\times8}(\hat{\beta}^{Aclp}-\beta^{*})).
\end{eqnarray*}
The variable selection result is described by the notation (ANC, ANIC) where ANC denotes the average number of non-zero components of estimate vector $(\hat{\beta}^{Aclp}_{1}, \hat{\beta}^{Aclp}_{2}, \hat{\beta}^{Aclp}_{5})$ and ANIC the average number of non-zero components of estimate vector $(\hat{\beta}^{Aclp}_{3}, \hat{\beta}^{Aclp}_{4}, \hat{\beta}^{Aclp}_{6}, \hat{\beta}^{Aclp}_{7}, \hat{\beta}^{Aclp}_{8})$. Simulation results are collected in Table 6.1.
\begin{table}
\begin{center}
Table 6.1. Simulation results for models with various error distributions
\begin{tabular}{cccccccccccccccccccccccccccccccccccccccccccc}
\hline
&Items& LPS &CCPA&  &  &  &  &  \\
\hline
E1 & p&	1	 & 1 & 	1.001&	1.1	 & 1.5	&1.9 &  2   \\
                &EE	&0.3587	&0.3408&	0.3410&	0.3374	&0.3173&	0.3109&	0.3086\\
                &(ANC, ANIC)	&(3, 0.82)&	(3, 1.1)&	(3, 1.1)	&(3, 1.09)	&(3, 1.09)	&(3, 1.16)&	(3, 1.16)\\
\hline
E2 & p&	1&	1	&1.001	&1.1	&1.5	&1.9	&2\\
                     &EE&	0.0554&	0.0498	&0.0496&	0.0528	&0.0717&	0.1746&0.2341\\
                     &(ANC, ANIC)	&(3, 1.09)&	(3, 1.19)	&(3, 1.18)	&(3, 1.22)	&(3, 1.36)	&(3, 1.39)	&(3, 1.49)\\
\hline
E3 &p	&1&	1&	1.001	&1.1	&1.5	&1.9	&2\\
                &EE&	0.1150&	0.0974&	0.0987	&0.1290	&1.0721&	271.4711&	587.2584\\
                &(ANC, ANIC)	&(3, 0.73)	&(3, 1.04)&	(3, 1.08)	&(3, 1.33)	&(3, 1.58)	&(3, 1.61)&	(3, 1.91)\\
\hline
E4 &p	&1	&1	&1.001	&1.1	&1.5	&1.9&	2.5\\
&EE	&0.0125&	0.0115	&0.0115	&0.0114	&0.0107&	0.0095	&0.0087\\
&(ANC, ANIC)	&(3, 0.82)&	(3, 1.05)&	(3, 1.06)	&(3, 1.05)&	(3, 1.12)	&(3, 0.98)&	(3, 0.96)\\
\hline
\end{tabular}
\end{center}
\end{table}

The third column contains the results obtained by using the standard linear program solver (LPS) when $p=1$ (corresponding to CQR-oracle) and columns 4-9 collect those results got by the algorithm in Section 5. Across all examples, there are some phenomena in common. For $p=1$, i.e. when calculating CQR-oracle, CCPA is not only successfully applied to deal with the issue but also tends to yield smaller estimation error than LPS. The variable selection results show that LPS is apt to give a little smaller estimation of coefficients in absolute sense than CCPA. Moreover, the results for $p=1$ and $p=1.001$ are very close when using CCPA, which shows that the numerical experiments agree with the near quantile theory. Specifically, in example one, the smallest estimation error unsurprisingly appears when $p=2$. In example 3, when $p\geq1.5$ Assumption \ref{ass2.2} does not hold and hence the asymptotic variance will diverge, which is supported by simulation results as well.
In the case of the generalized error distribution, we find that the estimation error keeps decreasing when $p$ increases and the change is substantial, for example $p=2.5$, the estimation error decreases 29.6\% compared with $p=1$.
\subsection{Simulation for near quantile regression}
The data generating process is
\begin{eqnarray}\label{eq1.24+6}
\textbf{y}=\boldsymbol{X}'\beta_{0}+\textbf{u}
\end{eqnarray}
where $\beta_{0}=(\beta_{0},\beta_{1},\beta_{2},\beta_{3})'=(0.5,3,1.5,2)'$, predictor vector $(1, \textbf{x}')'$ and $\textbf{x}$ comes from a multivariate normal distribution $N(\boldsymbol{0}, (0.5^{|i-j|})_{3\times3})$.
For $\textbf{u}$, we consider three error distribution examples: E1. $N(0,3)$, E2. T-distribution with 3 degrees of freedom and E3. Chi-square distribution with 3 degrees of freedom. We substitute $\hat{\boldsymbol{\beta}}_{T, p}(\tau)$ and $\hat{\boldsymbol{\Sigma}}_{0}$ calculated by \eqref{eq1.38} and \eqref{eq1.70} into and mainly examine the following expression
\begin{eqnarray}\label{eq1.24+6+88}
\sqrt{T}(\hat{\boldsymbol{\beta}}_{T, p}(\tau)-\beta_{0})\hat{\boldsymbol{\Sigma}}_{0}^{-1/2}
\end{eqnarray}
in finite sample cases.
We set $p$, $\tau$ and $T$ take values in sets $\{1.5, 1.1, 1.01, 1.001\}$, $\{\tau=0.01, 0.05, 0.5\}$ and $\{T=100, 1000, 5000\}$, respectively. The quantile-quantile plots are obtained by 100 repetitions and the distributions of the last three components of \eqref{eq1.24+6+88} have been depicted in the plots. Here only the $\tau=0.5$ case
is showed for the normal distribution and T-distribution errors and the other cases are collected in the ``Supplementary materials".
These finite sample distribution results exemplify that the normalized estimators $\hat{\boldsymbol{\beta}}_{T, p, 1}(\tau)$, $\hat{\boldsymbol{\beta}}_{T, p, 2}(\tau)$ and $\hat{\boldsymbol{\beta}}_{T, p, 3}(\tau)$ of $\beta_{1}$, $\beta_{2}$ and $\beta_{3}$ indeed converge to the standard normal random variable in distribution as $T\rightarrow\infty,
p\rightarrow1+$ simultaneously. Specifically, in the normal error case, when $T=100$, the $L^p$-quantile
regression estimator when $p=1.1$ already approximates its quantile counterpart well. In the three error cases, when $T=1000$, the $L^p$-quantile regression estimator when $p=1.001$ excellently approximates its quantile counterpart.

When $\tau$ deviates from $0.5$, especially when taking extremal levels, for example $\tau=0.01$, the convergence rate will slow down because of the sparsity of tail data, see ``Supplementary materials" for much more figures. But it doesn't really matter. 
The difficulty lies in estimating $f(0)$ when estimating the asymptotic covariance matrix of quantile regression. So we just set $\tau=0.5$ when using the denominator in \eqref{eq1.70} to estimate $f(0)$.

\begin{figure}
\centering
\includegraphics[width=2.1 in,height=2in , angle=0]{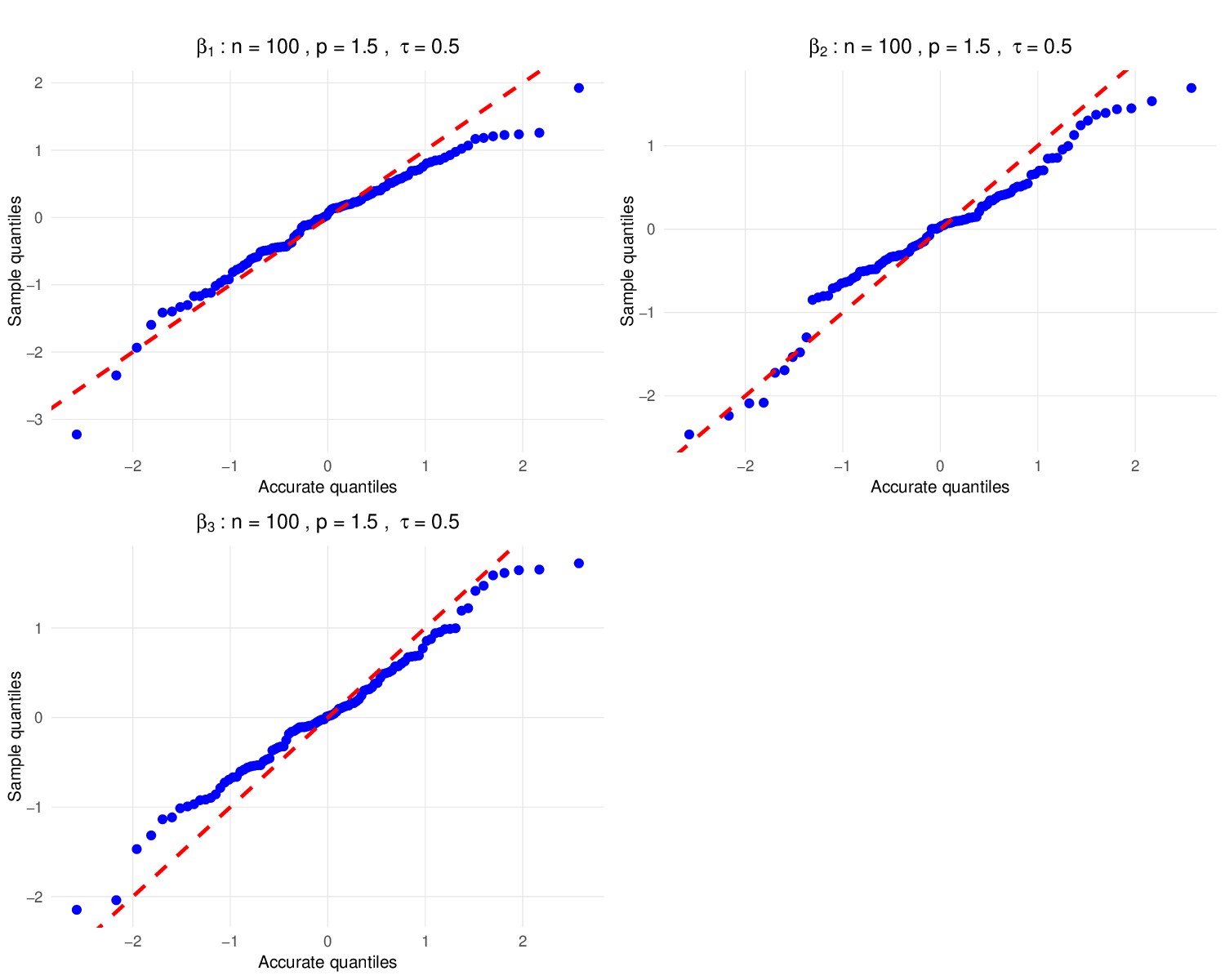}
\includegraphics[width=2.1 in,height=2in , angle=0]{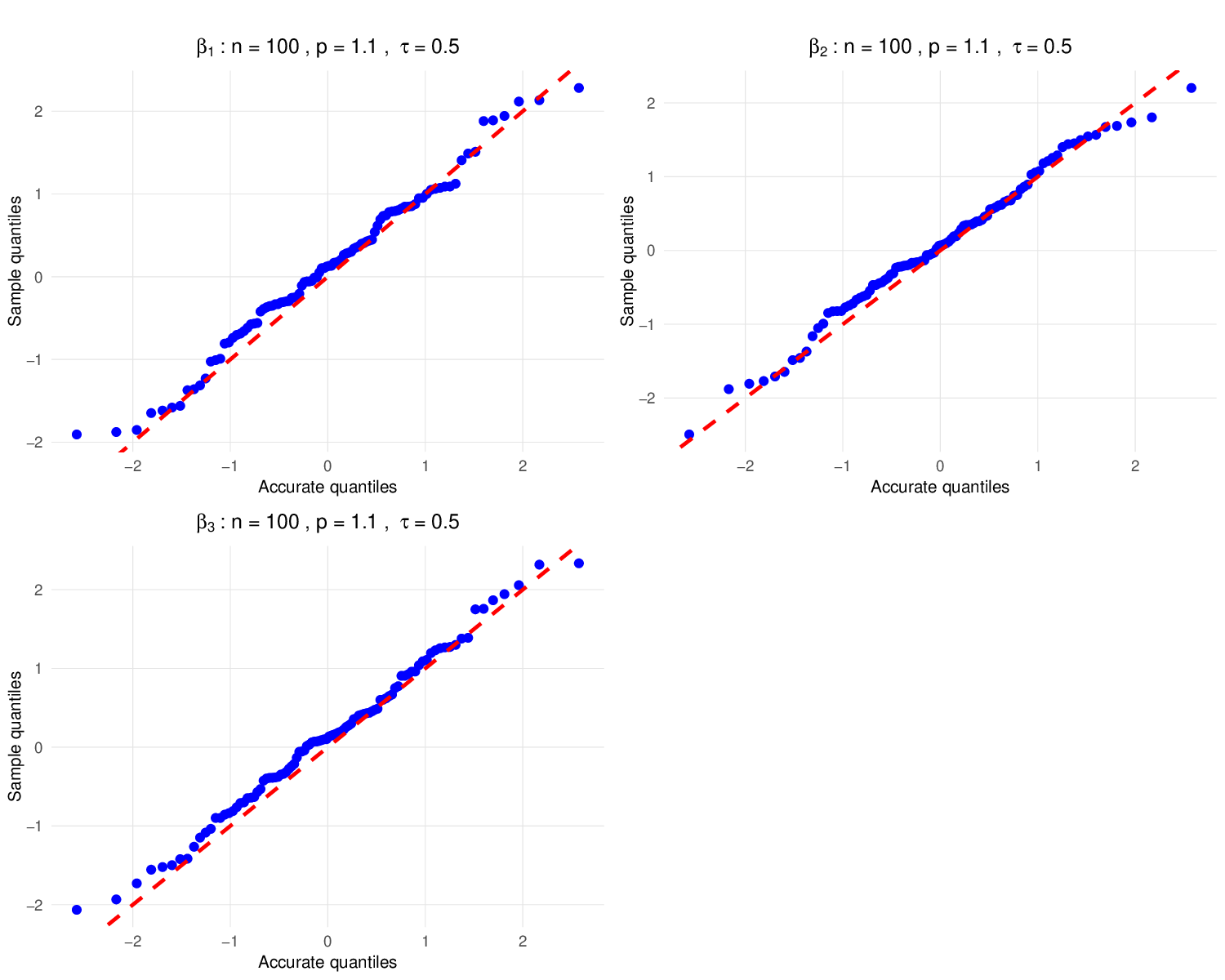}
\includegraphics[width=2.1 in,height=2in , angle=0]{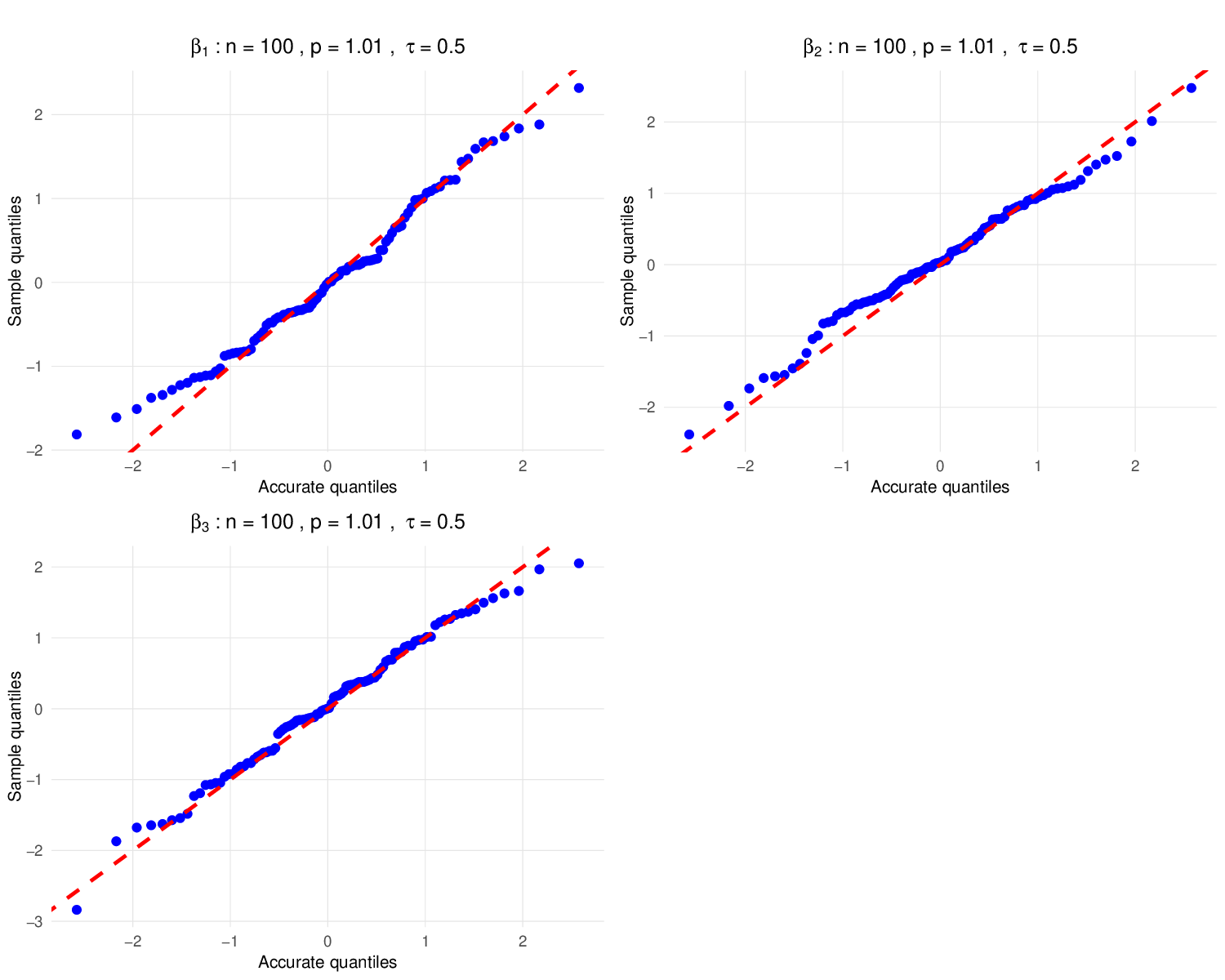}
\includegraphics[width=2.1 in,height=2in , angle=0]{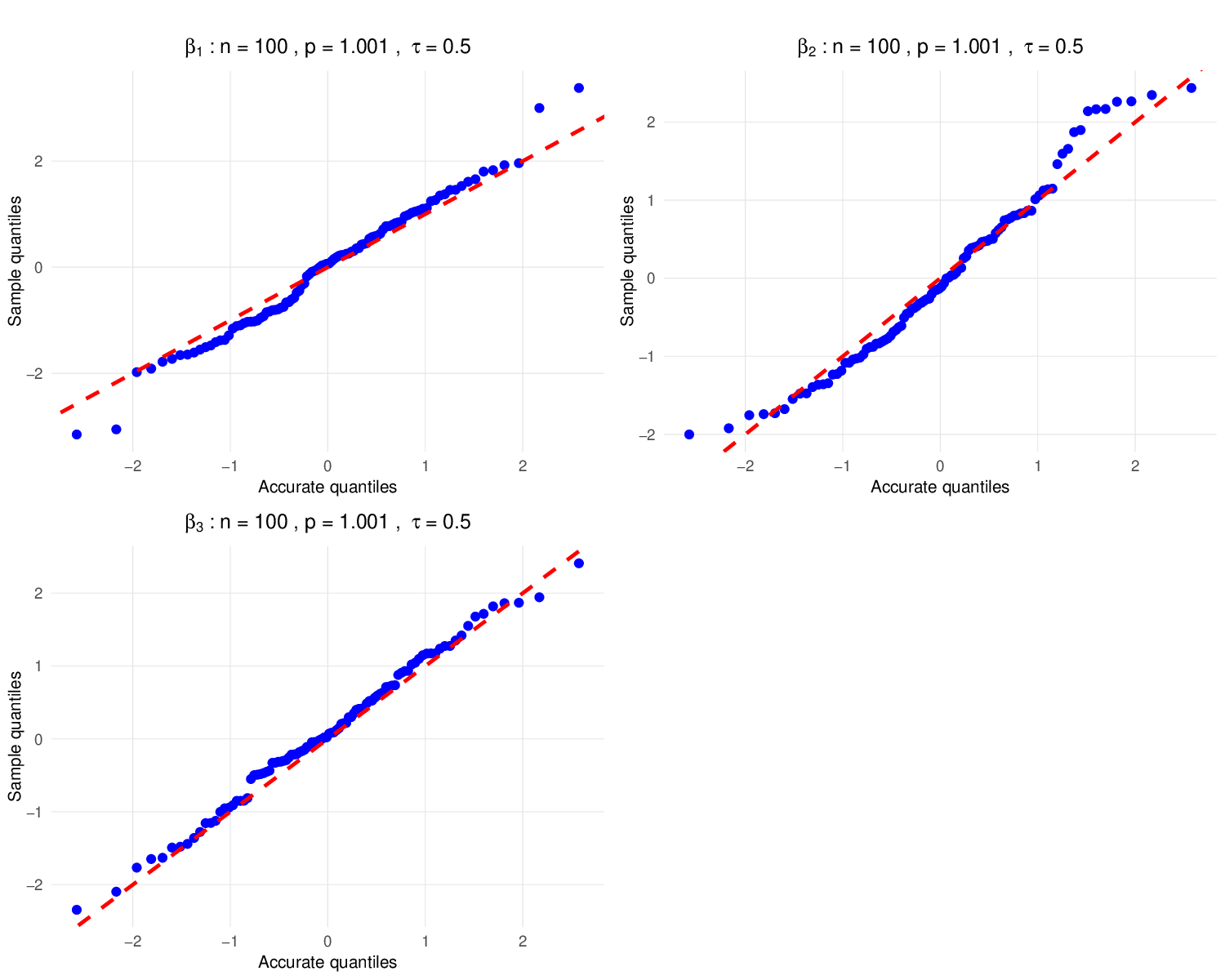}
\caption{The Normal error, the case of $T=100$, Upper left panel: $p=1.5$. Upper right panel: $p=1.1$. Lower left panel: $p=1.01$. Lower right panel: $p=1.001$.
} \label{fig2}
\end{figure}
\begin{figure}
\centering
\includegraphics[width=2.1 in,height=2in , angle=0]{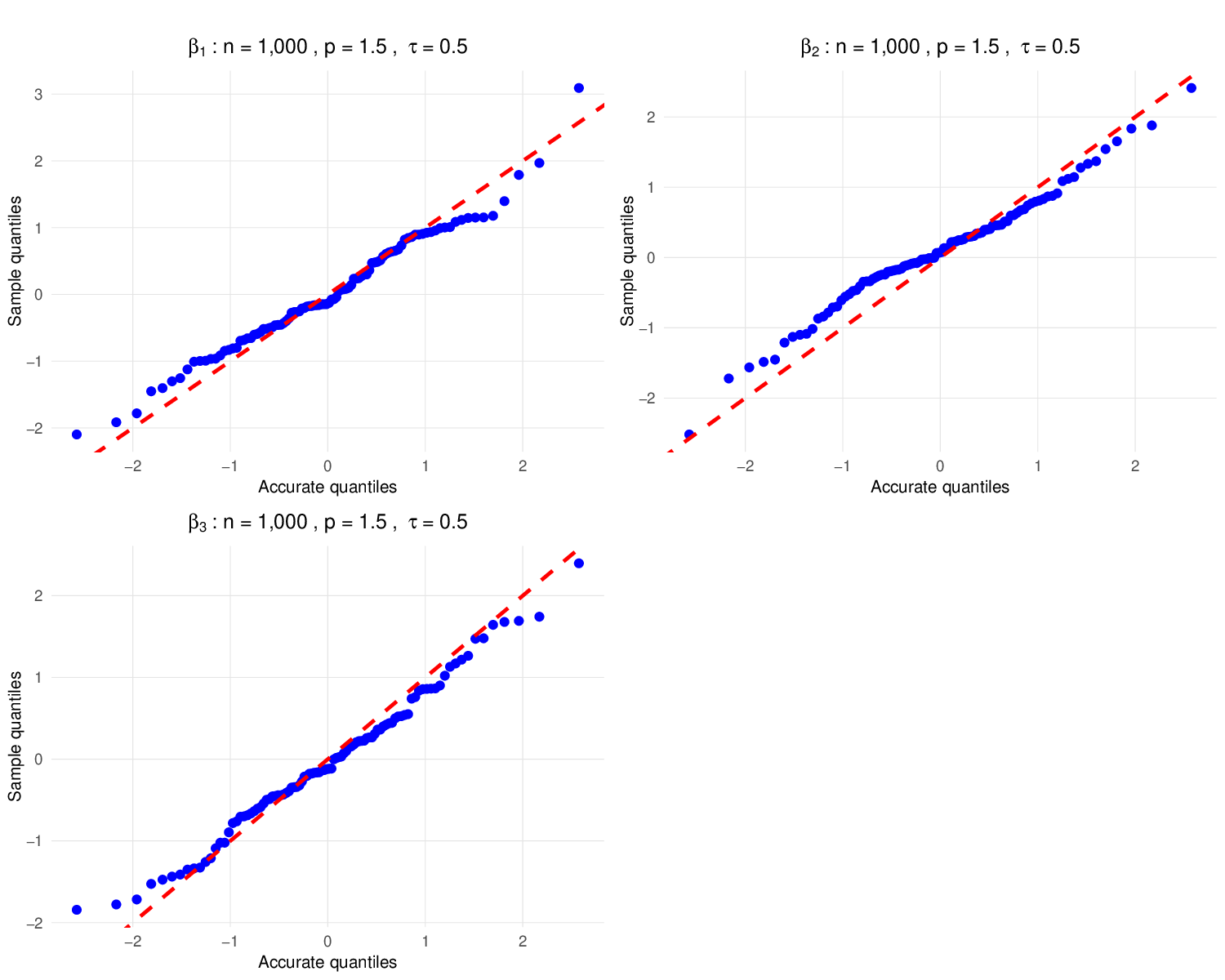}
\includegraphics[width=2.1 in,height=2in , angle=0]{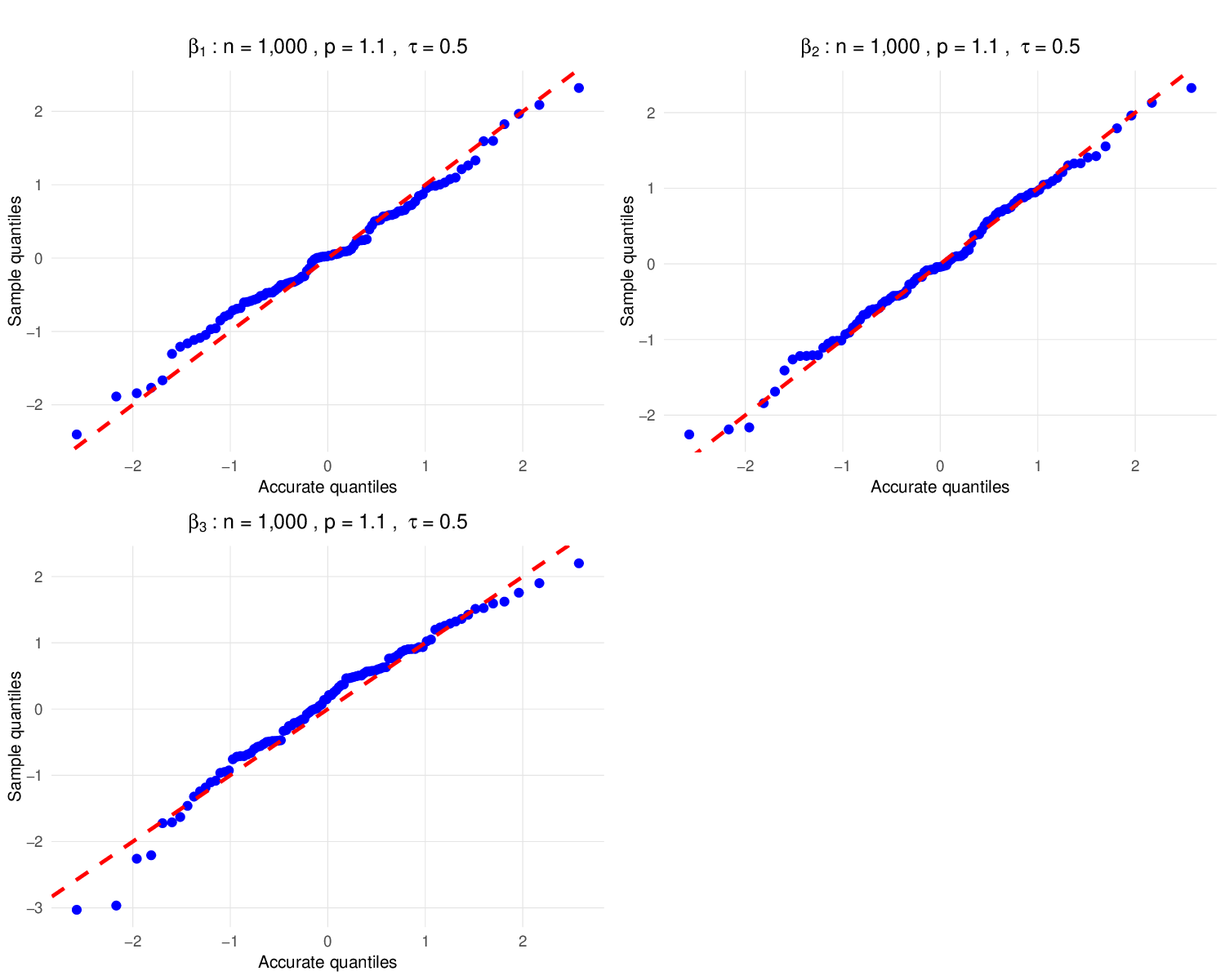}
\includegraphics[width=2.1 in,height=2in , angle=0]{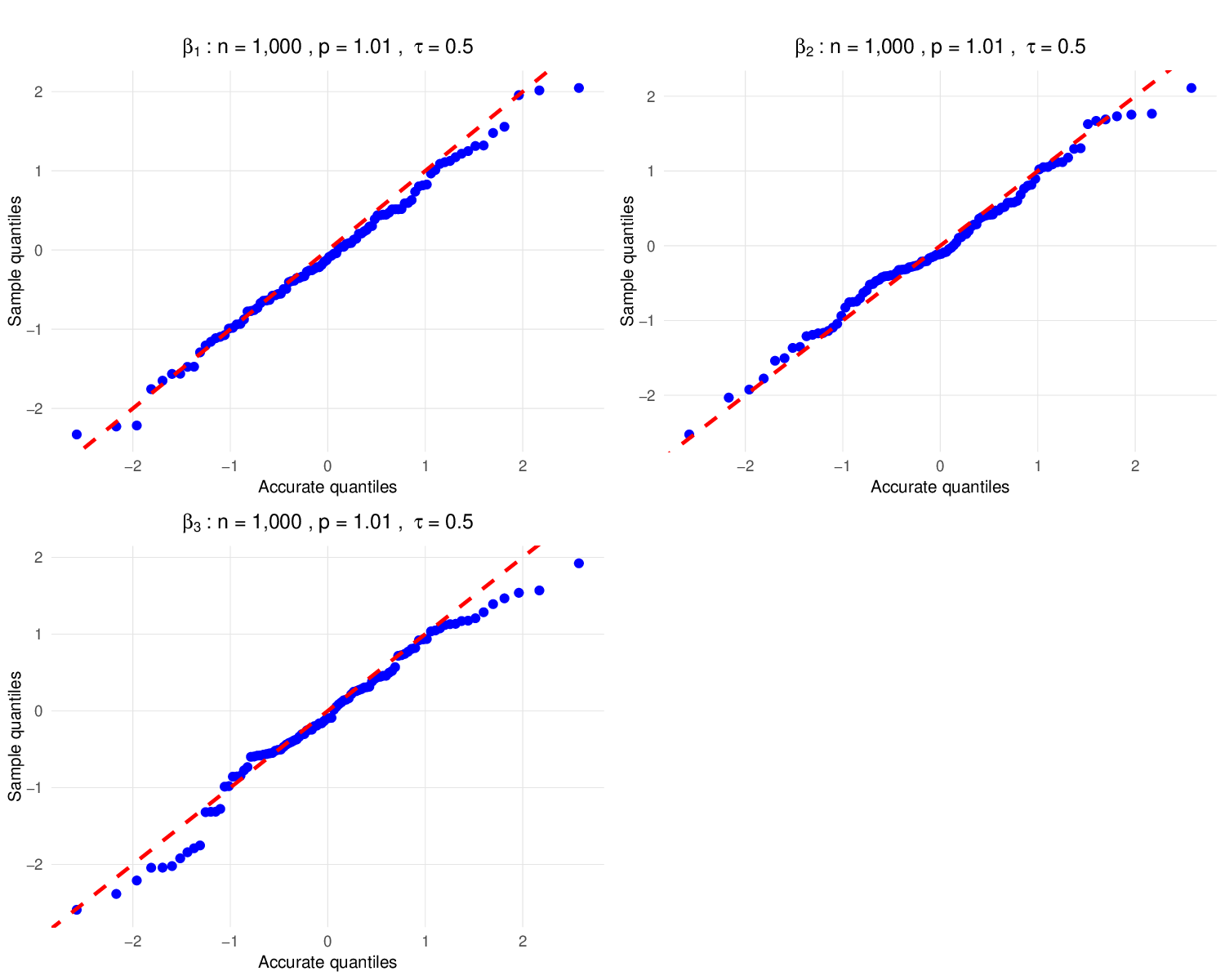}
\includegraphics[width=2.1 in,height=2in , angle=0]{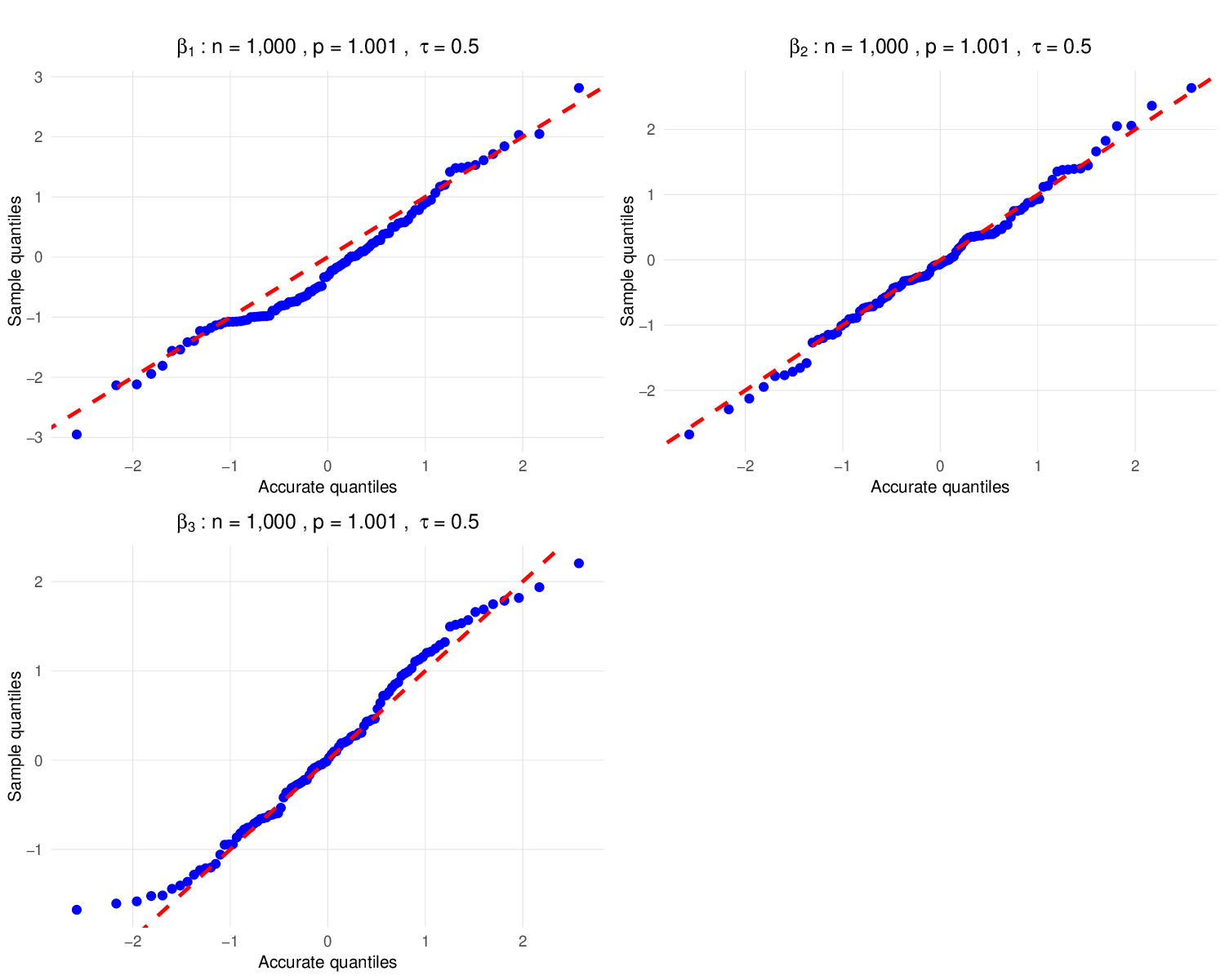}
\caption{The Normal error, the case of $T=1000$, Upper left panel: $p=1.5$. Upper right panel: $p=1.1$. Lower left panel: $p=1.01$. Lower right panel: $p=1.001$.
} \label{fig3}
\end{figure}
\begin{figure}
\centering
\includegraphics[width=2.1 in,height=2in , angle=0]{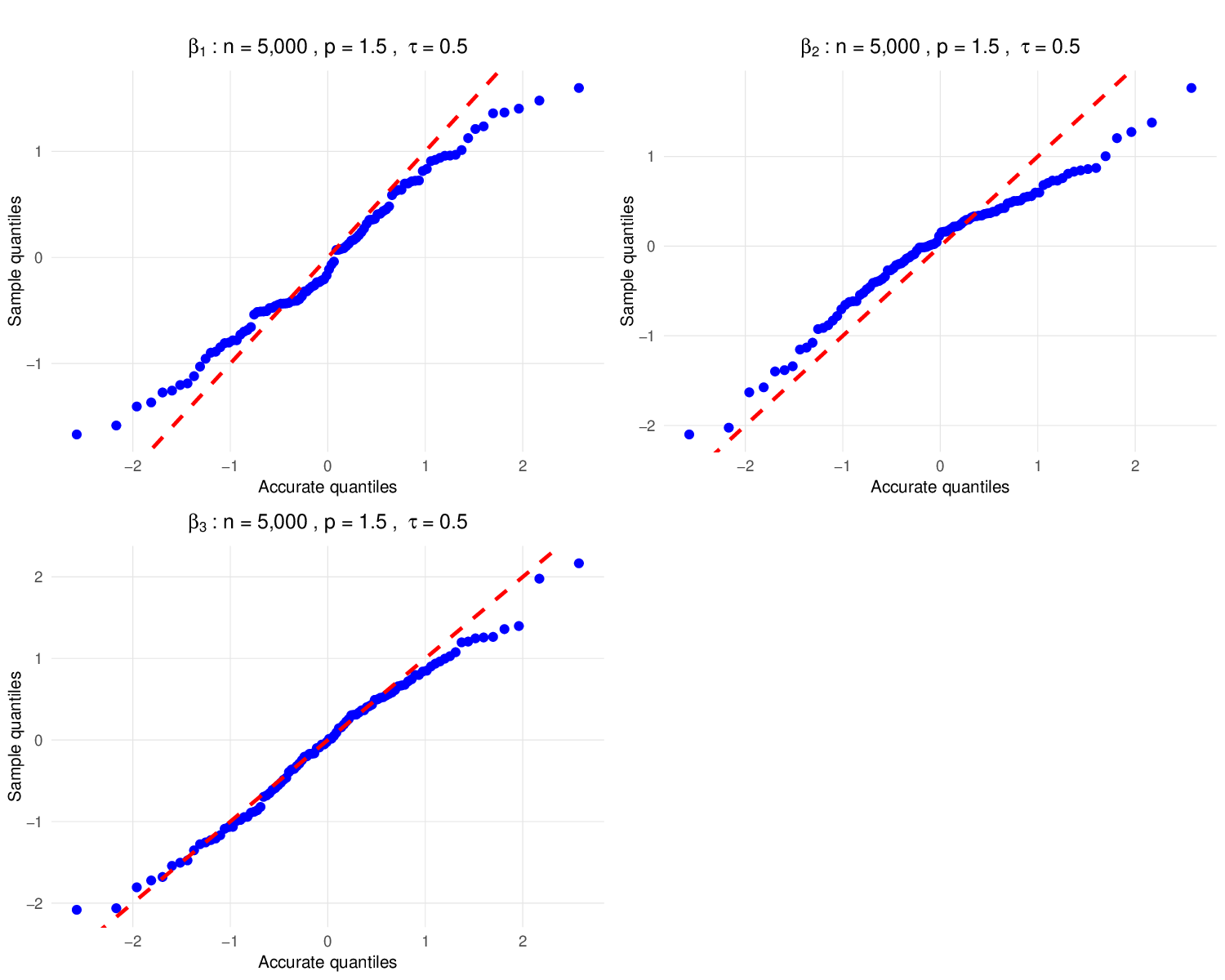}
\includegraphics[width=2.1 in,height=2in , angle=0]{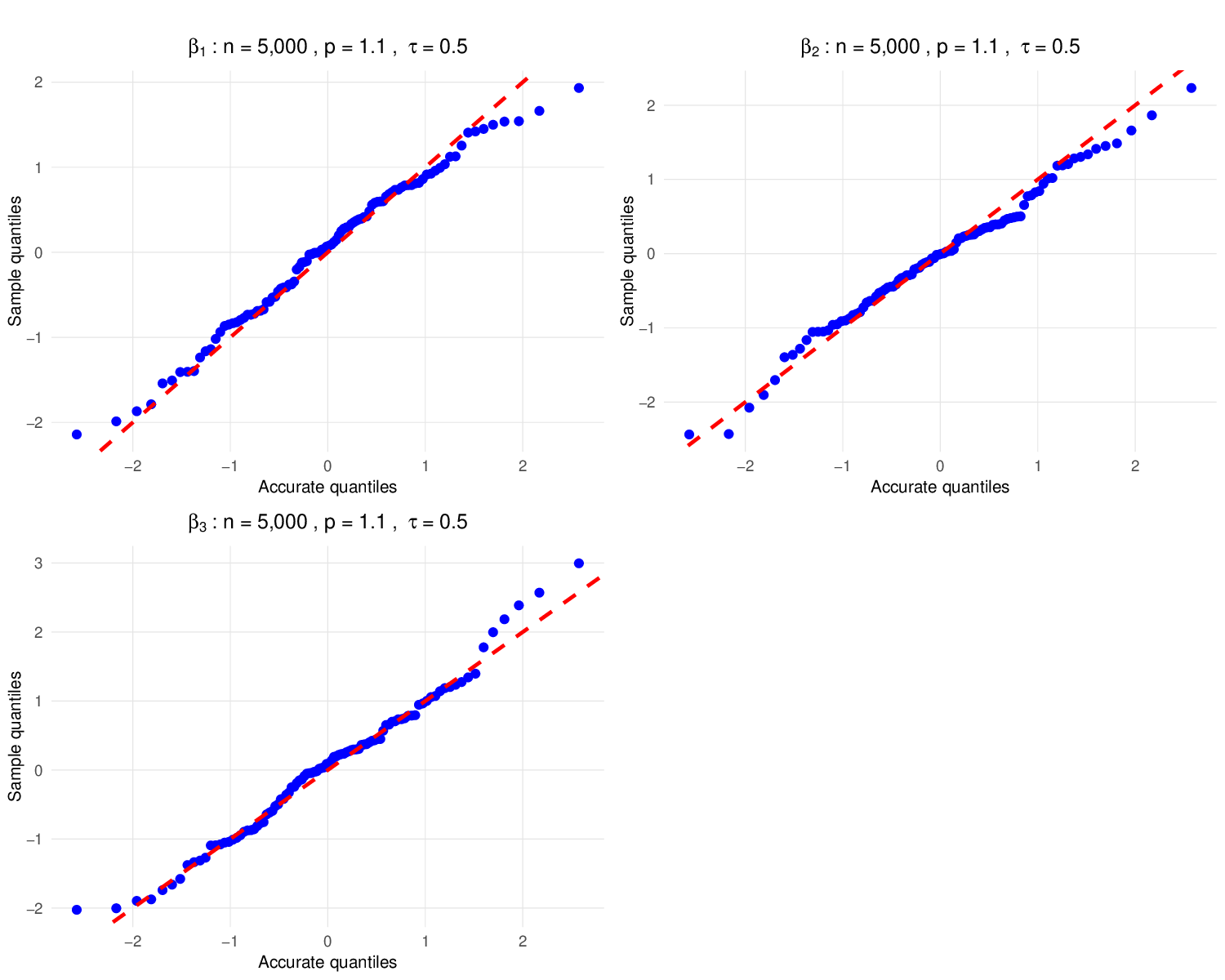}
\includegraphics[width=2.1 in,height=2in , angle=0]{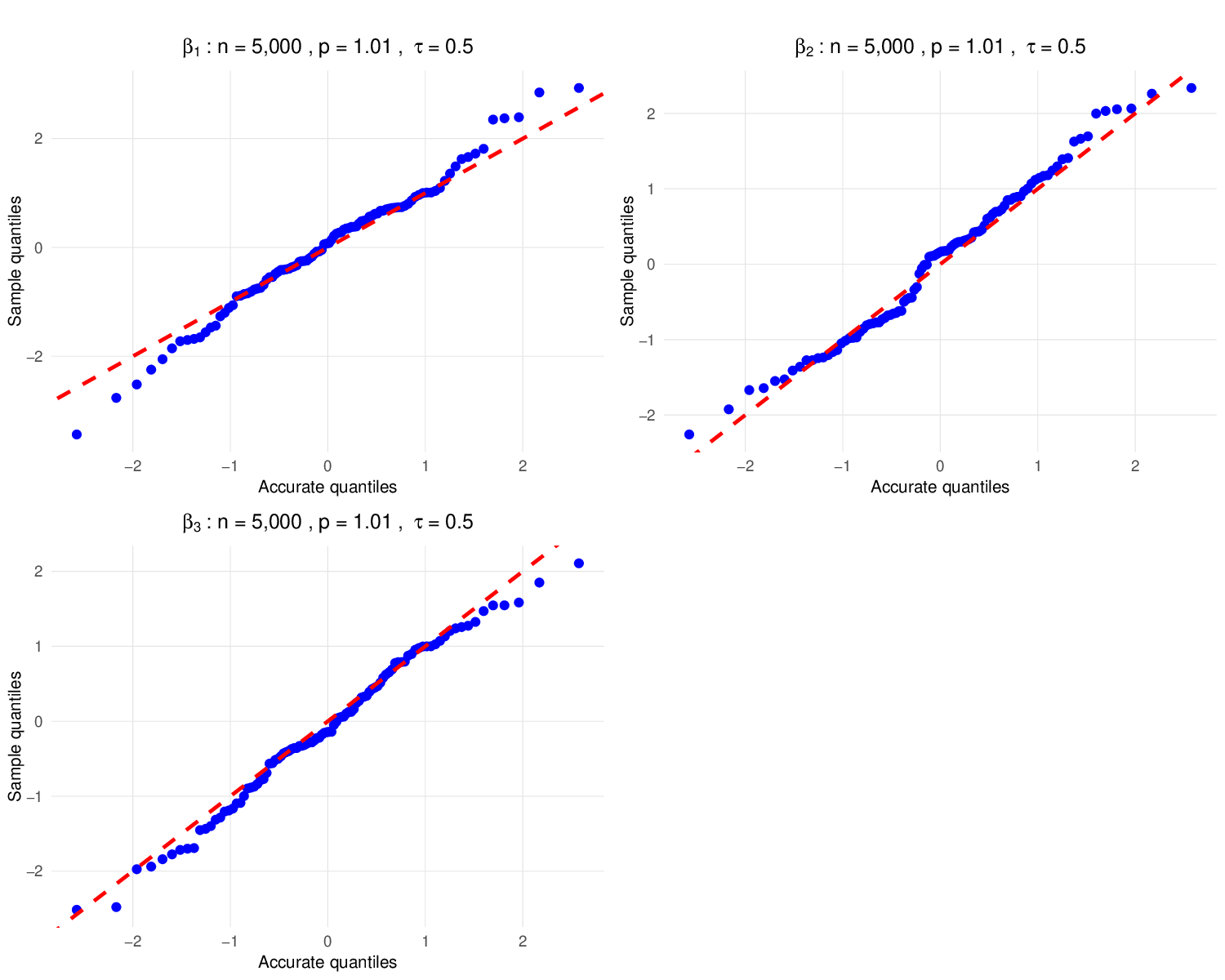}
\includegraphics[width=2.1 in,height=2in , angle=0]{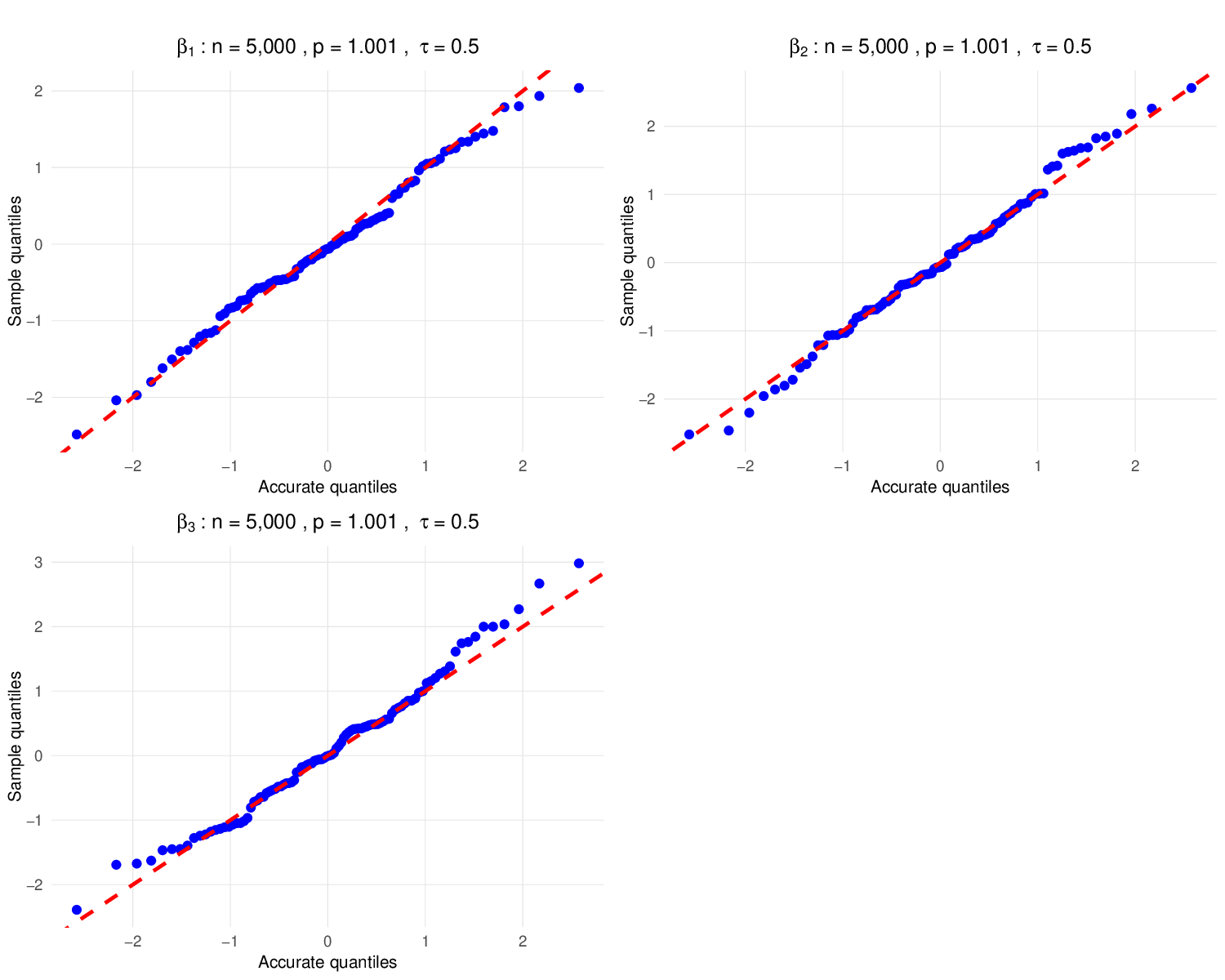}
\caption{The Normal error, the case of $T=5000$, Upper left panel: $p=1.5$. Upper right panel: $p=1.1$. Lower left panel: $p=1.01$. Lower right panel: $p=1.001$.
} \label{fig4}
\end{figure}
\begin{figure}
\centering
\includegraphics[width=2.1 in,height=2in , angle=0]{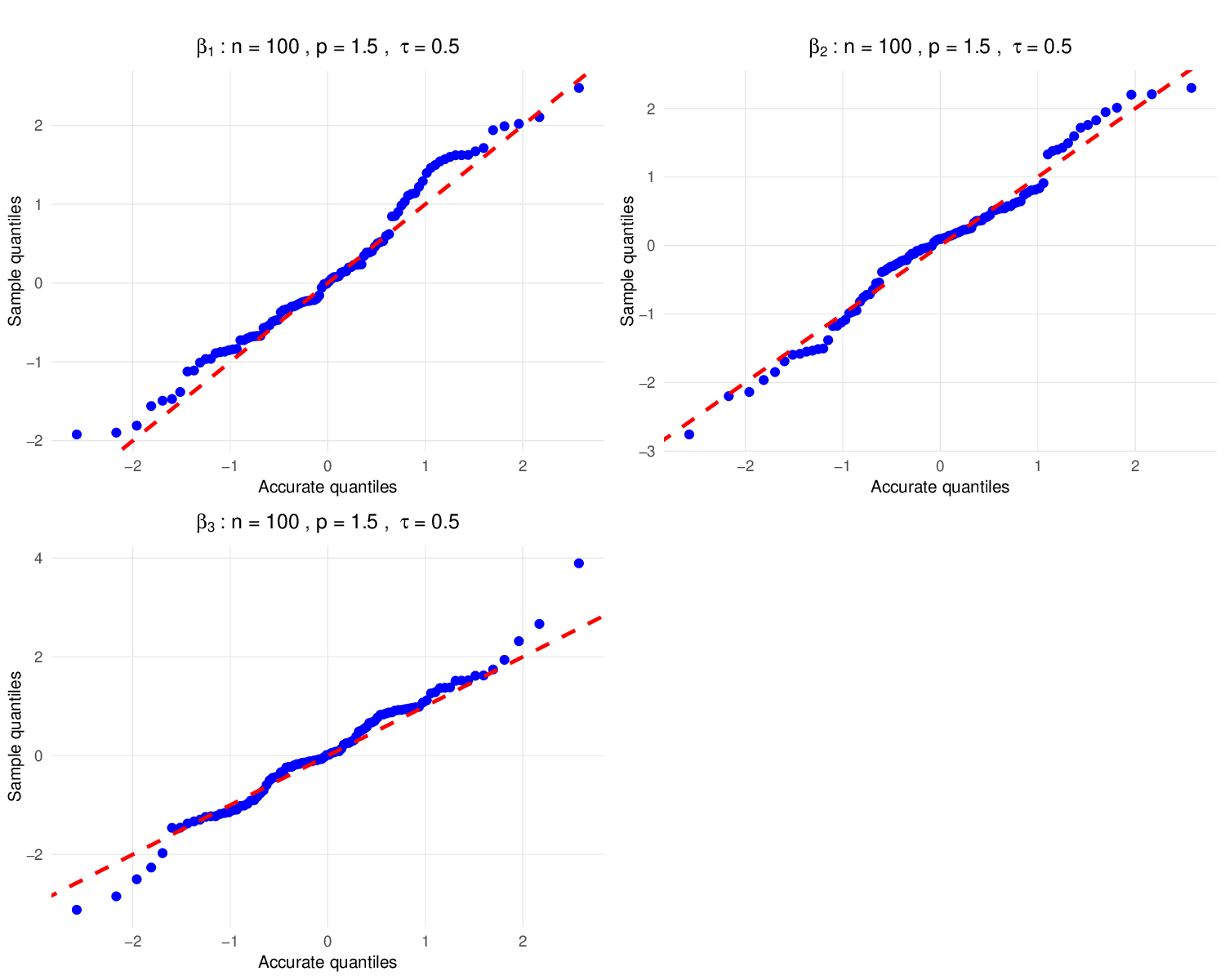}
\includegraphics[width=2.1 in,height=2in , angle=0]{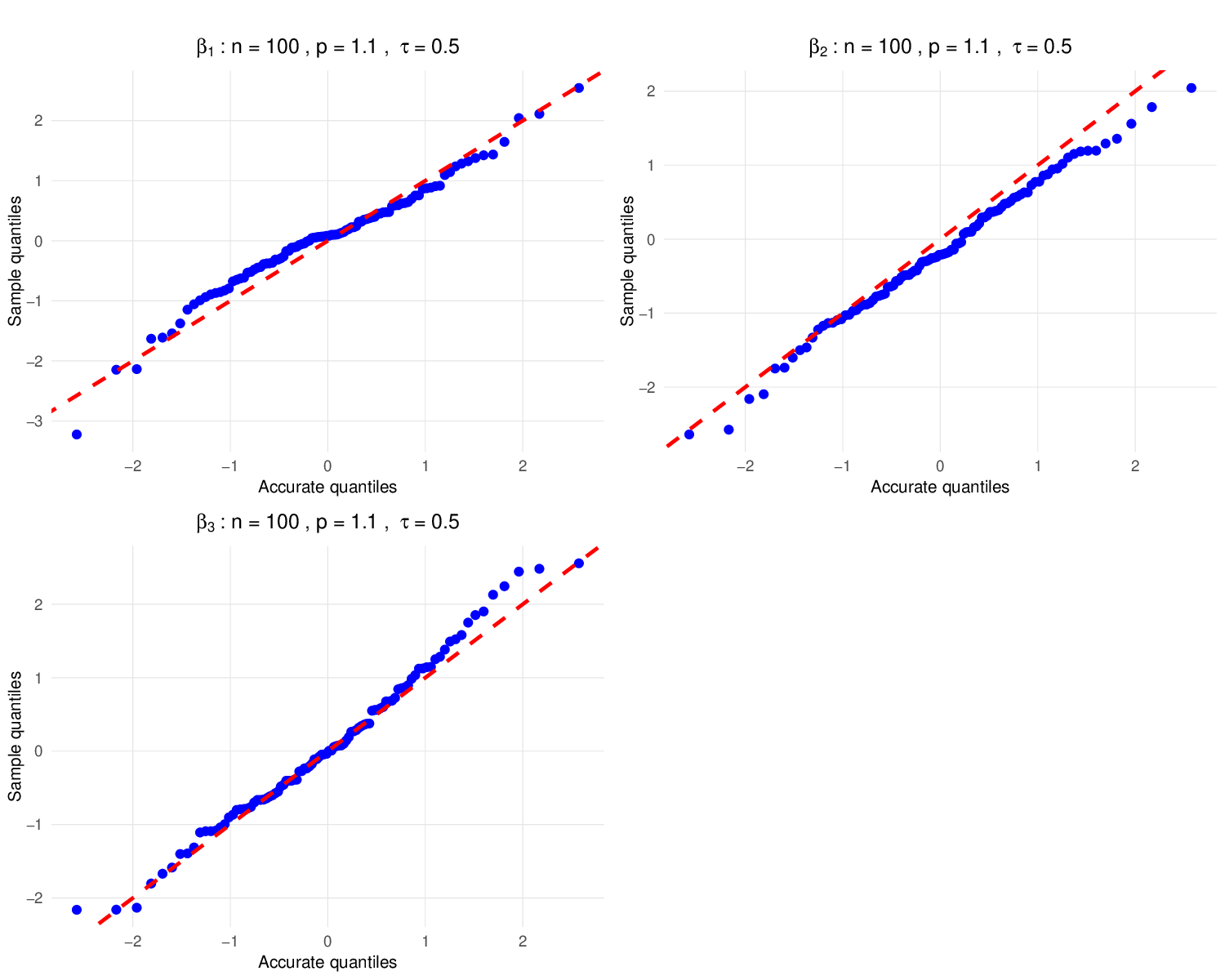}
\includegraphics[width=2.1 in,height=2in , angle=0]{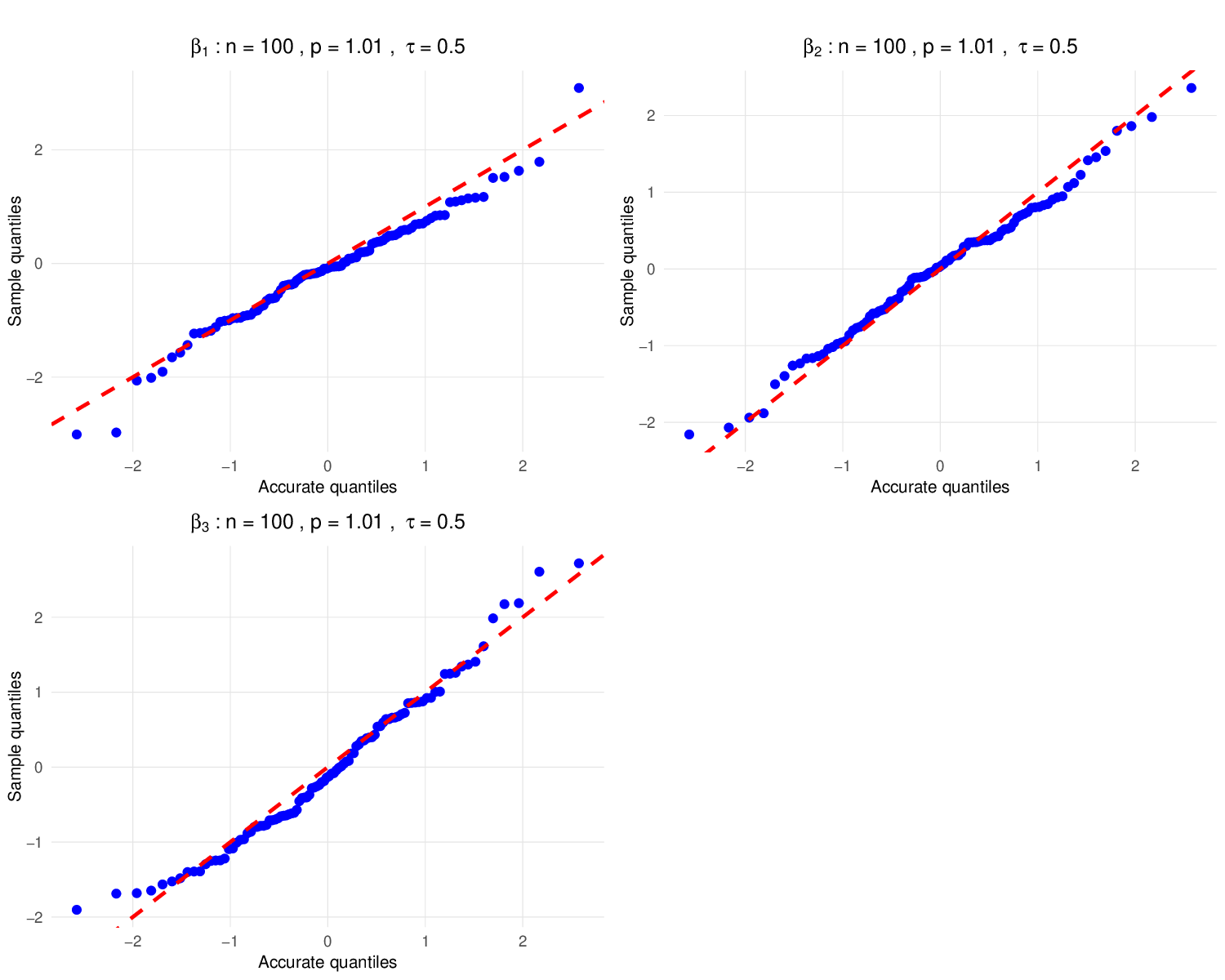}
\includegraphics[width=2.1 in,height=2in , angle=0]{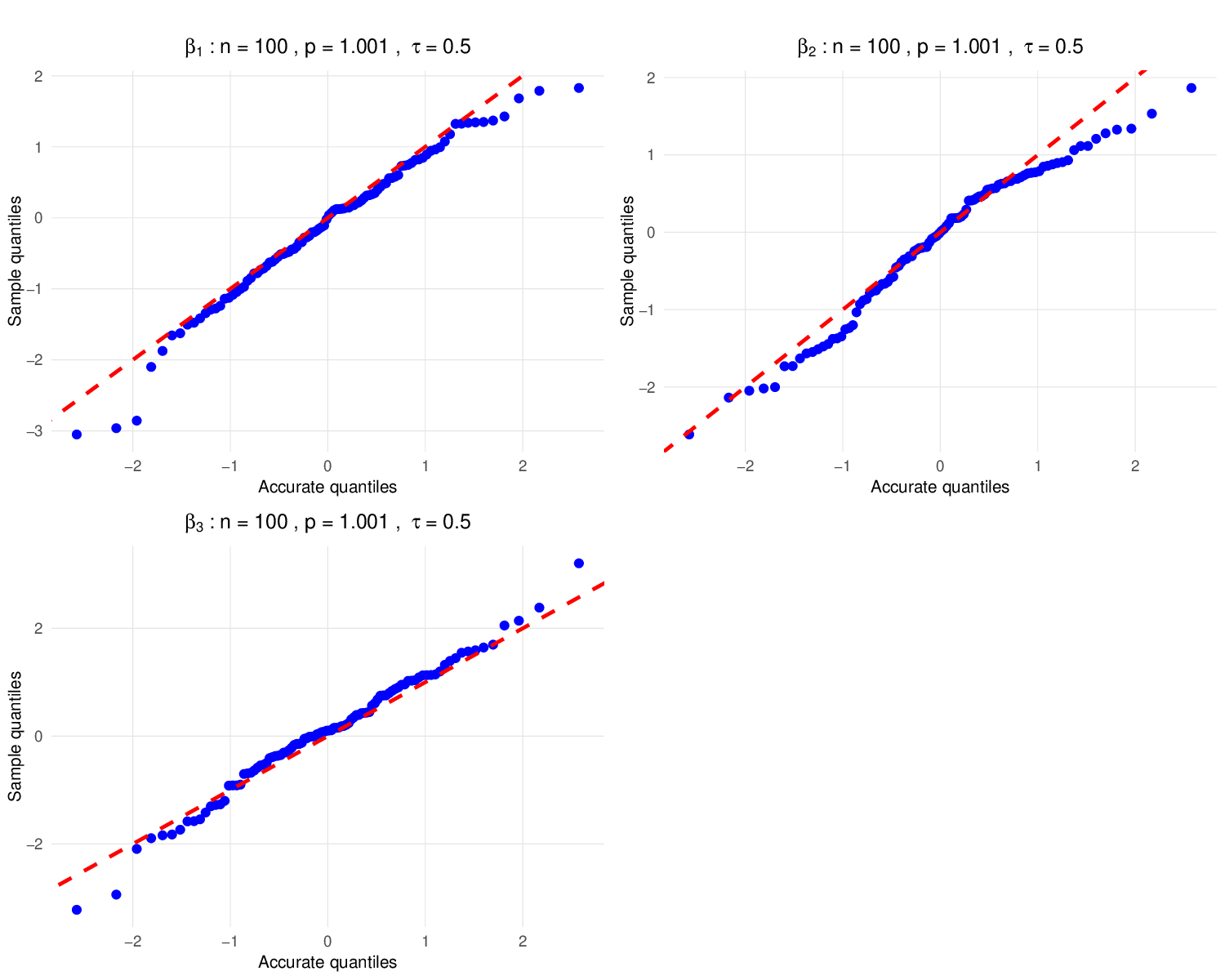}
\caption{The T-distribution error, the case of $T=100$, Upper left panel: $p=1.5$. Upper right panel: $p=1.1$. Lower left panel: $p=1.01$. Lower right panel: $p=1.001$.
} \label{fig5}
\end{figure}
\begin{figure}
\centering
\includegraphics[width=2.1 in,height=2in , angle=0]{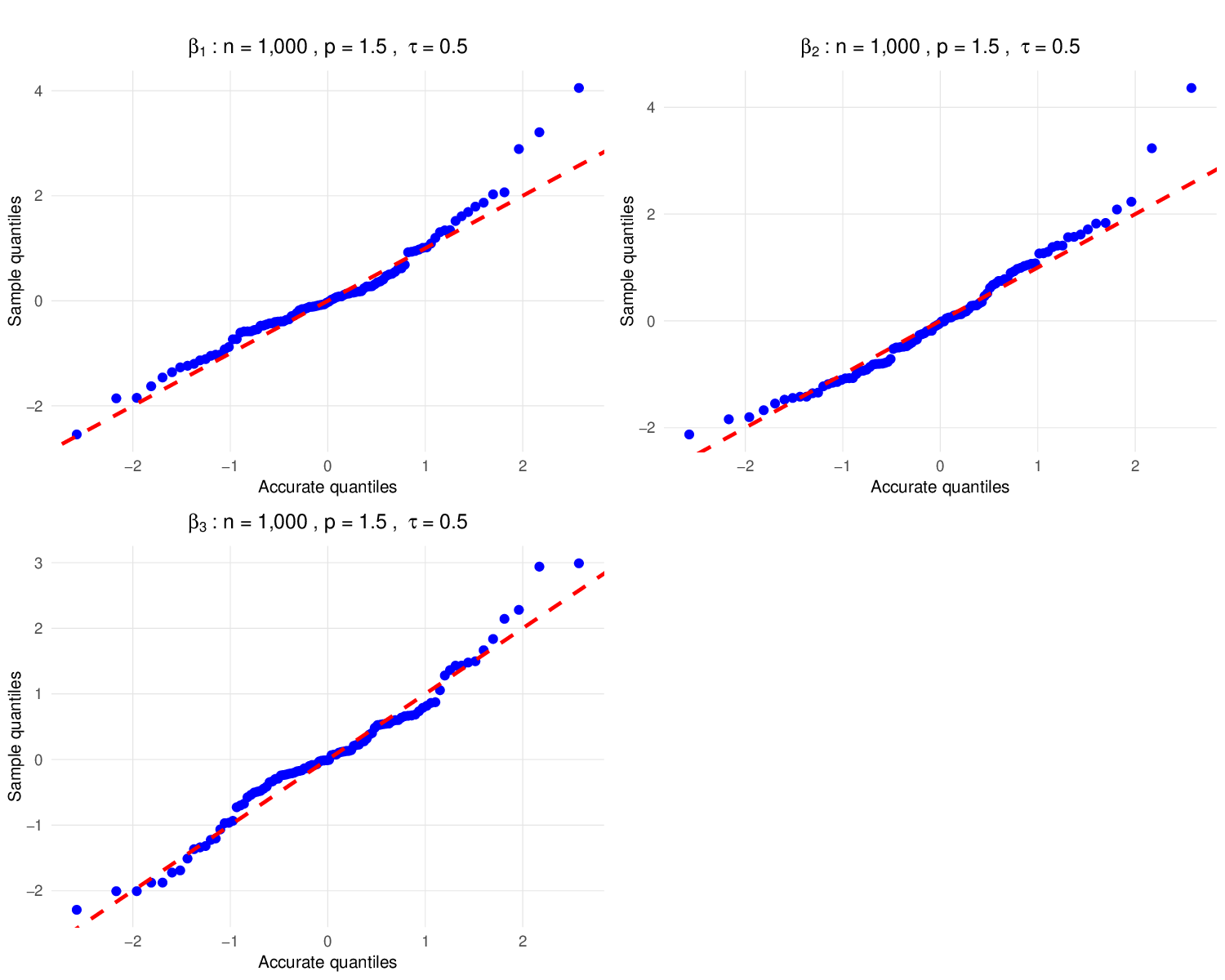}
\includegraphics[width=2.1 in,height=2in , angle=0]{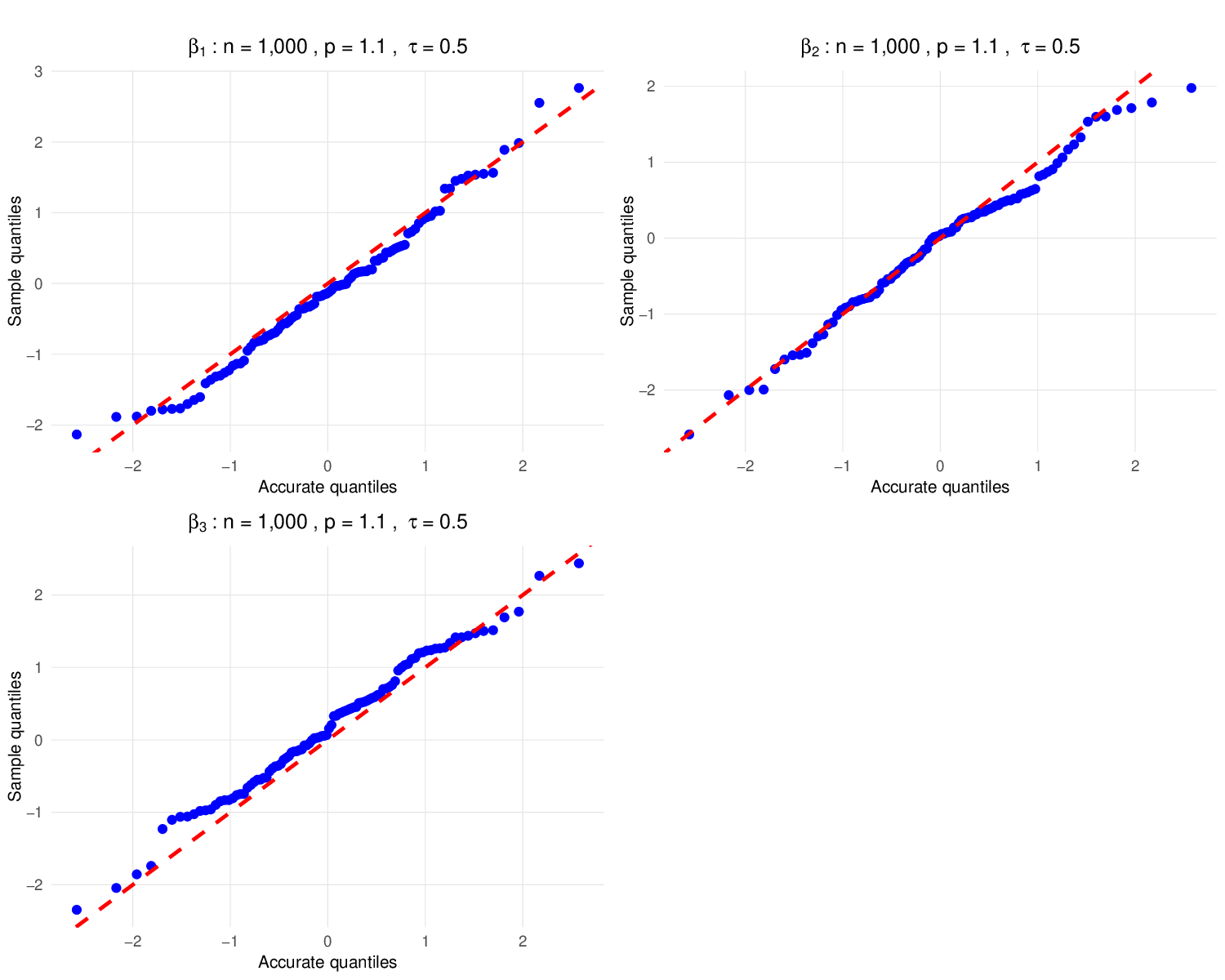}
\includegraphics[width=2.1 in,height=2in , angle=0]{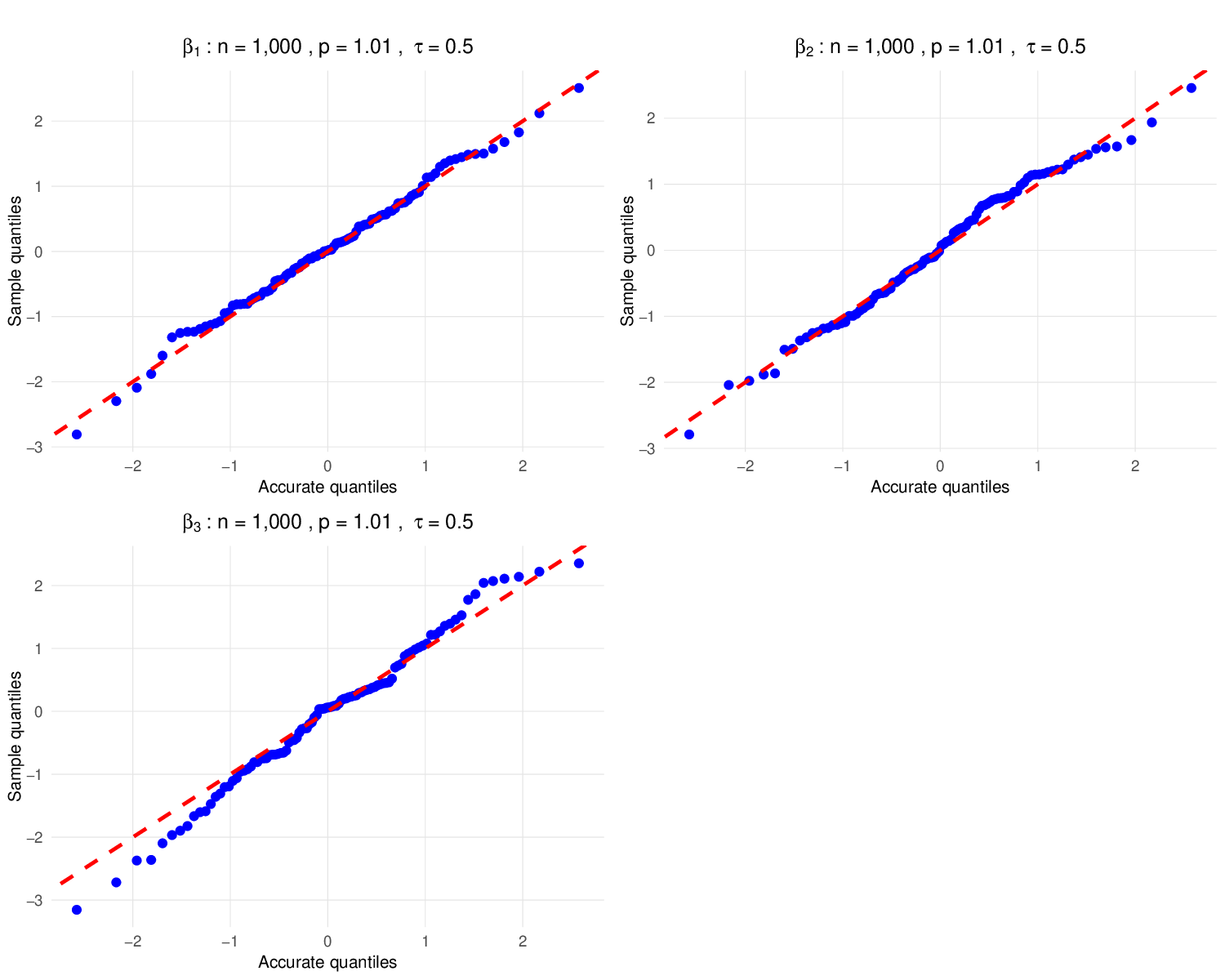}
\includegraphics[width=2.1 in,height=2in , angle=0]{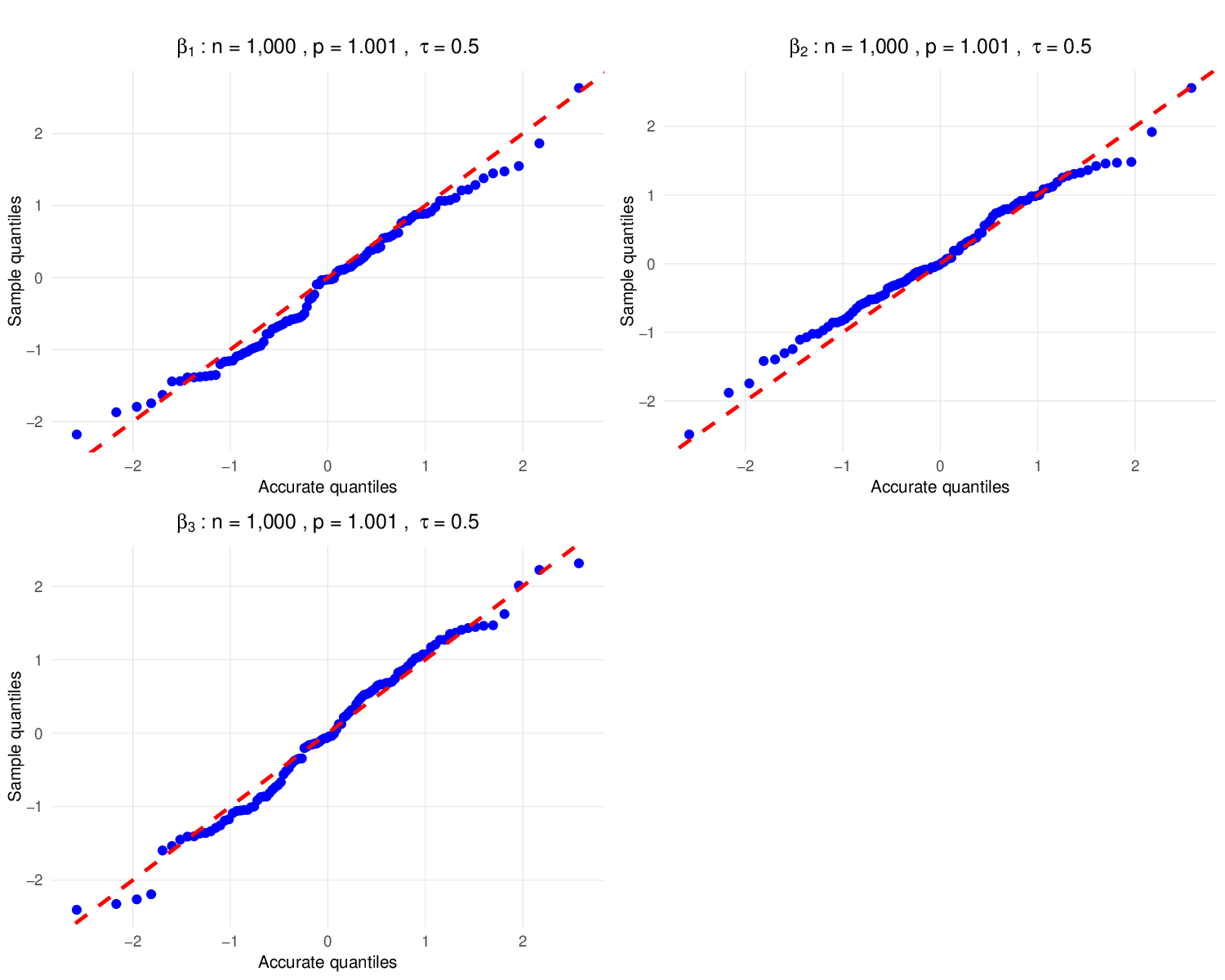}
\caption{The T-distribution error, the case of $T=1000$, Upper left panel: $p=1.5$. Upper right panel: $p=1.1$. Lower left panel: $p=1.01$. Lower right panel: $p=1.001$.
} \label{fig6}
\end{figure}
\begin{figure}
\centering
\includegraphics[width=2.1 in,height=2in , angle=0]{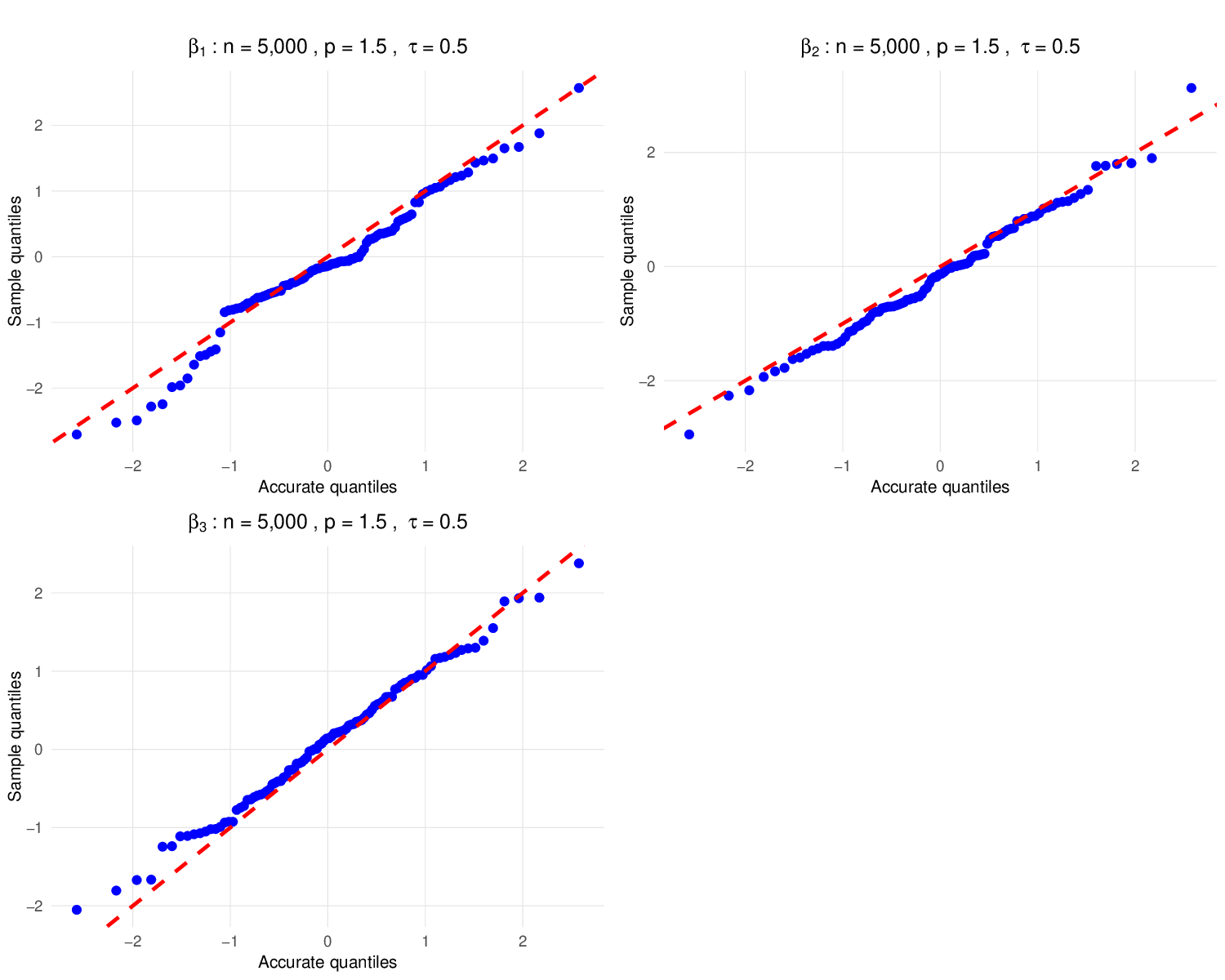}
\includegraphics[width=2.1 in,height=2in , angle=0]{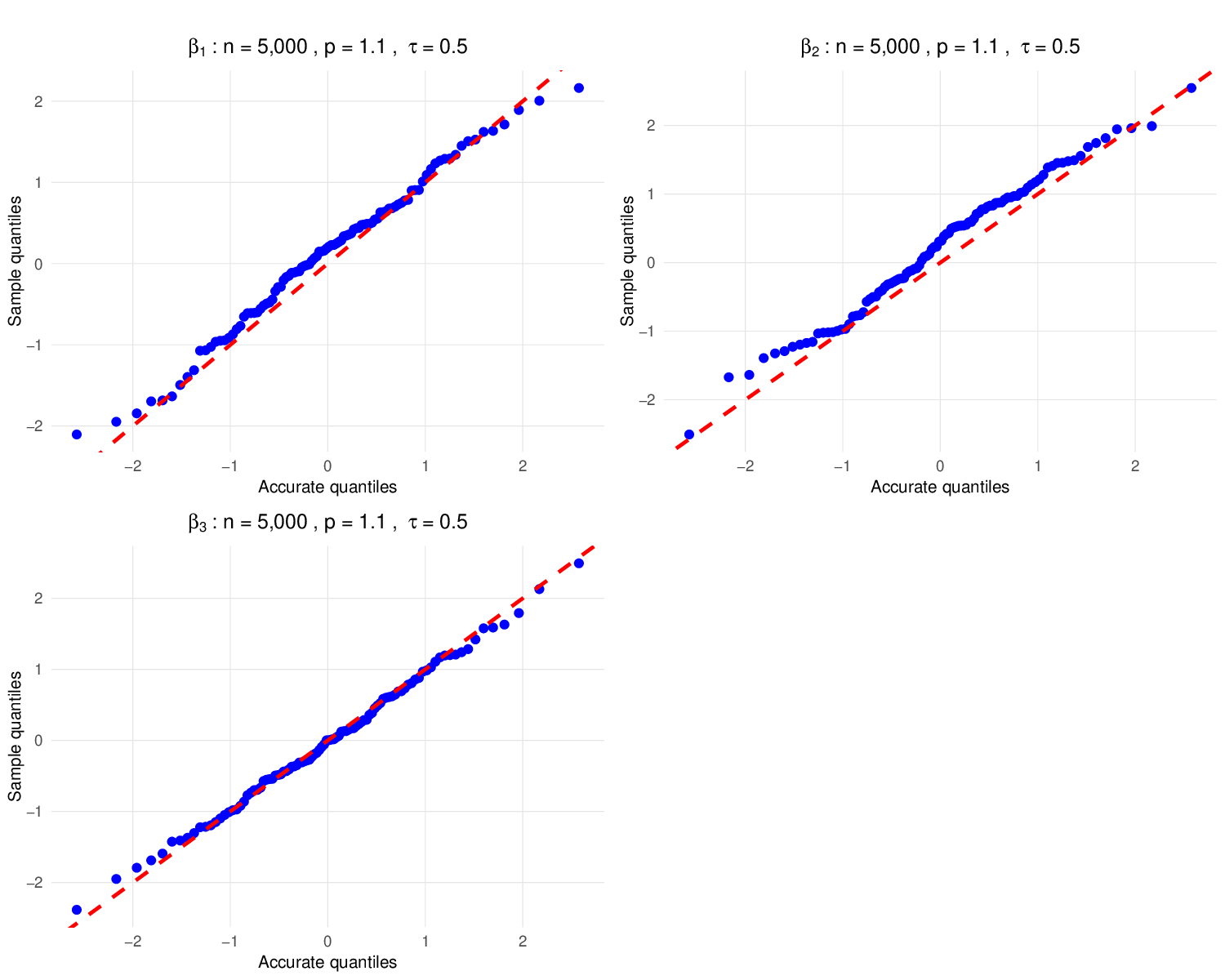}
\includegraphics[width=2.1 in,height=2in , angle=0]{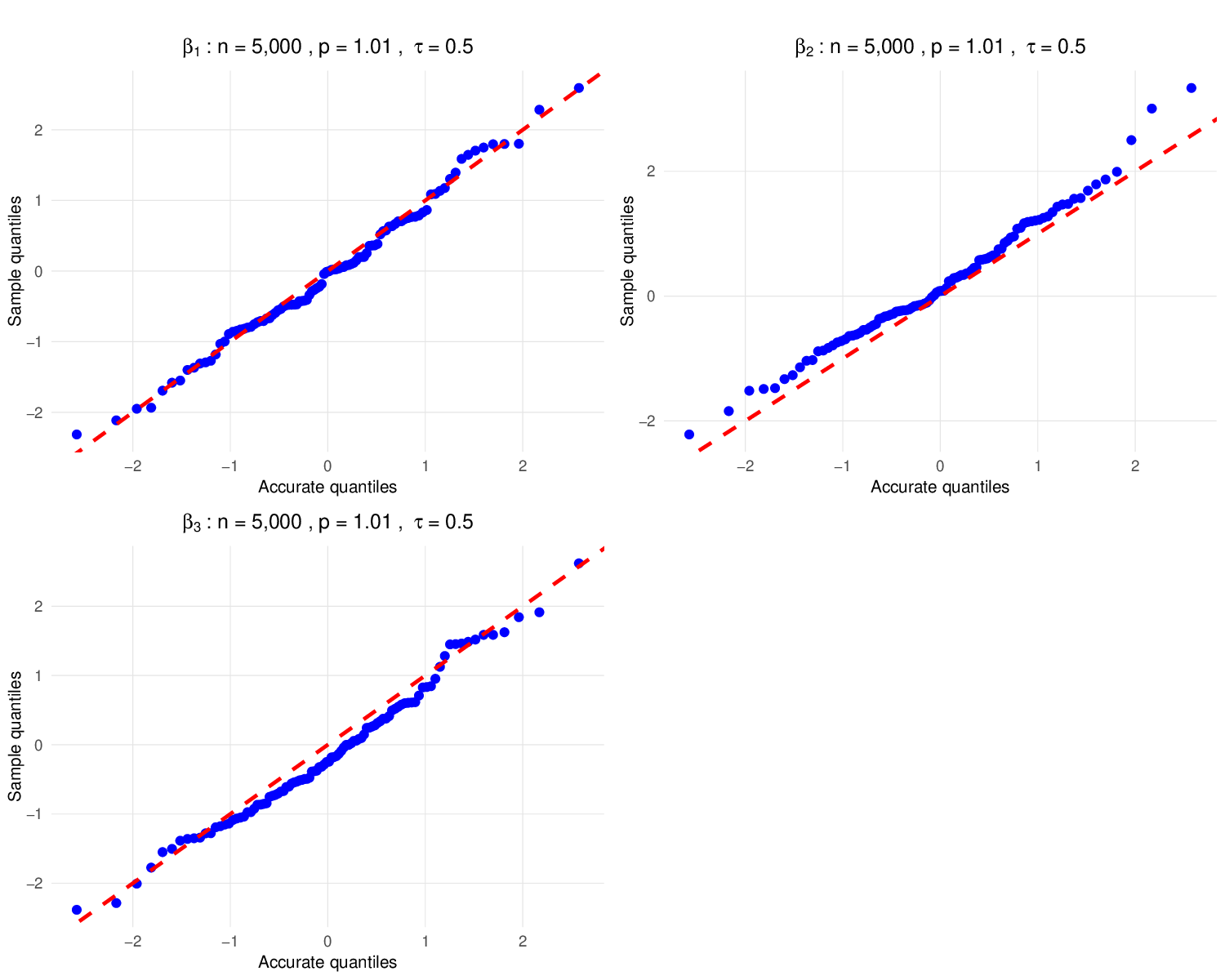}
\includegraphics[width=2.1 in,height=2in , angle=0]{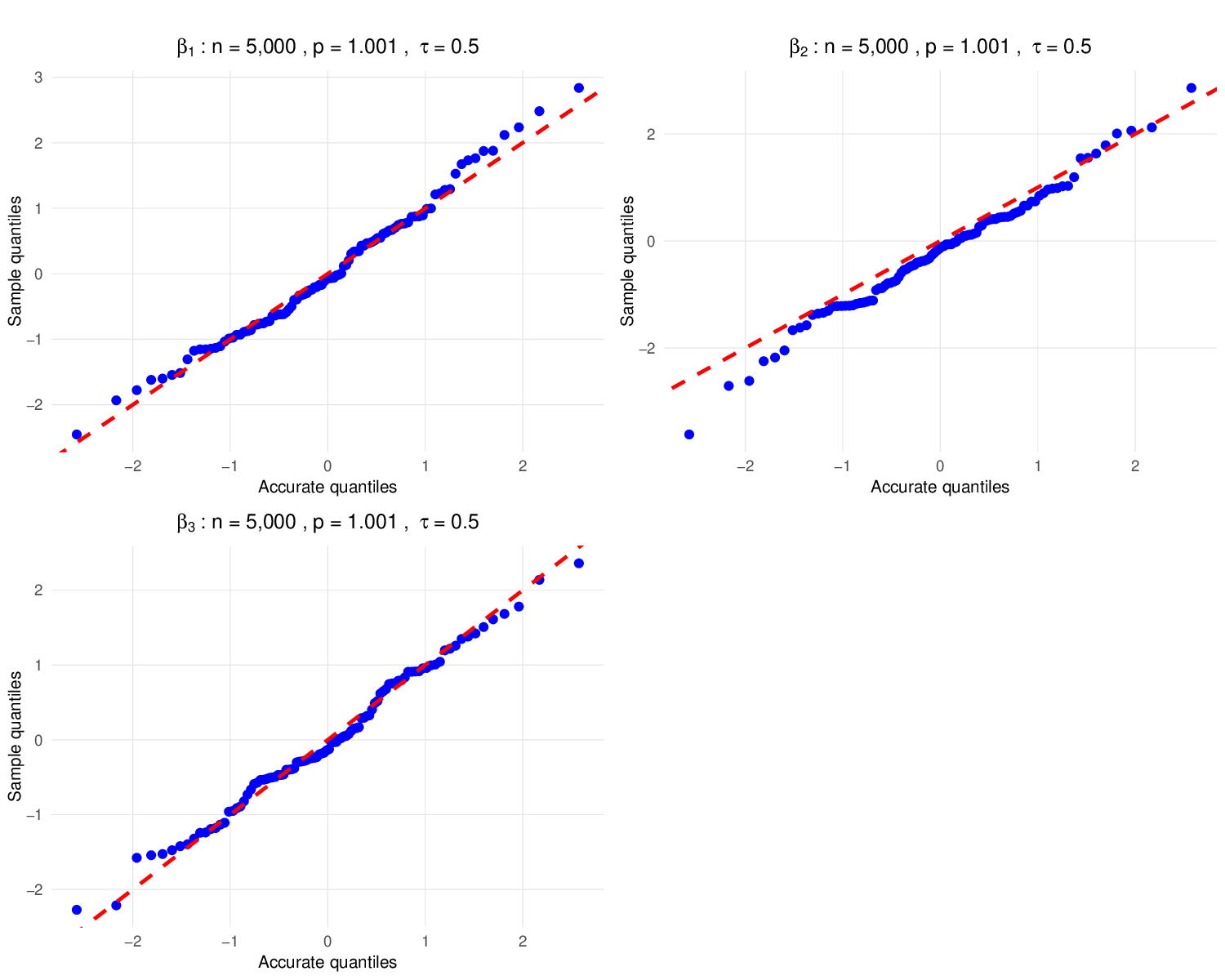}
\caption{The T-distribution error, the case of $T=5000$, Upper left panel: $p=1.5$. Upper right panel: $p=1.1$. Lower left panel: $p=1.01$. Lower right panel: $p=1.001$.
} \label{fig7}
\end{figure}

\section{A real example}
In this section, we apply the proposed method and algorithm to the housing market data in Harrison and Rubinfeld (1978)\cite{Harrison}. We use the augmented and corrected version of it, which is available online at http://lib.stat.cmu.edu/datasets/boston.
The data includes 506 observations, corrected median value of owner-occupied homes (CMEDV) as one response variable, and 15 non-constant predictor variables. They are longitude (LON), latitude (LAT), crime rate (CRIM), proportion of area zoned with large
lots (ZN), proportion of non-retail business acres per town (INDUS), Charles River as a dummy variable (=1 if tract bounds river; 0 otherwise) (CHAS), nitric oxides concentration (NOX), average number of rooms per dwelling (RM), proportion of owner-occupied units built prior to 1940 (AGE), weighted distances to five Boston employment centres (DIS), index of accessibility to radial high-
ways (RAD), property tax rate (TAX), pupil-teacher ratio by town (PTRATIO), black population proportion town (B), and lower status population proportion (LSTAT). Similar to the setting in Wu and Liu (2009)\cite{Wu}, we drop the categorical variable RAD and standardize the response variable and predictor variables except CHAS. Ultimately, we consider the standardized CMEDV as the response and the variable CHAS, the 13 standardized predictor variables and their squares as predictors (27 variables). We apply the CLpQR with the adaptive penalty to the latest data with $p$ taking value in the set $\{1, 1.1, 1.3, 1.5, 1.8, 2, 2.1\}$.

In order to compare estimation error and variable selection results for various $p$ cases we run the regression 10 times in each case. In each repetition, the data is randomly split into the training, tuning and testing data sets with size 200, 150, and 156.
We select the tuning parameter by minimizing the objective function in \eqref{eq1.4} and separately use $L^{p}$, $L^{1}$, and $L^{2}$ distance to calculate the test error. The later two distance are more important as only under the same distance, we can accurately choose suitable $p$ for the CLpQR. Empirical results are summarized in Table 7.1. Results show that $p=1.3$ is a good choice if one cares more about the stability of variable selection; $p=2$ or so is desirable if one concerns more the average precision. Moreover there is a difference between $L^{1}$ and $L^{2}$ distance when using them to calculate the deviation of estimation error: The former generates the smallest deviation when $p=1.3$, while the latter does that when $p=2.1$.
\begin{table}
\begin{center}
Table 7.1. Empirical results for the Corrected Boston House Price Data
\begin{tabular}{cccccccccccccccccccccccccccccccccccccccccccc}
\hline
& &$L^{p}$ distance  & $L^{1}$ distance   & $L^{2}$ distance  \\
p &no. of zeros & test error & test error & 	test error \\
\hline
1 &12.5 (2.418)& 0.3672 (0.0251) & 0.3672 (0.0251)  & 0.2832 (0.0636) \\
1.1&11.8	(2.481)& 0.3422	(0.0256)& 0.3637	(0.0240) & 0.2779	(0.0610) \\
1.3&10.1 (1.813)& 0.3124	(0.0245)& 0.3677	(0.0171) & 0.2747	(0.0568) \\
1.5&10.2 (2.227)& 0.2891	(0.0292) & 0.3691 (0.0187)	& 0.2688 (0.0512)	 \\
1.8&11.8	(5.344)& 0.2641	(0.0410) & 0.3655	(0.0285)& 0.2595	(0.0480)\\
2&12.3	(4.712)& 0.2492	(0.0449) & 0.3635	(0.0272) & 	0.2492	(0.0449) \\
2.1&13.2	(5.400)& 0.2395	(0.0445) & 0.3607	(0.0253)  & 	0.2415	(0.0410)\\
\hline
\end{tabular}
\end{center}
\end{table}
\section{Proofs}
{\bf Proof of Theorem \ref{thm2.1}.} Let $\sqrt{T}(\hat{\boldsymbol{\beta}}^{clp}-\boldsymbol{\beta}^{*})=\textbf{u}_{T}$ and
$\sqrt{T}(\hat{b}_{k}-b^{*}_{\tau_{k}})=u_{T, k}$. Then $(u_{T, 1},\cdots, u_{T, K}, \textbf{u}_{T})$
is the minimizer of the following criterion function
\begin{eqnarray*}\label{eq1.7}
Q_{T}=\sum_{k=1}^{K}\sum^{T}_{t=1}
\Big[\boldsymbol{\eta}_{\tau_{k},p}\Big(\varepsilon_{t}-b_{\tau_{k}}^{*}-\frac{u_{k}+\textbf{x}'_{t}\textbf{u}}{\sqrt{T}}\Big)
-\boldsymbol{\eta}_{\tau_{k},p}(\varepsilon_{t}-b_{\tau_{k}}^{*})\Big]
\end{eqnarray*}
over $u_{1}, \cdots, u_{K}, u$. Write $Q_{T}$ as
\begin{eqnarray*}\label{eq1.8}
\nonumber Q_{T}&=&\sum_{k=1}^{K}\sum^{T}_{t=1}
\Bigg[-\frac{u_{k}+\textbf{x}'_{t}\textbf{u}}{\sqrt{T}}\boldsymbol{\varphi}_{\tau_{k},p}(\varepsilon_{t}-b_{\tau_{k}}^{*})
\\
& &-\int_{0}^{(u_{k}+\textbf{x}'_{t}\textbf{u})/\sqrt{T}}(\boldsymbol{\varphi}_{\tau_{k},p}(\varepsilon_{t}-b_{\tau_{k}}^{*}-t)
-\boldsymbol{\varphi}_{\tau_{k},p}(\varepsilon_{t}-b_{\tau_{k}}^{*}))dt\Bigg].
\end{eqnarray*}
Define $Z_{T, k}=\frac{1}{\sqrt{T}}\sum^{T}_{t=1}\boldsymbol{\varphi}_{\tau_{k},p}(\varepsilon_{t}-b_{\tau_{k}}^{*})$,
$\textbf{Z}_{T}=\frac{1}{\sqrt{T}}\sum^{T}_{t=1}\textbf{x}'_{t}[\sum_{k=1}^{K}\boldsymbol{\varphi}_{\tau_{k},p}(\varepsilon_{t}
-b_{\tau_{k}}^{*})]$, $B_{T, k}=\sum^{T}_{t=1}\int_{0}^{(u_{k}+\textbf{x}'_{t}\textbf{u})/\sqrt{T}}(\boldsymbol{\varphi}_{\tau_{k},p}(\varepsilon_{t}-b_{\tau_{k}}^{*}-t)
-\boldsymbol{\varphi}_{\tau_{k},p}(\varepsilon_{t}-b_{\tau_{k}}^{*}))dt$. So we have
\begin{eqnarray*}\label{eq1.9}
Q_{T}=-\sum^{K}_{k=1}Z_{T, k}u_{k}-\textbf{Z}'_{T}\textbf{u}-\sum^{K}_{k=1}B_{T, k}.
\end{eqnarray*}
Under Assumption \ref{ass2.2}, using the Cram\'{e}r-Wald method and CLT, we get
\begin{eqnarray*}\label{eq1.10}
(Z_{T, 1}, \cdots, Z_{T, K}, \textbf{Z}'_{T})'\stackrel{\textit{D}}{\longrightarrow}
(Z_{1}, \cdots, Z_{K}, \textbf{Z}')'\sim N(0, \boldsymbol{\Sigma}),
\end{eqnarray*}
where the asymptotic covariance matrix $\boldsymbol{\Sigma}$ can be easily gotten by the routine procedure.

Next, focus on the limit of $B_{T, k}$. By some calculation, we have
\begin{eqnarray}\label{eq1.11}
\nonumber E(B_{T, k})&=&\frac{1}{T}\sum_{t=1}^{T}\int^{u_{k}+\textbf{x}'_{t}\textbf{u}}_{0}\sqrt{T}E
(\boldsymbol{\varphi}_{\tau_{k},p}(\varepsilon_{t}-b_{\tau_{k}}^{*}-s/\sqrt{T})
-\boldsymbol{\varphi}_{\tau_{k},p}(\varepsilon_{t}-b_{\tau_{k}}^{*}))ds\\
&=&\frac{1}{T}\sum_{t=1}^{T}\int^{u_{k}+\textbf{x}'_{t}\textbf{u}}_{0}E
\boldsymbol{\psi}_{\tau_{k},p}(\varepsilon_{t}-b_{\tau_{k}}^{*}-\tilde{s}/\sqrt{T})(-s)ds,
\end{eqnarray}
where $\tilde{s}$ lies between 0 and $s$. Note that, for positive $\delta>0$, there is a $T$ large enough such that
\begin{eqnarray*}\label{eq1.12}
|\boldsymbol{\psi}_{\tau_{k},p}(\varepsilon_{t}-b_{\tau_{k}}^{*}-\tilde{s}/\sqrt{T})|
\leq p(p-1)||\varepsilon_{t}-b^{*}_{\tau_{k}}|-\delta/2|^{p-2}.
\end{eqnarray*}
For $\delta$ small enough, Assumption \ref{ass2.3} makes sure $E||\varepsilon_{t}-b^{*}_{\tau_{k}}|-\delta/2|^{p-2}<\infty$, which further yields
$E\boldsymbol{\psi}_{\tau_{k},p}(\varepsilon_{t}-b_{\tau_{k}}^{*}-\tilde{s}/\sqrt{T})\rightarrow E\boldsymbol{\psi}_{\tau_{k},p}(\varepsilon_{t}-b_{\tau_{k}}^{*})$.
In fact, the convergence is uniform with respect to $\tilde{s}$ being between 0 and $u_{k}+\textbf{x}'_{t}\textbf{u}$. Namely,  $E\boldsymbol{\psi}_{\tau_{k},p}(\varepsilon_{t}-b_{\tau_{k}}^{*}-\tilde{s}/\sqrt{T})
=E\boldsymbol{\psi}_{\tau_{k},p}(\varepsilon_{t}-b_{\tau_{k}}^{*})(1+o(1))$ uniformly in $\tilde{s}\in [0, u_{k}+\textbf{x}'_{t}\textbf{u}]$ or $\tilde{s}\in [u_{k}+\textbf{x}'_{t}\textbf{u}, 0]$. So, the expression \eqref{eq1.11} equals
\begin{eqnarray*}\label{eq1.13}
\frac{1}{T}\sum_{t=1}^{T}\Big[\int^{u_{k}+\textbf{x}'_{t}\textbf{u}}_{0}E
\boldsymbol{\psi}_{\tau_{k},p}(\varepsilon_{t}-b_{\tau_{k}}^{*})(-s)ds+o(1)\Big].
\end{eqnarray*}
Some calculation induces
\begin{eqnarray*}\label{eq1.13}
E(B_{T, k})\longrightarrow-\frac{1}{2}E
\boldsymbol{\psi}_{\tau_{k},p}(\varepsilon-b_{\tau_{k}}^{*})(u_{k}, \textbf{u}')\left( {{\begin{array}{*{20}c}
1 & 0\\
0 & \textbf{C}\\
\end{array}}}
\right)(u_{k}, \textbf{u}')'.
\end{eqnarray*}
And
\begin{eqnarray}\label{eq1.14}
\nonumber \mbox{Var}(B_{T, k})&\leq&\sum_{t=1}^{T}E\Big[\int_{0}^{(u_{k}+\textbf{x}'_{t}\textbf{u})
/\sqrt{T}}(\boldsymbol{\varphi}_{\tau_{k},p}(\varepsilon_{t}-b_{\tau_{k}}^{*}-t)
-\boldsymbol{\varphi}_{\tau_{k},p}(\varepsilon_{t}-b_{\tau_{k}}^{*}))dt\Big]^{2}\\
\nonumber&\leq &\sum_{t=1}^{T}E\Big[-\int_{0}^{(u_{k}+\textbf{x}'_{t}u)
/\sqrt{T}}(\boldsymbol{\varphi}_{\tau_{k},p}(\varepsilon_{t}-b_{\tau_{k}}^{*}-t)
-\boldsymbol{\varphi}_{\tau_{k},p}(\varepsilon_{t}-b_{\tau_{k}}^{*}))dt\Big]\\
\nonumber& &\Big(c\int_{0}^{(u_{k}+\textbf{x}'_{t}\textbf{u})
/\sqrt{T}}|t|^{p-1}dt\Big)\\
\nonumber&\leq &E(-B_{T, k})\frac{c}{p}\Bigg(\frac{\mbox{max}_{1\leq t\leq T}|u_{k}+\textbf{x}'_{t}\textbf{u}|}{\sqrt{T}}\Bigg)^{p}\\
&\rightarrow &0.
\end{eqnarray}
The second `$\leq$' above is based on the fact, implied by Lemma 6 in Daouia et al. (2019)\cite{Daouia},
\begin{eqnarray}\label{eq1.15}
\varphi_{\tau_{k},p}(\varepsilon_{t}-b_{\tau_{k}}^{*}-t)
-\varphi_{\tau_{k},p}(\varepsilon_{t}-b_{\tau_{k}}^{*})\leq c|t|^{p-1},
\end{eqnarray}
where $c$ is a positive constant and the last `$\rightarrow$' is due to $\mbox{max}_{1\leq t\leq T}|u_{k}+\textbf{x}'_{t}\textbf{u}|/\sqrt{T}\rightarrow0$, which can be derived from Assumption \ref{ass2.1},
see Pollard (1991)\cite{Pollard} for more details. Combining \eqref{eq1.14} and \eqref{eq1.15} shows
\begin{eqnarray*}\label{eq1.16}
B_{T, k}\stackrel{\textit{P}}{\longrightarrow}-\frac{1}{2}E
\boldsymbol{\psi}_{\tau_{k},p}(\varepsilon-b_{\tau_{k}}^{*})(u_{k}, \textbf{u}')\left( {{\begin{array}{*{20}c}
1 & 0\\
0 & \textbf{C}\\
\end{array}}}
\right)(u_{k}, \textbf{u}')'.
\end{eqnarray*}
So by Slutsky's Theorem, we have
\begin{eqnarray*}\label{eq1.17}
\nonumber Q_{T}&\stackrel{\textit{D}}{\longrightarrow}&-\sum^{K}_{k=1}Z_{k}u_{k}-\textbf{Z}'\textbf{u}
+\frac{1}{2}\sum^{K}_{k=1}E
\boldsymbol{\psi}_{\tau_{k},p}(\varepsilon-b_{\tau_{k}}^{*})(u_{k}, \textbf{u}')\left( {{\begin{array}{*{20}c}
1 & 0\\
0 & \textbf{C}\\
\end{array}}}
\right)(u_{k}, \textbf{u}')'\\
\nonumber &=&-\sum^{K}_{k=1}Z_{k}u_{k}-\textbf{Z}'\textbf{u}
+\frac{1}{2}\sum^{K}_{k=1}E\boldsymbol{\psi}_{\tau_{k},p}(\varepsilon-b_{\tau_{k}}^{*})u^{2}_{k}\\
& &+\frac{1}{2}\sum^{K}_{k=1}E\boldsymbol{\psi}_{\tau_{k},p}(\varepsilon-b_{\tau_{k}}^{*})\textbf{u}'\textbf{C}\textbf{u}.
\end{eqnarray*}
Using the convexity of $Q_{T}$ and Basic Corollary in Hjort and Pollard (1993)\cite{Hjort}, we get
\begin{eqnarray*}\label{eq1.18}
\textbf{u}_{T}\stackrel{\textit{D}}{\longrightarrow}\Bigg(\textbf{C}\sum^{K}_{k=1}
E\boldsymbol{\psi}_{\tau_{k},p}(\varepsilon-b_{\tau_{k}}^{*})\Bigg)^{-1}\textbf{Z}\sim N\Bigg(0, \Bigg(\sum^{K}_{k=1}
E\boldsymbol{\psi}_{\tau_{k},p}(\varepsilon-b_{\tau_{k}}^{*})\Bigg)^{-2}\textbf{C}^{-1}\boldsymbol{\Sigma}_{\textbf{Z}}\textbf{C}^{-1}\Bigg),
\end{eqnarray*}
where
\begin{eqnarray*}\label{eq1.19}
\boldsymbol{\Sigma}_{\textbf{Z}}=\textbf{C}\sum^{K}_{k^{'}=1}\sum^{K}_{k=1}
E[\boldsymbol{\varphi}_{\tau_{k'},p}(\varepsilon-b^{*}_{\tau_{k'}})\boldsymbol{\varphi}_{\tau_{k},p}(\varepsilon-b^{*}_{\tau_{k}})]. \ \Box
\end{eqnarray*}

{\bf Proof of Theorem \ref{thm2.2}.} Divide by $K^{2}$ the numerator and denominator of the fraction in \eqref{eq1.6}. We first consider the resulting denominator and have, for $\tau_{k}=k/(K+1)$,
\begin{eqnarray}\label{eq1.21}
\nonumber \frac{1}{K}\sum_{k=1}^{K}E\boldsymbol{\psi}_{\tau_{k},p}(\varepsilon-b^{*}_{\tau_{k}})
&=&\frac{1}{K}\sum_{k=1}^{K}E(p(p-1)|\tau_{k}-I(\varepsilon<b^{*}_{\tau_{k}})|
|\varepsilon-b^{*}_{\tau_{k}}|^{p-2})\\
\nonumber &\stackrel{K\rightarrow\infty}{\longrightarrow}& \int^{1}_{0}E(p(p-1)|s-I(\varepsilon<F^{-1}_{\varepsilon, p}(s))||\varepsilon-F^{-1}_{\varepsilon, p}(s)|^{p-2})ds\\
\nonumber &=&E_{U_{1}}(E(p(p-1)|U_{1}-I(\varepsilon<F^{-1}_{\varepsilon, p}(U_{1}))||\varepsilon-F^{-1}_{\varepsilon, p}(U_{1})|^{p-2})),\\
& &
\end{eqnarray}
where $U_{1}$ is a random variable obeying the uniform distribution on $[0, 1]$.
Define $\varepsilon_{a}=F^{-1}_{\varepsilon, p}(U_{1})$, the expression \eqref{eq1.21} is further written as
\begin{eqnarray}\label{eq1.22}
p(p-1)E_{\varepsilon_{a}}E_{\varepsilon}
(|F_{\varepsilon, p}(\varepsilon_{a})-I(\varepsilon<\varepsilon_{a})||\varepsilon-\varepsilon_{a}|^{p-2}).
\end{eqnarray}
Second, focus on the numerator and we have
\begin{eqnarray}\label{eq1.23}
\nonumber & & \frac{1}{K^{2}}\sum^{K}_{k'=1}\sum^{K}_{k=1}
E[\boldsymbol{\varphi}_{\tau_{k'},p}(\varepsilon-b^{*}_{\tau_{k'}})\boldsymbol{\varphi}_{\tau_{k},p}(\varepsilon-b^{*}_{\tau_{k}})]\\
&=&\nonumber \frac{1}{K^{2}}\sum^{K}_{k'=1}\sum^{K}_{k=1}
E(p^{2}(\tau_{k'}-I(\varepsilon<b^{*}_{\tau_{k'}}))
(\tau_{k}-I(\varepsilon<b^{*}_{\tau_{k}}))|\varepsilon-b^{*}_{\tau_{k'}}|^{p-1}
|\varepsilon-b^{*}_{\tau_{k}}|^{p-1})\\
\nonumber &\stackrel{K\rightarrow\infty}{\longrightarrow}&
p^{2}\int_{0}^{1}\int_{0}^{1}E((s-I(\varepsilon<F^{-1}_{\varepsilon, p}(s))
(t-I(\varepsilon<F^{-1}_{\varepsilon, p}(t)))\\
\nonumber & & |\varepsilon-F^{-1}_{\varepsilon, p}(s)|^{p-1}
|\varepsilon-F^{-1}_{\varepsilon, p}(t)|^{p-1})dsdt\\
\nonumber&=& p^{2}\int_{0}^{1}\int_{0}^{1}E((s-I(\varepsilon<F^{-1}_{\varepsilon, p}(s))
(t-I(\varepsilon<F^{-1}_{\varepsilon, p}(t)))\\
\nonumber & & |\varepsilon-F^{-1}_{\varepsilon, p}(s)|^{p-1}
|\varepsilon-F^{-1}_{\varepsilon, p}(t)|^{p-1})dsdt\\
\nonumber &=&p^{2}E_{\varepsilon_{c}}E_{\varepsilon_{b}}E_{\varepsilon}((F_{\varepsilon, p}(\varepsilon_{c})-I(\varepsilon<\varepsilon_{c}))
(F_{\varepsilon, p}(\varepsilon_{b})-I(\varepsilon<\varepsilon_{b}))
|\varepsilon-\varepsilon_{b}|^{p-1}|\varepsilon-\varepsilon_{c}|^{p-1}),\\
&&
\end{eqnarray}
where $\varepsilon_{b}=F^{-1}_{\varepsilon, p}(U_{2})$, $\varepsilon_{c}=F^{-1}_{\varepsilon, p}(U_{3})$, $U_{2}$ and $U_{3}$ are two random variables obeying the uniform distribution on $[0, 1]$. The $U_{i}$, $i=1, 2, 3$ are mutually independent. Combining \eqref{eq1.22} and \eqref{eq1.23} completes the proof. $\Box$

{\bf Proof of Theorem \ref{thm3.1}}. Let $\sqrt{T}(\hat{\boldsymbol{\beta}}^{Aclp}-\boldsymbol{\beta}^{*})=\textbf{u}_{T}$ and
$\sqrt{T}(\hat{b}_{k}-b^{*}_{\tau_{k}})=u_{T, k}$. We can get $(u_{T, 1}, u_{T, 2}, \cdots, u_{T, K}, \textbf{u}_{T})$ by minimizing the following criterion function
\begin{eqnarray}\label{eq1.25}
\nonumber Q_{T}&=&\sum_{k=1}^{K}\sum^{T}_{t=1}
\Big[\boldsymbol{\eta}_{\tau_{k},p}\Big(\varepsilon_{t}-b_{\tau_{k}}^{*}-\frac{u_{k}+\textbf{x}'_{t}\textbf{u}}{\sqrt{T}}\Big)
-\boldsymbol{\eta}_{\tau_{k},p}(\varepsilon_{t}-b_{\tau_{k}}^{*})\Big]\\
\nonumber & & +\sum^{m}_{j=1}\frac{\lambda_{T}}
{\sqrt{T}|\hat{\beta}^{clp}_{j}|^{2}}\sqrt{T}\Big[\Big|\beta^{*}_{j}+\frac{u_{j}}{\sqrt{T}}\Big|
-|\beta^{*}_{j}|\Big].
\end{eqnarray}
As in the proof of Theorem \ref{thm2.1}, the function can be written as
\begin{eqnarray}\label{eq1.26}
\nonumber Q_{T}=-\sum_{k=1}^{K}Z_{T, k}u_{k}-\textbf{Z}'_{T}\textbf{u}-\sum_{K}^{k=1}B_{T, k}
+\sum^{m}_{j=1}\frac{\lambda_{T}}
{\sqrt{T}|\hat{\beta}^{clp}_{j}|^{2}}\Big[\Big|\beta^{*}_{j}+\frac{u_{j}}{\sqrt{T}}\Big|
-|\beta^{*}_{j}|\Big].
\end{eqnarray}
About the penalty term in the above expression, if $\beta^{*}_{j}\neq0$, then $|\hat{\beta}^{clp}_{j}|^{2}\rightarrow|\beta^{*}_{j}|^{2}$ in probability and ${\sqrt{T}|\hat{\beta}^{clp}_{j}|^{2}}\Big[\Big|\beta^{*}_{j}+\frac{u_{j}}{\sqrt{T}}\Big|
-|\beta^{*}_{j}|\Big]\rightarrow u_{j}\mbox{sgn}(\beta^{*}_{j})$. Slutsky's theorem makes sure
$\frac{\lambda_{T}}
{\sqrt{T}|\hat{\beta}^{clp}_{j}|^{2}}\sqrt{T}\Big[\Big|\beta^{*}_{j}+\frac{u_{j}}{\sqrt{T}}\Big|
-|\beta^{*}_{j}|\Big]\rightarrow0$ in probability. If $\beta^{*}_{j}=0$ then $\sqrt{T}\Big[\Big|\beta^{*}_{j}+\frac{u_{j}}{\sqrt{T}}\Big|
-|\beta^{*}_{j}|\Big]=|u_{j}|$ and $\frac{\lambda_{T}}
{\sqrt{T}|\hat{\beta}^{clp}_{j}|^{2}}=\frac{\sqrt{T}\lambda_{T}}
{(\sqrt{T}|\hat{\beta}^{clp}_{j}|)^{2}}\rightarrow\infty$ in probability. So we have
\begin{equation*}\label{eq1.27}
\frac{\lambda_{T}}
{\sqrt{T}|\hat{\beta}^{clp}_{j}|^{2}}\sqrt{T}\Big[\Big|\beta^{*}_{j}+\frac{u_{j}}{\sqrt{T}}\Big|
-|\beta^{*}_{j}|\Big]\stackrel{\textit{P}}{\longrightarrow}V(\beta_{j}, u_{j})=
\begin{cases}
0,&\mbox{if}\ \beta^{*}_{j}\neq 0,\\
0,&\mbox{if}\ \beta^{*}_{j}=0\ \mbox{and}\ u_{j}=0,\\
\infty,& \mbox{if}\ \beta^{*}_{j}=0\ \mbox{and}\ u_{j}\neq0.
\end{cases}
\end{equation*}
Additionally, using the same argument in the proof of Theorem \ref{thm2.1}, we have
\begin{eqnarray*}\label{eq1.28}
\nonumber Q_{T}&\stackrel{\textit{D}}{\longrightarrow}&-\sum^{K}_{k=1}Z_{k}u_{k}-\textbf{Z}'\textbf{u}
+\frac{1}{2}\sum^{K}_{k=1}E\boldsymbol{\psi}_{\tau_{k},p}(\varepsilon-b_{\tau_{k}}^{*})u^{2}_{k}\\
& &+\frac{1}{2}\sum^{K}_{k=1}E\boldsymbol{\psi}_{\tau_{k},p}(\varepsilon-b_{\tau_{k}}^{*})\textbf{u}'C\textbf{u}+\sum^{m}_{j=1}V(\beta_{j}, u_{j}).
\end{eqnarray*}

Write $\textbf{u}=(\textbf{u}'_{1}, \textbf{u}'_{2})'$ where $\textbf{u}_{1}$ contains the first $q$ elements of $\textbf{u}$ which corresponds to the $q$ non-zero $\beta_{j}^{*}$, $j\in\mathcal{A}$ in terms of indice. Using the same arguments in Knight (1998)\cite{Knight} and thoughts in Theorem \ref{thm2.1}, we have
\begin{eqnarray}\label{eq1.29}
\hat{\textbf{u}}_{2, T}\stackrel{\textit{D}}{\longrightarrow}0
\end{eqnarray}
and
\begin{eqnarray}\label{eq1.30}
\nonumber\hat{\textbf{u}}_{1, T}&\stackrel{\textit{D}}{\longrightarrow}&\Bigg(\textbf{C}_{\mathcal{A}\mathcal{A}}\sum^{K}_{k=1}
E\boldsymbol{\psi}_{\tau_{k},p}(\varepsilon-b_{\tau_{k}}^{*})\Bigg)^{-1}\textbf{Z}\\
& &\sim N\Bigg(0, \Bigg(\sum^{K}_{k=1}
E\boldsymbol{\psi}_{\tau_{k},p}(\varepsilon-b_{\tau_{k}}^{*})\Bigg)^{-2}\textbf{C}_{\mathcal{A}\mathcal{A}}^{-1}
\boldsymbol{\Sigma}_{\textbf{Z}1}\textbf{C}_{\mathcal{A}\mathcal{A}}^{-1}\Bigg)£¬
\end{eqnarray}
where
\begin{eqnarray*}\label{eq1.31}
\boldsymbol{\Sigma}_{\textbf{Z}1}=\textbf{C}_{\mathcal{A}\mathcal{A}}\sum^{K}_{k^{'}=1}\sum^{K}_{k=1}
E[\boldsymbol{\varphi}_{\tau_{k'},p}(\varepsilon-b^{*}_{\tau_{k'}})\boldsymbol{\varphi}_{\tau_{k},p}(\varepsilon-b^{*}_{\tau_{k}})].
\end{eqnarray*}
Hence, the desired asymptotic normality holds.

Next we focus on the consistent selection property. Define $\hat{\mathcal{A}}_{T}=\{j: \hat{\beta}^{Aclp}_{j}\neq0\}$. $\forall j\in \mathcal{A}$, the asymptotic normality implies
$P(j\in\hat{\mathcal{A}}_{T})\rightarrow1$.
We only need to show that $\forall j\notin \mathcal{A}$, $P(j\in\hat{\mathcal{A}}_{T})\rightarrow0$.
When $j'\in\hat{\mathcal{A}}_{T}$, according to the KKT optimality conditions, we have
\begin{eqnarray*}\label{eq1.33}
\sum^{K}_{k=1}\sum^{T}_{t=1}\boldsymbol{\eta}'_{\tau_{k},p}(y_{t}-b_{k}-\textbf{x}'_{t}\hat{\boldsymbol{\beta}}^{Aclp})x_{t,j'}
=\frac{\lambda(t)}{|\hat{\beta}_{j'}^{clp}|^2}.
\end{eqnarray*}
By the $c_{p}$-inequality, the left-hand side of the above equation is not larger than
\begin{eqnarray*}\label{eq1.34}
c_{p}p\sum^{K}_{k=1}\sum^{T}_{t=1}(|(\varepsilon_{t}-b_{k})|x_{t,j'}|^{1/(p-1)}|^{p-1}+
|\textbf{x}'_{t}(\boldsymbol{\beta}^{*}-\hat{\boldsymbol{\beta}}^{Aclp})|x_{t,j'}|^{1/(p-1)}|^{p-1}).
\end{eqnarray*}
By \eqref{eq1.29}, \eqref{eq1.30} and Slutsky's theorem, we have
\begin{eqnarray*}\label{eq1.34+1}
\sum^{K}_{k=1}\frac{1}{T^{3/2}}\sum^{T}_{t=1}|\textbf{x}'_{t}|x_{t,j'}|^{1/(p-1)}\sqrt{T}
(\boldsymbol{\beta}^{*}-\hat{\boldsymbol{\beta}}^{Aclp})|^{p-1}
\rightarrow0
\end{eqnarray*}
and
\begin{eqnarray*}\label{eq1.34+2}
& &\sum^{K}_{k=1}\frac{1}{T}\sum^{T}_{t=1}|\varepsilon_{t}-b_{k}|^{p-1}|x_{t,j'}|\rightarrow
\sum^{K}_{k=1}E(|\varepsilon-b_{k}|^{p-1})\frac{1}{T}\sum^{T}_{t=1}|x_{t,j'}|\\
& &\leq\sum^{K}_{k=1}E(|\varepsilon-b_{k}|^{p-1})\sqrt{\frac{1}{T}\sum^{T}_{t=1}|x_{t,j'}|^{2}}
\rightarrow\sum^{K}_{k=1}E(|\varepsilon-b_{k}|^{p-1})\sqrt{C_{j'j'}}.
\end{eqnarray*}
But the condition of the theorem shows
\begin{eqnarray*}\label{eq1.35}
\frac{\lambda(t)T^{\frac{p-2}{2}}}{|\sqrt{T}\hat{\beta}_{j'}^{clp}|^2}\rightarrow\infty.
\end{eqnarray*}
So
\begin{eqnarray*}\label{eq1.36}
P(j\in\hat{\mathcal{A}}_{T})\leq P\Bigg(\sum^{K}_{k=1}\sum^{T}_{t=1}\boldsymbol{\eta}'_{\tau_{k},p}(y_{t}-b_{k}-\textbf{x}'_{t}\hat{\boldsymbol{\beta}}^{Aclp})x_{t,j'}
=\frac{\lambda(t)}{|\hat{\beta}_{j'}^{clp}|^2}\Bigg)\rightarrow0
\end{eqnarray*} $\Box$

We need the following lemmas to complete the proof of Theorem \ref{thm4.1}. Define
\begin{eqnarray*}\label{eq1.41}
Z_{T, p}(\boldsymbol{\delta})=\sum^{T}_{t=1}(\boldsymbol{\eta}_{\tau,p}(u_{t}-\textbf{x}'_{t}\boldsymbol{\delta}/\sqrt{T})
-\boldsymbol{\eta}_{\tau,p}(u_{t})),
\end{eqnarray*}
where $u_{t}=y_{t}-\textbf{x}'_{t}\boldsymbol{\beta}_{0}$.

\begin{lemma}\label{lem1}
Under model \eqref{eq1.37} and Assumption \ref{ass4.1} we have
\begin{eqnarray*}\label{eq1.45}
|Z_{T, 1}(\boldsymbol{\delta})-\frac{1}{2}f(0)\boldsymbol{\delta}'\textbf{D}_{0}\boldsymbol{\delta}-\textbf{W}'_{T}
\boldsymbol{\delta}|\stackrel{\textit{P}}{\longrightarrow}0.
\end{eqnarray*}
\end{lemma}

Proof. According to Zou and Yuan (2008)\cite{ZOUY}, we have
\begin{eqnarray*}\label{eq1.46}
\boldsymbol{\eta}_{\tau, 1}(r-s)-\boldsymbol{\eta}_{\tau, 1}(r)=s(I(r<0)-\tau)+\int^{s}_{0}(I(r\leq t)-I(r\leq0))dt.
\end{eqnarray*}
Using this identity, we write
\begin{eqnarray}\label{eq1.47}
\nonumber Z_{T, 1}(\boldsymbol{\delta})&=&\sum^{T}_{t=1}\frac{\textbf{x}'_{t}\boldsymbol{\delta}}{\sqrt{T}}(I(u_{t}<0)-\tau)
+\sum^{T}_{t=1}\int^{\textbf{x}'_{t}\boldsymbol{\delta}/\sqrt{T}}_{0}(I(u_{t}\leq t)-I(u_{t}\leq 0))dt\\
&=:&\textbf{W}'_{T}\boldsymbol{\delta}+B_{T}.
\end{eqnarray}
Further, we have, with $F$ being the cumulative distribution function of $u_{t}$,
\begin{eqnarray*}\label{eq1.48}
E(B_{T})&=&\sum^{T}_{t=1}\int^{\textbf{x}'_{t}\boldsymbol{\delta}/\sqrt{T}}_{0}(F(t)-F(0))dt\\
&=&\frac{1}{T}\sum^{T}_{t=1}\int^{\textbf{x}'_{t}\boldsymbol{\delta}}_{0}\sqrt{T}(F(t/\sqrt{T})-F(0))dt\\
&=&\frac{1}{T}\sum^{T}_{t=1}\int^{\textbf{x}'_{t}\boldsymbol{\delta}}_{0}f(rt/\sqrt{T})tdt,
\end{eqnarray*}
where $|r|<1$. Based on the property that $f(u)$ is continuous in a neighborhood of 0, it clear that
$f(rt/\sqrt{T})$ converges to $f(0)$ uniformly in $|rt|\in [0, \textbf{x}'_{t}\boldsymbol{\delta}]$ and thus
$E(B_{T})\rightarrow(f(0)\boldsymbol{\delta}'\textbf{D}_{0}\boldsymbol{\delta})/2$. Using the same argument as in the proof of Theorem 2.1 in Zou and Yuan (2008)\cite{ZOUY}, we can show $\mbox{Var}(B_{T})\rightarrow0$ and hence $B_{T}\rightarrow(f(0)\boldsymbol{\delta}'\textbf{D}_{0}\boldsymbol{\delta})/2$ in probability. Combining this and \eqref{eq1.47}, the desired
result is obtained. $\Box$

\begin{lemma}\label{lem2}
Under Assumptions \ref{ass4.2}-\ref{ass4.4}, when $p\rightarrow1+$, we have the following two convergence results.
\begin{eqnarray}\label{eq1.49}
E\boldsymbol{\psi}_{\tau, p}(u+\alpha_{p}+q_{\tau})\longrightarrow f(q_{\tau}),
\end{eqnarray}
where $\alpha_{p}\rightarrow0$ as $p\rightarrow1+$, $q_{\tau}$ is the $\tau$th-quantile of $u$,
and the definition of $\boldsymbol{\psi}_{\tau, p}(s)$ can be found in Theorem \ref{ass2.1}.
Moreover,
\begin{eqnarray}\label{eq1.49+0}
E\boldsymbol{\psi}_{\tau, p}(u+\alpha+q_{\tau})\longrightarrow f(q_{\tau}-\alpha),
\end{eqnarray}
where $\alpha$ is a constant.
\end{lemma}

Proof. First we focus on the proof of the limit in \eqref{eq1.49}. Without the loss of generality, we consider the case of $q_{\tau}=0$. According to Assumption \ref{ass4.3}, it is easily to show $E(|u+c|^{p-2})<\infty$ for a suitable constant $c$. So we can write
\begin{eqnarray*}\label{eq1.49+1}
E\boldsymbol{\psi}_{\tau, p}(u+\alpha_{p})&=&p(p-1)E(|\tau-I(u<0)||u+\alpha_{p}|^{p-2})\\
&=&p(p-1)\int_{0}^{\infty}(\tau x^{p-2}f(x-\alpha_{p})+(1-\tau)x^{p-2}f(-x-\alpha_{p}))dx\\
&=&p\int_{0}^{\infty}(\tau f(x-\alpha_{p})+(1-\tau)f(-x-\alpha_{p}))dx^{p-1}\\
&=&p(\tau f(x-\alpha_{p})+(1-\tau)f(-x-\alpha_{p}))x^{p-1}|^{\infty}_{0}\\
& &-p\int_{0}^{\infty}x^{p-1}(\tau f^{(1)}(x-\alpha_{p})-(1-\tau)f^{(1)}(-x-\alpha_{p}))dx\\
&=&-p\int_{0}^{\infty}x^{p-1}(\tau f^{(1)}(x-\alpha_{p})-(1-\tau)f^{(1)}(-x-\alpha_{p}))dx.
\end{eqnarray*}
The last equality is based on Assumption \ref{ass4.2}. In fact we have $E(|u_{t}|^{p-1})<\infty$ and thus $f(x)|x|^{p-1}\rightarrow0$ and further $f(x-\alpha_{p})|x|^{p-1}\rightarrow0$ as $x\rightarrow\pm\infty$. We have
\begin{eqnarray*}\label{eq1.49+1+1}
p\tau\int_{0}^{\infty}x^{p-1}f^{(1)}(x-\alpha_{p})dx=p\tau\int_{-\alpha_{p}}^{\infty}
(x+\alpha_{p})^{p-1}f^{(1)}(x)dx,
\end{eqnarray*}
and for $p_{0}\in (1, \triangle)$, when $p\leq p_{0}$
\begin{equation*}\label{eq1.49+2}
(x+\alpha_{p})^{p-1}\leq
\begin{cases}
2x^{p_{0}-1},& x>1,\\
2,& \max\{-\alpha_{p}, 0\}<x\leq 1,\\
1,& \min\{-\alpha_{p}, 0\}<x\leq 0.
\end{cases}
\end{equation*}
Using Assumption \ref{ass4.4}, Heine's theorem and the Lebesque control-convergent theorem, we get
$p\int_{0}^{\infty}x^{p-1}\tau f^{(1)}(x-\alpha_{p})dx\rightarrow\int_{0}^{\infty}\tau f^{(1)}(x)dx$ as $p\rightarrow1+$. Similarly, $p\int_{0}^{\infty}x^{p-1}(1-\tau)f^{(1)}(-x-\alpha_{p}))dx\rightarrow\int_{0}^{\infty}(1-\tau)f^{(1)}(-x))dx$.
So we have
\begin{equation*}\label{eq1.49+3}
E\boldsymbol{\psi}_{\tau, p}(u+\alpha_{p})\longrightarrow -\int_{0}^{\infty}(\tau f^{(1)}(x)-(1-\tau)f^{(1)}(-x))dx=f(0).
\end{equation*}
The proof for \eqref{eq1.49+0} is the same as that for \eqref{eq1.49} and so we omit it. $\Box$

\begin{lemma}\label{lem3}
Under the model \eqref{eq1.37} and Assumptions \ref{ass4.2}-\ref{ass4.4}, for any $\varsigma>0$ and $\varepsilon>0$, $\exists p_{0}>0$ and $N>0$, when $0< p-1\leq p_{0}$ and $T>N$, we have
\begin{eqnarray*}\label{eq1.50}
P(|Z_{T, p}(\delta)-Z_{T, 1}(\delta)|\geq \varepsilon)<\varsigma.
\end{eqnarray*}
\end{lemma}

Proof. Using the arguments in the proofs of Theorem \ref{thm2.1} and Lemma \ref{lem1}, we have
\begin{eqnarray*}\label{eq1.51}
Z_{T, p}(\delta)-Z_{T, 1}(\delta)&=&-\frac{1}{\sqrt{T}}\sum^{T}_{t=1}\textbf{x}'_{t}\boldsymbol{\delta}(\boldsymbol{\varphi}_{\tau, p}(u_{t}) +I(u_{t}<0)-\tau)\\
& &-\sum^{T}_{t=1}\int^{\textbf{x}'_{t}\boldsymbol{\delta}/\sqrt{T}}_{0}(\boldsymbol{\varphi}_{\tau, p}(u_{t}-s)-\boldsymbol{\varphi}_{\tau, p}(u_{t}))ds\\
& &-\sum^{T}_{t=1}\int^{\textbf{x}'_{t}\boldsymbol{\delta}/\sqrt{T}}_{0}(I(u_{t}\leq s)-I(u_{t}\leq 0))ds\\
&=:&I+II+III,
\end{eqnarray*}
and
\begin{eqnarray*}\label{eq1.52}
A_{T}:=E(II+III)\longrightarrow\frac{1}{2}\boldsymbol{\delta}'\textbf{D}_{0}\boldsymbol{\delta}(E\boldsymbol{\psi}_{\tau, p}(u)-f(0))=:A.
\end{eqnarray*}
Let $\varepsilon_{1}<\varepsilon$. According to Lemma \ref{lem2}, $|\boldsymbol{\delta}'\textbf{D}_{0}\boldsymbol{\delta}(E\boldsymbol{\psi}_{\tau, p}(u)-f(0))/2|\rightarrow0$ as $p\rightarrow1+$, so there are a $p_{01}$ and  a positive $N_{1}$ such that if $0< p-1\leq p_{01}$ and $T>N_{1}$,
$|A_{T}|<\varepsilon_{1}/2$.
So we have
\begin{eqnarray}\label{eq1.53}
\nonumber P(|II+III|\geq\varepsilon_{1})&=&P(II+III+\varepsilon_{1}/2\leq -\varepsilon_{1}/2)+
P(II+III-\varepsilon_{1}/2\geq\varepsilon_{1}/2)\\
\nonumber &\leq&P(II+III-A_{T}\leq -\varepsilon_{1}/2)+
P(II+III-A_{T}\geq\varepsilon_{1}/2)\\
\nonumber &=&P(|II+III-A_{T}|\geq\varepsilon_{1}/2)\leq \frac{D(II+III)}{(\varepsilon_{1}/2)^{2}}\\
\nonumber &\leq&\frac{8}{\varepsilon_{1}^{2}}(D(II)+D(III))\\
\nonumber &\leq&\frac{8}{\varepsilon_{1}^{2}}\Big(\frac{4}{T}\sum^{T}_{t=1}
\int^{\textbf{x}'_{t}\boldsymbol{\delta}}_{0}\sqrt{T}(F(t/\sqrt{T})-F(0))dt\cdot\max_{1\leq t\leq T}\Big\{\frac{\textbf{x}'_{t}\boldsymbol{\delta}}{\sqrt{T}}\Big\}\\
& & +\frac{1}{T}\sum^{T}_{t=1}
\int^{\textbf{x}'_{t}\boldsymbol{\delta}}_{0}E\boldsymbol{\psi}_{\tau, p}(u-\tilde{t}/\sqrt{T})tdt\cdot\frac{c}{p}\Big(\max_{1\leq t\leq T}\Big\{\frac{\textbf{x}'_{t}\boldsymbol{\delta}}{\sqrt{T}}\Big\}\Big)^{p}\Big).
\end{eqnarray}
According to the proof of Lemma \ref{lem1}, the first term of the right-hand side in \eqref{eq1.53} converges to zero.
According to Lemma \ref{lem2}, when $0< p-1\leq p_{02}$, $E\boldsymbol{\psi}_{\tau, p}(u-\tilde{t}/\sqrt{T})\leq (f(0)+c)(1+o(1))$ with $o(1)$ holds uniformly for $\tilde{t}$ between $\textbf{x}'_{t}\boldsymbol{\delta}$ and 0. Then using the same argument in the proof of Theorem \ref{thm2.1}, the second term of the right-hand side in \eqref{eq1.53} also converges to zero. So, there exists $N_{2}$ such that when $T>N_{2}$,
\begin{eqnarray}\label{eq1.53+1}
P(|II+III|\geq\varepsilon_{1})<\varsigma/2.
\end{eqnarray}
Then, using Markov's inequality and noting $0<\varepsilon_{1}<\varepsilon$ we have
\begin{eqnarray}\label{eq1.54}
\nonumber & & P(|Z_{T, p}(\boldsymbol{\delta})-Z_{T, 1}(\boldsymbol{\delta})|\geq \varepsilon)\\
\nonumber & &\leq P(|I+II+III|\geq \varepsilon, |II+III|<\varepsilon_{1})+P(|II+III|\geq\varepsilon_{1})\\
\nonumber & &\leq P(|I|\geq \varepsilon-\varepsilon_{1})+ P(|II+III|\geq\varepsilon_{1})\\
& &\leq \frac{\boldsymbol{\delta}'\frac{1}{T}\sum^{T}_{t=1}\textbf{x}_{i}\textbf{x}'_{i}\boldsymbol{\delta} E(\boldsymbol{\varphi}_{\tau, p}(u_{t}) +I(u_{t}<0)-\tau)^{2}}{(\varepsilon-\varepsilon_{1})^{2}}+P(|II+III|\geq\varepsilon_{1}).
\end{eqnarray}
From the definition of $\boldsymbol{\varphi}_{\tau, p}(s)$ in Theorem \ref{thm2.1}, $\boldsymbol{\varphi}_{\tau, p}(u_{t})+I(u_{t}<0)-\tau$ converges to 0 almost surely. Combining this and Assumption \ref{ass4.1}, there are $p_{03}$ and $N_{3}$ such that when $0\leq p-1\leq p_{03}$
and $T>N_{3}$ we have the first term in \eqref{eq1.54} is not larger than $\varsigma/2$. Combining this,  \eqref{eq1.53+1} and \eqref{eq1.54}, letting $p_{0}=\min\{p_{01}, p_{02}, p_{03}\}$ and $N=\max\{N_{1}, N_{2}, N_{3}\}$,
we complete the proof. $\Box$

{\bf Proof of Theorem \ref{thm4.1}}. Clearly,
$$\hat{\boldsymbol{\delta}}_{T, p}=\sqrt{T}(\hat{\boldsymbol{\beta}}_{T, p}-\boldsymbol{\beta}(\tau))=
\arg\min_{\boldsymbol{\delta}}Z_{T, p}(\boldsymbol{\delta}).$$
We firstly need to prove, for each compact set $K\in R^{d}$,
\begin{eqnarray}\label{eq1.42}
\lim_{T\rightarrow\infty
\atop p\rightarrow1+}\sup_{\boldsymbol{\delta}\in K}\Big|Z_{T, p}(\boldsymbol{\delta})-\frac{1}{2}f(0)\boldsymbol{\delta}'\textbf{D}_{0}\boldsymbol{\delta}-\textbf{W}'_{T}\boldsymbol{\delta}\Big|=0
\end{eqnarray}
in probability, where
\begin{eqnarray*}\label{eq1.43}
\textbf{W}_{T}=\frac{1}{\sqrt{T}}\sum^{T}_{t=1}\textbf{x}_{t}(I(u_{t}<0)-\tau)\stackrel{\textit{D}}{\longrightarrow}
N(0, \tau(1-\tau)\textbf{D}_{0}).
\end{eqnarray*}
The expression in the left-hand side of \eqref{eq1.42} is not larger than
\begin{eqnarray}\label{eq1.44}
|Z_{T, p}(\boldsymbol{\delta})-Z_{T, 1}(\boldsymbol{\delta})|+\Big|Z_{T, 1}(\boldsymbol{\delta})-\frac{1}{2}f(0)\boldsymbol{\delta}'\textbf{D}_{0}\boldsymbol{\delta}-\textbf{W}'_{T}\boldsymbol{\delta}\Big|.
\end{eqnarray}
Considering the second term in \eqref{eq1.44}, Lemmas \ref{lem1} and \ref{lem3} together, we have
\begin{eqnarray*}\label{eq1.55}
\lim_{T\rightarrow\infty
\atop p\rightarrow1+}Z_{T, p}(\boldsymbol{\delta})-\frac{1}{2}f(0)\boldsymbol{\delta}'\textbf{D}_{0}\boldsymbol{\delta}-\textbf{W}'_{T}\boldsymbol{\delta}=0.
\end{eqnarray*}
Based on this, we can use the train of thought in the proof of the convexity lemma in Pollard (1991)\cite{Pollard} to prove \eqref{eq1.42} as $Z_{T, p}(\boldsymbol{\delta})$ is the convex function of $\boldsymbol{\delta}$. Although the argument in Section 6 of Pollard (1991)\cite{Pollard} involves the limit only related to sample size, essentially, it has nothing to do with how the limit is calculated. The detailed argument is omitted.

Define $\boldsymbol{\eta}_{T}=\textbf{D}_{0}^{-1}\textbf{W}_{T}/f(0)$. It is sufficient to prove for each $\zeta>0$ that
\begin{eqnarray*}\label{eq1.56}
\lim_{T\rightarrow\infty
\atop p\rightarrow1+}P(|\hat{\boldsymbol{\delta}}_{T, p}-\boldsymbol{\eta}_{T}|>\zeta)\rightarrow 0.
\end{eqnarray*}
To this end, we further write
\begin{eqnarray*}\label{eq1.57}
Z_{T, p}(\boldsymbol{\delta})=\frac{1}{2}f(0)(\boldsymbol{\delta}-\boldsymbol{\eta}_{T})'
\textbf{D}_{0}(\boldsymbol{\delta}-\boldsymbol{\eta}_{T})-\frac{1}{2}f(0)\boldsymbol{\eta}'_{T}\textbf{D}_{0}\boldsymbol{\eta}_{T}
+r_{T}(\boldsymbol{\delta}),
\end{eqnarray*}
where, for each compact set $K$ in $R^{d}$,
\begin{eqnarray*}\label{eq1.58}
\lim_{T\rightarrow\infty
\atop p\rightarrow1+}\sup_{\boldsymbol{\delta}\in K}|r_{T}(\boldsymbol{\delta})|=0 \ \mbox{in \ probability}.
\end{eqnarray*}
Let $B(T)$ be a closed ball with center $\boldsymbol{\eta}_{T}$ and radius $\zeta$. The random boundedness of $\boldsymbol{\eta}_{T}$ makes sure that there is the compact set $K$ that contains $B(T)$ with probability
arbitrarily close to one, so we have
\begin{eqnarray*}\label{eq1.59}
\lim_{T\rightarrow\infty
\atop p\rightarrow1+}\triangle_{T}=0 \ \mbox{in \ probability},
\end{eqnarray*}
where $\triangle_{T}=\sup_{\boldsymbol{\delta}\in B(T)}|r_{T}(\boldsymbol{\delta})|$.

Next examine the property of $Z_{T, p}(\delta)$ outside $B(T)$. Denote any point outside $B(T)$ by
$\boldsymbol{\delta}=\boldsymbol{\eta}_{T}+\alpha\boldsymbol{\upsilon}$, with $\alpha>\zeta$ and $\boldsymbol{\upsilon}$ a $d$-dimensional unit vector.
$\boldsymbol{\delta}^{*}$ stands for the boundary point of $B(T)$ that just lies on the line segment from $\boldsymbol{\eta}_{T}$ to
$\boldsymbol{\delta}$, namely $\boldsymbol{\delta}^{*}=\boldsymbol{\eta}_{T}+\zeta\boldsymbol{\upsilon}$. Convexity of $Z_{T, p}(\boldsymbol{\delta})$ and definition of $\triangle_{T}$ yield
\begin{eqnarray*}\label{eq1.60}
\frac{\zeta}{\alpha}Z_{T, p}(\boldsymbol{\delta})+\Big(1-\frac{\zeta}{\alpha}\Big)Z_{T, p}(\boldsymbol{\eta}_{T})&\geq& Z_{T, p}(\boldsymbol{\delta}^{*})\\
&\geq& \frac{1}{2}f(0)(\zeta\boldsymbol{\upsilon})'\textbf{D}_{0}(\zeta\boldsymbol{\upsilon})-\frac{1}{2}f(0)
\boldsymbol{\eta}'_{T}\textbf{D}_{0}\boldsymbol{\eta}_{T}-\triangle_{T}\\
&\geq& \frac{1}{2}f(0)(\zeta\boldsymbol{\upsilon})'\textbf{D}_{0}(\zeta\boldsymbol{\upsilon})+Z_{T, p}(\boldsymbol{\eta}_{T})-2\triangle_{T}.
\end{eqnarray*}
So we have
\begin{eqnarray*}\label{eq1.61}
\inf_{|\boldsymbol{\delta}-\boldsymbol{\eta}_{T}|>\zeta}Z_{T, p}(\boldsymbol{\delta})\geq Z_{T, p}(\boldsymbol{\eta}_{T})+\frac{\alpha}{\zeta}\Big(\frac{1}{2}f(0)(\zeta\boldsymbol{\upsilon})'
\textbf{D}_{0}(\zeta\boldsymbol{\upsilon})-2\triangle_{T}\Big).
\end{eqnarray*}
With probability tending to one, $\frac{1}{2}f(0)(\zeta\boldsymbol{\upsilon})'\textbf{D}_{0}(\zeta\boldsymbol{\upsilon})>2\triangle_{T}$,
thus the minimum of $Z_{T, p}(\boldsymbol{\delta})$ cannot appear at any $\boldsymbol{\delta}$ outside $B(T)$ and in other words $|\hat{\boldsymbol{\delta}}_{T, p}-\boldsymbol{\eta}_{T}|\leq\zeta$ with probability tending to 1 as $T\rightarrow\infty$ and $ p\rightarrow1+$ simultaneously. The proof of Theorem 4.1 is completed. $\Box $

The proof of Theorem \ref{thm4.2} needs the following lemma.
\begin{lemma}\label{lem4}
If $E|\varepsilon|^{p-1}<\infty$, for $p\in (1, \triangle)$, we have, as $p\rightarrow1+$, the
$\tau$th $L^{p}$-quantile of $\varepsilon$ converges to its $\tau$th quantile, namely,
\begin{eqnarray*}\label{eq1.62}
q^{lp}_{\varepsilon}(\tau)\rightarrow q_{\varepsilon}(\tau),
\end{eqnarray*}
where $q^{lp}_{\varepsilon}(\tau)=\max_{s}\{E(|\varepsilon-s|^{p}|\tau-I(\varepsilon<s)|)
-E(|\varepsilon|^{p}|\tau-I(\varepsilon<0)|)\}$ and
$q_{\varepsilon}(\tau)=\max_{s}\{E(|\varepsilon-s||\tau-I(\varepsilon<s)|)
-E(|\varepsilon||\tau-I(\varepsilon<0)|)\}$.
\end{lemma}

Proof. Firstly, we have, $|r|\leq1$,
\begin{eqnarray*}\label{eq1.63}
& & E(|\varepsilon-s|^{p}|\tau-I(\varepsilon<s)|)
-E(|\varepsilon|^{p}|\tau-I(\varepsilon<0)|)\\
&=&E(p|\tau-I(\varepsilon<rs)||\varepsilon-rs|^{p-1}\mbox{sign}(\varepsilon-rs)(-s)),
\end{eqnarray*}
and, when $p_{0}\in (1, \triangle)$ and $p_{0}>p$
\begin{eqnarray*}\label{eq1.64}
& & |p|\tau-I(\varepsilon<rs)||\varepsilon-rs|^{p-1}\mbox{sign}(\varepsilon-rs)(-s))|\\
&\leq & g(\varepsilon)=
\begin{cases}
p|s||\varepsilon-rs|^{p_{0}-1},& |\varepsilon-rs|>1,\\
p|s|,& 0<|\varepsilon-rs|\leq1.
\end{cases}
\end{eqnarray*}
Then based on $E|\varepsilon|^{p-1}<\infty$ and the $c_{p}$-inequality, we have $Eg<\infty$ and hence
\begin{eqnarray*}\label{eq1.65}
& &Q^{lp}(s):=E(|\varepsilon-s|^{p}|\tau-I(\varepsilon<s)|)
-E(|\varepsilon|^{p}|\tau-I(\varepsilon<0)|)\\
& & \longrightarrow Q(s):=E(|\varepsilon-s||\tau-I(\varepsilon<s)|)
-E(|\varepsilon||\tau-I(\varepsilon<0)|)
\end{eqnarray*}
by Heine's theorem and the Lebesque control-convergent theorem. Defining $r_{p}(s)=Q^{lp}(s)-Q(s)$ and using Theorem 10.8 in Rockafellar (1970)\cite{Rockafellar} or the same argument in the proof of the convexity lemma as in Pollard (1991)\cite{Pollard} but for the nonstochastic case, we further get, as $p\rightarrow1+$,
\begin{eqnarray}\label{eq1.66}
\sup_{s\in B}|r_{p}(s)|\rightarrow0,
\end{eqnarray}
where $B$ is any compact subset of $R$.

Next, we show that, for any $\varsigma>0$, there will be a $\epsilon>0$ such that if $0<p-1<\epsilon$,
$q^{lp}_{\varepsilon}(\tau)\in (q_{\varepsilon}(\tau), \varsigma)$. Let $t$ be any point outside $U(q_{\varepsilon}(\tau), \varsigma)$ and may write $t=q_{\varepsilon}(\tau)+\kappa e$ with $e$ a unit vector and $\kappa>\varsigma$. The intersection of the line segment from $q_{\varepsilon}(\tau)$ to $t$ and the boundary of $U(q_{\varepsilon}(\tau)$ is $q_{\varepsilon}(\tau)+\varsigma e$, which can be written as
$(1-\frac{\varsigma}{\kappa})q_{\varepsilon}(\tau)+\frac{\varsigma}{\kappa}t$. Using the convexity of $Q^{lp}(s)$, we get
\begin{eqnarray*}\label{eq1.67}
(1-\frac{\varsigma}{\kappa})Q^{lp}(q_{\varepsilon}(\tau))+\frac{\varsigma}{\kappa}Q^{lp}(t)\geq Q^{lp}(q_{\varepsilon}(\tau)+\varsigma e),
\end{eqnarray*}
and hence
\begin{eqnarray*}\label{eq1.68}
& & \frac{\varsigma}{\kappa}(Q^{lp}(t)-Q^{lp}(q_{\varepsilon}(\tau)))\geq Q^{lp}(q_{\varepsilon}(\tau)+\varsigma e)-Q^{lp}(q_{\varepsilon}(\tau))\\
&=&Q(q_{\varepsilon}(\tau)+\varsigma e)-Q(q_{\varepsilon}(\tau))
+r_{p}(q_{\varepsilon}(\tau)+\varsigma e)-r_{p}(q_{\varepsilon}(\tau))\\
&\geq&h(\varsigma)-2\nabla_{p}(\varsigma),
\end{eqnarray*}
where
\begin{eqnarray*}\label{eq1.69}
h(\varsigma)&=&\inf_{|t-q_{\varepsilon}(\tau)|=\varsigma}(Q(t)-Q(q_{\varepsilon}(\tau)))\\
\nabla_{p}(\varsigma)&=&\sup_{|t-q_{\varepsilon}(\tau)|\leq\varsigma}|Q^{lp}(t)-Q(t)|.
\end{eqnarray*}
According to \eqref{eq1.66}, there must be a $\epsilon>0$ such that if $0< p-1<\epsilon$, $h(\varsigma)-2\nabla_{p}(\varsigma)>0$. So $Q^{lp}(t)>Q^{lp}(q_{\varepsilon}(\tau))$ if $t\notin (q_{\varepsilon}(\tau), \varsigma)$ and thus $q^{lp}_{\varepsilon}(\tau)\in (q_{\varepsilon}(\tau), \varsigma)$. The arbitrariness of $\varsigma$ shows the desired result. \ $\Box$

{\bf Proof of Theorem 4.2}. We mainly need to prove
\begin{eqnarray}\label{eq1.72}
\frac{1}{T}\sum^{T}_{t=1}\boldsymbol{\psi}_{\tau, p}(y_{t}-\textbf{x}'_{t}\hat{\boldsymbol{\beta}}_{T, p}(\tau))
\stackrel{\textit{P}}{\longrightarrow}E\boldsymbol{\psi}_{\tau, p}(u-q^{lp}_{u}(\tau)).
\end{eqnarray}
Write
\begin{eqnarray*}\label{eq1.73}
& & \Big|\frac{1}{T}\sum^{T}_{t=1}\boldsymbol{\psi}_{\tau, p}(y_{t}-\textbf{x}'_{t}\hat{\boldsymbol{\beta}}_{T, p}(\tau))
-E\boldsymbol{\psi}_{\tau, p}(u-q^{lp}_{u}(\tau))\Big|\\
&\leq&\frac{1}{T}\Big|\sum^{T}_{t=1}\boldsymbol{\psi}_{\tau, p}(y_{t}-\textbf{x}'_{t}\hat{\boldsymbol{\beta}}_{T, p}(\tau))
-\sum^{T}_{t=1}E\boldsymbol{\psi}_{\tau, p}(y_{t}-\textbf{x}'_{t}\hat{\boldsymbol{\beta}}_{T, p}(\tau))\Big|\\
& & +\frac{1}{T}\sum^{T}_{t=1}|E\boldsymbol{\psi}_{\tau, p}(y_{t}-\textbf{x}'_{t}\hat{\boldsymbol{\beta}}_{T, p}(\tau))-E\boldsymbol{\psi}_{\tau, p}(u-q^{lp}_{u}(\tau))|\\
&\leq&\sup_{\boldsymbol{\delta}\in U[\boldsymbol{\beta}_{p}(\tau), \textbf{r}_{1}]}\frac{1}{T}\Big|\sum^{T}_{t=1}\boldsymbol{\psi}_{\tau, p}(y_{t}-\textbf{x}'_{t}\boldsymbol{\delta})-E\boldsymbol{\psi}_{\tau, p}(y_{t}-\textbf{x}'_{t}\boldsymbol{\delta})\Big|+\frac{1}{T}\sum^{T}_{t=1}o_{P}(1),
\end{eqnarray*}
where the last inequality is valid in probability according Assumption \ref{ass4.6} and the result $\hat{\boldsymbol{\beta}}_{T, p}(\tau)\stackrel{\textit{P}}{\longrightarrow}\boldsymbol{\beta}_{p}(\tau)$ which can be obtained by Theorem \ref{thm2.1} as the assumptions in Section 4 satisfies the requirement of Theorem \ref{thm2.1}. Using Assumption \ref{ass4.5}, we obtain \eqref{eq1.72}.
Then using Lemmas \ref{lem2} and \ref{lem4}, we get $E\boldsymbol{\psi}_{\tau, p}(u-q^{lp}_{u}(\tau))\rightarrow f(0)$. Based on this, using Assumption \ref{ass4.1} finally completes the proof. $\Box$
\section{Conclusion}
In this article we have proposed composite $L^p$-quantile regression and have established the relevant asymptotic theory.
We have further considered the oracle theory of penalized composite $L^p$-quantile regression.
In order to smooth the objective function of quantile regression, we have proposed near quantile regression.
The simulation and empirical analysis have both demonstrated the merits of our proposed methodology.
Of note, the provided algorithm could be effectively used to fit quantile regression in high-dimensional regime, which could
help improve the status of quantile regression in the machine learning field. As to why and when the algorithm works in modelling quantiles, we believe that a rigorous theoretical analysis is necessary. This is an open problem for future research.
In addition, based on the near quantile regression, there are many interesting problems to be further explored.

\textbf{Acknowledgements}\\

The author's research was supported by the
Opening Project of Sichuan Province University Key Laboratory of
Bridge Non-destruction Detecting and Engineering Computing
(2024QYY02). The author is very grateful for the help of Yu Chen, Xiao Guo and Jie Hu at University of Science and Technology of China.

{\hspace {-9.0mm}}

\end{document}